%% file: Bayesian_SIAC_main.tex
\documentclass[final,onefignum,onetabnum,letterpaper]{siamonline250211}
\usepackage[scale=0.8]{geometry} 

\input{Bayesian_SIAC_shared}

\usepackage{lineno}
\usepackage{comment}

\ifpdf
\hypersetup{
	pdftitle={Bayesian SIAC filter},
 	pdfauthor={Glaubitz}
}
\fi

\begin{document}

\maketitle

\begin{abstract}
	\input{0_abstract}
\end{abstract}



\begin{Code}
    \url{https://github.com/RomanStuhlmacher/paper-2025-Bayesian-SIAC-Filter} 
\end{Code}

\begin{DOI}
    \url{https://doi.org/10.1137/25M1800433}
\end{DOI}

\input{1_intro} 
\input{2_prelim}
\input{3_Bayesian} 
\input{4_MAP}
\input{5_MCMC}
\input{6_numerics} 
\input{7_summary}

\section*{Acknowledgements}
\input{acknowledgements}

\appendix 
\input{appA_notation}
\input{appB_SIACimpl}

\bibliographystyle{siamplain}
\bibliography{literature}

\end{document}

%% file: Bayesian_SIAC_shared.tex
\usepackage{soul}
\usepackage{lipsum}
\usepackage{amsfonts}
\usepackage{epstopdf}
\ifpdf
  \DeclareGraphicsExtensions{.eps,.pdf,.png,.jpg}
\else
  \DeclareGraphicsExtensions{.eps}
\fi

\usepackage{datetime}
\newdateformat{monthyeardate}{%
	\monthname[\THEMONTH] \THEDAY, \THEYEAR}

\usepackage{academicons}
\usepackage{xcolor}
\renewcommand{\orcid}[1]{\href{https://orcid.org/#1}{\textcolor[HTML]{A6CE39}{orcid.org/#1}}}

\usepackage{amsmath}
\allowdisplaybreaks
\usepackage{amssymb}
\usepackage{commath}
\usepackage{mathtools}
\usepackage{bbm}
\usepackage{bm}

\usepackage{color}
\usepackage{graphicx}
\usepackage[small]{caption}
\usepackage{subcaption}

\usepackage{relsize}
\usepackage{adjustbox}
\usepackage{algorithm}
\usepackage[noend]{algpseudocode}
\usepackage{booktabs}
\usepackage{tikz}
\usepackage{cancel}
\usepackage{pifont}
\usetikzlibrary{bayesnet} 

\usepackage{enumitem}
\setlist[enumerate]{leftmargin=.5in}
\setlist[itemize]{leftmargin=.5in}

\newsiamremark{remark}{Remark}
\newsiamremark{example}{Example}
\newsiamremark{hypothesis}{Hypothesis}
\crefname{hypothesis}{Hypothesis}{Hypotheses}
\newsiamthm{claim}{Claim}

\headers{The Bayesian SIAC filter}{J.\ Glaubitz, T.\ Li, J.\ Ryan, and R.\ Stuhlmacher}

\title{The Bayesian SIAC filter\thanks{\monthyeardate\today 
\corresponding{Jan Glaubitz} 
}}
\author{ 
Jan Glaubitz\thanks{
Department of Mathematics, Link\"oping University, Sweden
(\email{jan.glaubitz@liu.se}, \orcid{0000-0002-3434-5563})
} 
\and 
Tongtong Li\thanks{
Department of Mathematics and Statistics, University of Maryland, Baltimore County, MD, USA 
(\email{tongtong.li@umbc.edu}, \orcid{0000-0002-7664-4764})
}
\and 
Jennifer Ryan\thanks{
Department of Mathematics, KTH Royal Institute of Technology, Sweden
(\email{jryan@kth.se}, \orcid{0000-0002-6252-8199}; \email{romanst@kth.se}, \orcid{0009-0004-4842-0765})
} 
\and 
Roman Stuhlmacher\footnotemark[4]
}

\usepackage{amsopn}

\usepackage{comment}

\usepackage[normalem]{ulem}

\usepackage[normalem]{ulem}

\DeclareMathAlphabet\mathbfcal{OMS}{cmsy}{b}{n}

\DeclareMathOperator*{\argmin}{arg\,min}

\DeclareMathOperator*{\kernel}{kernel}

\newcommand{\intd}{\, \mathrm{d}} 

\newcommand{\N}{\mathbb{N}}

\newcommand{\R}{\mathbb{R}}

\newcommand{\cgam}[1]{c_{#1}}
\newcommand{\KH}[4]{K_{#1}^{(#2,#3)}\left(#4\right)}

\newcommand{\FullK}[3]{\sum_{\gamma = 0}^{#1}\, \cgam{\gamma} {#2}({#3}-x_\gamma)}
\newcommand{\BS}[1]{B^{(#1)}}


%% file: 0_abstract.tex
We propose the \emph{Bayesian {S}moothness-{I}ncreasing {A}ccuracy-{C}onserving (SIAC) filter}---a hierarchical Bayesian generalization of the existing deterministic SIAC filter. 
The SIAC filter is a powerful numerical tool for removing high-frequency noise from data or numerical solutions without degrading accuracy.
However, current SIAC methodology is limited to (i) nodal/modal data (noisy direct function values/coefficients of a piecewise polynomial function approximation) and (ii) deterministic point estimates that do not account for uncertainty propagation of input data to the SIAC reconstruction.
The proposed Bayesian SIAC filter overcomes these limitations by (i) supporting \emph{general (non-nodal) data models} and (ii) enabling \emph{rigorous uncertainty quantification (UQ)}, thereby broadening the applicability of SIAC filtering.
We also develop structure-exploiting algorithms for efficient maximum a posteriori (MAP) estimation and Markov chain Monte Carlo (MCMC) sampling, with a focus on linear data models with additive Gaussian noise.
Computational experiments demonstrate the effectiveness of the Bayesian SIAC filter across several applications, including signal denoising, image deblurring, and post-processing of numerical solutions to hyperbolic conservation laws. 
The results show that the Bayesian approach produces point estimates with accuracy comparable to, and in some cases exceeding, that of the deterministic SIAC filter. 
In addition, it extends naturally to general data models and provides built-in UQ.

%% file: 1_intro.tex
\section{Introduction} 
\label{sec:intro} 

The {S}moothness-{I}ncreasing {A}ccuracy-{C}onserving (SIAC) filter is a potent numerical tool designed to remove high-frequency noise from data or numerical solutions without compromising their accuracy. 
It does so by accelerating the decay of the Fourier coefficients, similar to more traditional filters that are designed to manage oscillations near discontinuities \cite{PICKLO2024JCP}.
It operates by convolving a quantity of interest $u$ (e.g., nodal data or a numerical approximation) with a carefully constructed kernel $K$, producing the filtered function
\begin{equation}\label{eq:SIAC_filter}
    \mathcal{F}[u](x) = \int K( x - y ) u(y) \intd y.
\end{equation}
The kernel $K$ is specifically chosen to enhance the smoothness of $u$ while preserving---or even improving---its convergence properties and usually consists of B-spline functions. 
Introduced initially as a post-processing technique for finite element methods \cite{bramble1977higher,cockburn2003enhanced}, the SIAC filter has been most prominently developed in the context of discontinuous Galerkin (DG) methods and hyperbolic conservation laws \cite{ryan2005extension,curtis2008postprocessing,mirzaee2010quantification,mirzaee2012efficient,ji2014superconvergent}.
Over the past two decades, the methodology has been extended in several directions, including adaptations for boundary treatment and non-uniform meshes \cite{van2011position,mirzaee2011smoothness,ryan2015one,li2019smoothness}, shock-capturing techniques for managing Gibbs oscillations near a discontinuity \cite{wissink2018shock,bohm2019multi}, entropy correction \cite{EdohPicklo}, multi-resolution analysis \cite{SIAC-MRA,picklo}, as well as Vlasov-Maxwell equations \cite{galindoolartevm}, and denoising particle-in cell methods \cite{PICKLO2024JCP}.
Also see \cite{docampo2020enhancing} and references therein for a review of the technique. 

Notably, the existing SIAC methodology primarily focuses on reconstructive accuracy and assumes that the input data are either nodal data or modal coefficients of a piecewise polynomial function approximation.\footnote{Modal coefficients can be transformed into nodal data by evaluating the function approximation at grid points.}
However, taking a broader perspective reveals two important and largely overlooked challenges:
\begin{enumerate}
    \item[(i)] 
    quantifying the uncertainty in SIAC-filtered data, and

    \item[(ii)] 
    extending SIAC filtering to indirect data.
\end{enumerate}
In practice, uncertainty arises from both input data and the data-generating process. 
This uncertainty then propagates through the reconstruction, affecting subsequent predictions and decision-making. 
Quantifying such uncertainty is often critical in real-world applications. 
Moreover, extending the SIAC framework to handle indirect data is highly desirable. 
In many real-world applications, the available data are both noisy and non-nodal. 
That is, the observations correspond \emph{not} to direct evaluations of the underlying function $u$, but to noisy, non-trivial transformations (linear or nonlinear) such as weighted averages or integral transforms. 
Such non-nodal data commonly arise from indirect measurement processes \cite{vogel2002computational,kaipio2006statistical,hansen2006deblurring,hansen2010discrete,hansen2021computed}. 
Prominent examples include blurred measurements in deconvolution problems, Fourier data in magnetic resonance imaging (MRI) and synthetic aperture radar (SAR), and Radon data in computed tomography (CT). 
We formalize these general data models in \Cref{sec:Bayesian}.

\subsection*{Our contribution: A Bayesian perspective}

We propose the \emph{Bayesian SIAC filter}, which places SIAC filtering within a probabilistic framework. 
This formulation enables two key advancements: (i) the ability to quantify uncertainty in SIAC-filtered data and (ii) the extension of SIAC filtering to general (non-nodal) data models.
To achieve this, we reinterpret SIAC filtering as an inverse problem, in which observational data are used to recover the true underlying function, assuming it possesses smoothness that is preserved by the SIAC filter. 
To quantify uncertainty as well as to automatically learn hyper-parameters (e.g., noise variance and regularization strength), we adopt a hierarchical Bayesian framework \cite{kaipio2006statistical,stuart2010inverse,calvetti2023bayesian}.
In this Bayesian setting, all relevant quantities are modeled as random variables equipped with prior distributions. 
The computational goal is then to infer the posterior distribution of the unknown function 
$u$ given the observational data. 
The posterior naturally combines a likelihood function, encoding assumptions about the data-generating process and noise characteristics; and
a prior distribution, capturing our beliefs about the structure of $u$.
The core idea behind the Bayesian SIAC filter is that the prior is directly informed by the foundational SIAC assumption: the true function 
$u$ is smooth, and its smoothness is preserved by SIAC filtering, i.e., $\mathcal{F}[u] \approx u$. 

We note that incorporating smoothing or filtering operators into Bayesian priors has appeared in the broader literature on Bayesian inverse problems \cite{bardsley2012mcmc,calvetti2023bayesian,LiGelb2025}.
In the present work, however, the prior is derived from the numerical properties of the SIAC filter, which is designed to preserve polynomial accuracy and consistency in high-order post-processing of discretizations. 
While alternative convolution kernels such as Gaussian or top-hat filters could in principle be used within a Bayesian formulation, doing so would lead to a fundamentally different model that no longer reflects the numerical motivations underlying SIAC filtering.

For efficient posterior inference, we focus on linear data models with Gaussian additive noise. 
This setting allows highly efficient maximum a posteriori (MAP) estimation using a block coordinate descent (BCD) approach \cite{wright2015coordinate,beck2017first}, inspired by structurally similar hierarchical Bayesian models in the context of generalized sparse Bayesian learning (GSBL) \cite{calvetti2020sparse,glaubitz2023generalized,xiao2023sequential,glaubitz2023leveraging,lindbloom2024generalized,glaubitz2025efficient,LiGelb2025}.
Additionally, by exploiting the conditional conjugacy between the proposed Gaussian prior and the gamma hyper-prior, we derive an efficient Gibbs sampler. 
This sampler, building on earlier work in imaging \cite{bardsley2012mcmc}, enables efficient Markov chain Monte Carlo (MCMC) sampling of the Bayesian SIAC posterior, thereby facilitating rigorous uncertainty quantification (UQ).
Notably, the approaches described above for efficient inference rely on well-established and readily available methods. 
We deliberately adopt mature optimization and sampling tools to ensure conceptual clarity, robustness, and ease of adoption. 
The primary contribution of this work is the formulation of a Bayesian extension of the SIAC filter that enables uncertainty quantification and accommodates general (non-nodal) data models.

We demonstrate the performance of the proposed Bayesian SIAC filter, encompassing MAP estimation, MCMC sampling, and UQ, on a variety of test problems, including signal denoising, post-processing of numerical DG solutions of the linear advection equation, and image deblurring.
Our numerical experiments show that the MAP and mean estimates of the Bayesian SIAC filter often yield point estimates with accuracy comparable to, and in some cases exceeding, that of the deterministic SIAC filter. 
Moreover, the Bayesian SIAC filter significantly extends the applicability of the traditional deterministic SIAC filter: it can be applied to general (non-nodal) data models and enables UQ.

At the same time, we emphasize the intended scope of the Bayesian SIAC filter. 
Our results indicate that it is most beneficial when (i) the observational data are strongly corrupted by noise and UQ is desired, or (ii) the data are indirect, rendering deterministic SIAC filters inapplicable.  
In contrast, for post-processing highly resolved DG solutions of smooth hyperbolic conservation laws, the classical deterministic SIAC filter remains preferable: 
In such cases, it exhibits superconvergence properties that the Bayesian SIAC filter does not replicate, as shown in our numerical experiments. 

In conclusion, we anticipate that the proposed Bayesian SIAC filter will pave the way for the SIAC methodology to be applied to a much broader class of problems in future work.

\subsection*{Outline}

We begin in \Cref{sec:prelim} by reviewing the necessary preliminaries on the existing deterministic SIAC filter.
In \Cref{sec:Bayesian}, we introduce the proposed Bayesian SIAC filter and formulate it within a hierarchical Bayesian framework.
Efficient MAP estimation and Gibbs sampling strategies for the Bayesian SIAC filter are then discussed in \Cref{sec:MAP} and \Cref{sec:MCMC}, respectively.
Next, in \Cref{sec:numerics}, we evaluate the performance of the Bayesian SIAC filter through a series of computational experiments.
Finally, \Cref{sec:summary} summarizes our findings and outlines directions for future work. 
Furthermore, we provide a table with a notation overview in \Cref{app:notation} and implementation details for the SIAC filter and its matrix representation in \Cref{sec:SIAC-imple}.

%% file: 2_prelim.tex
\section{Preliminaries on SIAC filtering} 
\label{sec:prelim} 

We begin by establishing the notation used in this manuscript. 
We denote scalars, scalar-valued functions, and operators between function spaces by non-bold letters (e.g., $x \in \R$, $H > 0$, $u: \R^d \to \R$, $K: \R^d \to \R$, and $\mathcal{F}: u \mapsto \mathcal{F}[u]$). 
Furthermore, we denote vectors and vector-valued functions by bold letters (e.g., $\mathbf{x}^N \in \R^N$ and $\mathbfcal{A}: \R^N \to \R^M$), including matrices $\mathbf{A} \in \R^{M \times N}$. 
Specifically, $u: \R^d \to \R$ denotes an underlying scalar-valued function on a $d$-dimensional domain, and $\mathbf{u} \in \R^N$ its finite-dimensional discretization. 
Furthermore, $\mathcal{F}$ denotes the SIAC filter (an operator between function spaces) and $\mathbf{F} \in \R^{N \times N}$ is its matrix representation. 

The SIAC filter was originally designed to post-process Galerkin finite element (FE) approximations, given by piecewise polynomial functions \cite{bramble1977higher,mock1978computation,cockburn2003enhanced,docampo2020enhancing}. It is based on the observation that piecewise polynomial approximations often converge with a higher order of accuracy in the negative-order norm than in the usual Sobolev norms. 
This suggests that one can construct a function that approximates the solution pointwise and converges with a higher order of accuracy than the FE approximation itself \cite{thomee1984galerkin}.

The SIAC filter \cref{eq:SIAC_filter} realizes this idea by convolving the given data $u$, such as a numerical approximation,  with a kernel $K$ based on a linear combination of function translates. 
Convolution enables the modification of the Fourier coefficients such that they decay faster---even if the data describes a phenomenon with high-frequency oscillations  \cite{PICKLO2024JCP}. 
The kernel $K$ is based on a linear combination of function translates. 
This convolution kernel can be specified by its basis functions $\Psi = \{ \psi(\cdot -x_{\gamma})\}_{\gamma=0}^{r}$ and shifts $x_{\gamma} \in \mathbb{R}^d$, yielding $K = K^{(r+1,\Psi)}$ with $\KH{}{r+1}{\Psi}{x} = \FullK{r}{\psi}{x}$. 
For structured discretizations, the $x_{\gamma}$ are spaced at a distance equal to the mesh size. 
There are many options for non-uniform meshes beyond using the maximum element size \cite{curtis2008postprocessing,ryan2015one}. 
The kernel coefficients $c_{\gamma}$ are determined by enforcing consistency and $r$ moment conditions \cite{mock1978computation,tanaka2019investigation}. The $x_{\gamma}$, together with the basis function $\psi$, can be chosen to take advantage of the underlying data and geometry. 
The kernel can achieve up to $r+1$ orders of accuracy, and up to $r+2$ if the kernel is symmetric. 
As $r$ can be chosen as desired, the limitation in extracting higher orders of accuracy that are contained in the negative-order norm therefore relies on the generation of the underlying data and any associated model problem. 
For example, for models described by transport equations and DG methods, the accuracy is $2k+1$ when piecewise polynomials of degree $k$ are used \cite{cockburn2003enhanced}. 

\begin{figure}[tb]
    \centering
    \includegraphics[width=.99\linewidth]{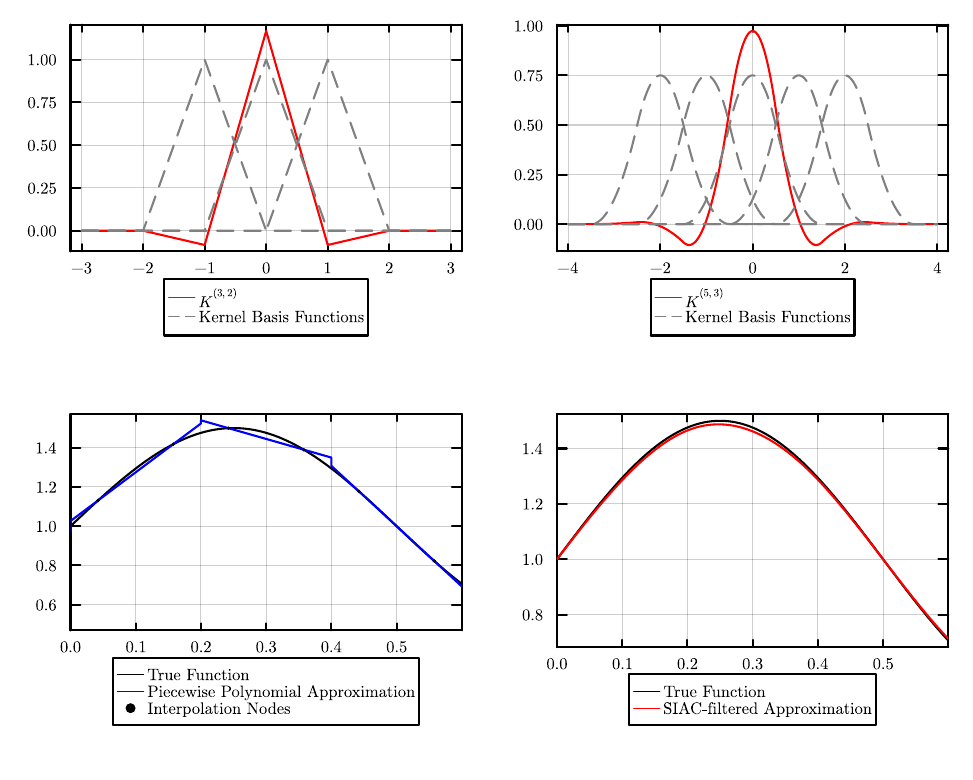}
    \caption{Top row: SIAC kernels (straight red line) and the shifted B-splines for $(r+1,\ell)=\{(3,2),\, (5,3)\}$ (dashed gray line). 
    Bottom row: Differences between an exact function and its piecewise polynomial approximation (left) and the SIAC-filtered approximation, using $K^{(5,3)}$ (right).
    }
    \label{fig:SIAC_kernel}
\end{figure}

For computational efficiency, B-splines are often used for the basis functions $\Psi$. 
This filter is the most common version of the SIAC filter. 
In one dimension, utilizing $r+1$ central, univariate B-Splines, $\BS{\ell},$ of order $\ell$, the kernel is $\KH{}{r+1}{\ell}{x} = \FullK{r}{\BS{\ell}}{x}$.
\cref{fig:SIAC_kernel} illustrates $K^{(3,2)}(x)$ and $K^{(5,3)}(x)$ with central B-splines. 
The kernel is then scaled by a factor $H$, i.e. $K^{(r+1,\ell)}_H(x) = \frac{1}{H}K^{(r+1,\ell)}(\frac{x}{H})$, that relies on the length scale of the discretization\footnote{For uniform quadrilateral meshes, the length scale is simply taken as the mesh size.}. 
The kernel thus becomes 
\begin{equation}\label{eq:Sym_SIAC_kernel}
    K^{(r+1,\ell)}_H(x) = \frac{1}{H} \sum_{\gamma=0}^r c_{\gamma} B^{(\ell)}\left( \frac{x - x_{\gamma}}{H} \right).
\end{equation}
Note that the support of the kernel is $\left[-\frac{r+\ell}{2}H,\frac{r+\ell}{2}H\right]$ and hence the symmetric filter should be confined to the interior of the domain unless utilizing periodic boundary conditions. 
Otherwise, a suitable modification can be made \cite{van2011position,ryan2015one}. 
\cref{fig:SIAC_kernel} also illustrates a piecewise polynomial approximation and its SIAC-filtered version.

The error estimate \cref{eq:Error} below essentially states when $\mathcal{F}[u] \approx u$ is satisfied and thus quantifies the assumption \cref{eq:prior_assumption} that motivates the formulation of a smoothness-promoting Bayesian prior based on the SIAC filter in \Cref{sub:Bayesian_prior}. 
Assuming sufficient local regularity of the approximated function $u$ (the ground truth), it is possible to obtain an error estimate using the SIAC-filtered function $u$: 
\begin{equation}\label{eq:Error}
    \| u - \mathcal{F}[u] \|_{L^2} \leq C_F H^{r+1}
\end{equation} 
with SIAC filter $\mathcal{F}$ as in \cref{eq:SIAC_filter}, kernel $K = K^{(r+1,\ell)}_H$ as in \cref{eq:Sym_SIAC_kernel}, and a constant $C_F$ that is independent of $H$.  
Notably, \cref{eq:Error} assumes $u \in \mathcal{H}^{s+1}(\Omega)$ with $0 < s \leq \ell$, where $\mathcal{H}^s(\Omega) = \{v \in L^1_{loc}(\Omega) \mid \partial^{\alpha}v \in L^2(\Omega) \ \forall \alpha \in \mathbb{N}^d \text{ with } |\alpha|\leq s \}$ is the usual $s$-th order Sobolev space.  
Importantly, the parameter $\ell$ controls the smoothness of $\mathcal{F}[u]$. 
For example, when utilizing B-splines of order $\ell,$ the global smoothness of $\mathcal{F}[u]$ is $\ell-2$. 
See \cite{bramble1977higher,cockburn2003enhanced} for more details. 

\begin{remark} (On the kernel scaling).
There are many studies investigating the effects of the scaling parameter $H$ and the associated kernel support size.
For instance, \cite{king2012smoothness,li2019smoothness} studied the effect of varying $H$ on the accuracy and convergence rates, with emphasis on nonuniform meshes. 
There are also several approaches to selecting a locally varying scaling parameter; see \cite{jallepalli2019adaptive,jallepalli2024non} for recent examples. 
In the present work, we restrict our attention to the classical choice, in which the scaling is set equal to the characteristic mesh size. 
This serves as a benchmark case for comparing the Bayesian approach with the deterministic SIAC filtering methodology.
\end{remark}

%% file: 3_Bayesian.tex
\section{The Bayesian SIAC filter} 
\label{sec:Bayesian} 

We now describe the proposed \emph{Bayesian SIAC filter} framework. 
To this end, consider a data-generating process of the form 
\begin{equation}\label{eq:data_model}
	\mathbf{b} = \mathbfcal{A}(\mathbf{u}) + \mathbf{e}, 
\end{equation} 
where $\mathbf{b} \in \R^M$ is the given observational data, $\mathbfcal{A}: \R^N \to \R^M$ is a known forward operator, $\mathbf{u} \in \R^N$ is the unknown parameter vector of interest, and $\mathbf{e} \in \R^M$ denotes an unknown additive noise component.

\subsection{A Bayesian perspective}
\label{sub:Bayesian_perspective}

If $\mathbfcal{A}$ is ill-conditioned or if the observational data $\mathbf{b}$ is scarce ($M < N$) or distorted by noise, the inverse problem \cref{eq:data_model} becomes ill-posed and pathologically hard to solve.
Prior knowledge about $\mathbf{u}$ is often leveraged to overcome the associated challenges. 
In this regard, using a Bayesian approach \cite{stuart2010inverse,calvetti2023bayesian}, which models the involved quantities as random variables, is known to be highly successful. 
Specifically, we model the available observation data $\mathbf{b}$ and the parameter vector of interest $\mathbf{u}$ as random variables. 
The goal of computation is to explore the \emph{posterior distribution}, whose density is given by Bayes' theorem as 
\begin{equation}\label{eq:Bayes} 
	\pi^{\mathbf{b}}( \mathbf{u}, \boldsymbol{\theta} ) 
		\propto f( \mathbf{u}, \boldsymbol{\theta}; \mathbf{b} ) \, 
		\pi^0( \mathbf{u}, \boldsymbol{\theta} ),
\end{equation} 
where ``$\propto$'' means that the two sides are equal up to a normalizing multiplicative constant that is independent of $\mathbf{u}$ and $\boldsymbol{\theta}$. 
The \emph{likelihood function} $f( \mathbf{u}, \boldsymbol{\theta}; \mathbf{b} )$ encodes our knowledge about the data-generating process and the noise characteristics. 
The \emph{prior} $\pi^0( \mathbf{u}, \boldsymbol{\theta} )$ models our a priori belief about the unknown parameters, usually formulated before any data is observed. 
Finally, $\boldsymbol{\theta}$ collects all hyper-parameters. 
For instance, we will later model the noise and prior precision, $\alpha, \beta>0$, as random variables so that $\boldsymbol{\theta} = [\alpha,\beta]$.

\subsection{The likelihood}
\label{sub:Bayesian_likelihood}

Consider the data model \cref{eq:data_model} and assume that $\mathbf{e}$ is a realization of independent and identically distributed (i.i.d.) zero-mean normal noise, $\mathbf{e} \sim \mathcal{N}(\mathbf{0}, \alpha^{-1} \mathbf{I})$, which is a typical modeling choice. 
Here, $\mathcal{N}(\boldsymbol{\mu}, \boldsymbol{\Sigma})$ denotes the multivariate normal distribution with mean $\boldsymbol{\mu}$ and covariance matrix $\boldsymbol{\Sigma}$.
At the same time, we also treat the noise precision (inverse variance) $\alpha>0$ as a random variable, as it is often unknown a priori. 
In this case, the likelihood function becomes  
\begin{equation}\label{eq:likelihood}
	f( \mathbf{u}, \boldsymbol{\theta}; \mathbf{b} ) 
		\propto \alpha^{M/2} \exp\left( -\frac{\alpha}{2} \| \mathbfcal{A}(\mathbf{u}) - \mathbf{b} \|_2^2 \right) \pi^0( \alpha ),
\end{equation}
where $\pi^0( \alpha )$ is the hyper-prior density for $\alpha$. 
We assume $\alpha \sim \mathcal{G}(c_{\alpha},d_{\alpha})$, where $\mathcal{G}(c,d)$ denotes the univariate gamma distribution with shape and rate parameters $c,d > 0$.
This yields the density $\pi^0( \alpha ) \propto \alpha^{c_{\alpha}-1} \exp\left( - d_{\alpha} \alpha \right)$, where $c_{\alpha},d_{\alpha} > 0$ and \cref{eq:likelihood} becomes 
\begin{equation}\label{eq:likelihood2}
	f( \mathbf{u}, \boldsymbol{\theta}; \mathbf{b} ) 
		\propto \alpha^{M/2 + c_{\alpha} - 1} \exp\left( -\frac{\alpha}{2} \| \mathbfcal{A}(\mathbf{u}) - \mathbf{b} \|_2^2 - d_{\alpha} \alpha \right).
\end{equation}
Notably, the gamma hyper-prior is a conjugate for the conditional likelihood in the following sense: 
If $\mathbf{b} | \mathbf{u}, \alpha \sim \mathcal{N}(\mathbfcal{A}(\mathbf{u}), \alpha^{-1} \mathbf{I})$ and $\alpha \sim \mathcal{G}(c, d)$, then $\alpha | \mathbf{u}, \mathbf{b} \sim \mathcal{G}(\tilde{c}_{\alpha}, \tilde{d}_{\alpha})$ for constants $\tilde{c}_{\alpha}$ and $\tilde{d}_{\alpha}$ that remain to be determined. 
In other words, the distribution of $\alpha$ conditioned on all other parameters is gamma-distributed.  
This will allow us to derive efficient approximate MAP estimation and sampling methods.  

At the same time, the gamma hyper-prior can be made approximately uninformative:
Recall that a gamma-distributed random variable, $\alpha \sim \mathcal{G}(c_{\alpha},d_{\alpha})$, has mean $\mathbb{E}[\alpha] = c_{\alpha}/d_{\alpha}$ and variance $\mathbb{V}[\alpha] = c_{\alpha}/d_{\alpha}^2$. 
In particular, $c_{\alpha} \to 1$ and $d_{\alpha} \to 0$ implies $\mathbb{E}[\alpha], \mathbb{V}[\alpha] \to \infty$, making the gamma hyper-prior uninformative. 
That is, $\alpha$ is mostly free from the moderating influence of the hyper-prior and allowed to ``freely" follow the observational data. 
In our later numerical tests, we use $c_{\alpha}=1$ and $d_{\alpha}= 10^{-3}$, which is similar to the choices in related works such as  \cite{tipping2001sparse,bardsley2012mcmc,glaubitz2023generalized,glaubitz2023leveraging}.

\begin{remark} 
Modeling the noise precision $\alpha$ as a random variable is not only a practical choice when $\alpha$ is unknown, but also a statistically beneficial one. 
Even when the exact noise precision is known, using it as a fixed value for $\alpha$ can yield suboptimal reconstructions \cite{zhang2011clarify}. 
Therefore, treating $\alpha$ hierarchically enables adaptive regularization that more accurately reflects the information in the data.
\end{remark}

\subsection{The SIAC prior}
\label{sub:Bayesian_prior}

The prior $\pi^0(\mathbf{u}, \boldsymbol{\theta})$ in \cref{eq:Bayes} models our a priori belief about the unknown parameter vector $\mathbf{u}$. 
Here, we assume that $\mathbf{u}$ represents the nodal values of a smooth function $u$ at certain grid points. 
As discussed in \Cref{sec:prelim}, specifically \cref{eq:Error}, we expect 
\begin{equation}\label{eq:prior_assumption}
    \mathcal{F}[u] \approx u
\end{equation}
if $u$ is sufficently smooth and $\mathcal{F}$ is an appropriate SIAC filter.  
In other words, applying the SIAC filter to a smooth function yields an output that more closely matches the original function. 
Conversely, if $u$ is contaminated by high-frequency artifacts, such as noise or numerical oscillations, the filter largely eliminates them. 
Consequently, enforcing \cref{eq:prior_assumption} promotes a smooth function devoid of spurious high-frequency components.

Let $\mathbf{F} \in \R^{N \times N}$ be the matrix representation of the SIAC filter $\mathcal{F}$ that uses $\mathbf{u}$, containing the nodal values of $u$ at some grid points, to approximate $\mathcal{F}[u]$ at the same grid points; See \Cref{sec:SIAC-imple} for details on $\mathbf{F}$.
We then model \cref{eq:prior_assumption} as 
\begin{equation}
	(\mathbf{F} - \mathbf{I}) \mathbf{u} \sim \mathcal{N}(\mathbf{0},\beta^{-1} \mathbf{I}),
\end{equation} 
where $\mathbf{I} \in \R^{N \times N}$ is the identity matrix and $\beta > 0$ is the so-called \emph{prior precision}. 
Note that the larger $\beta$, the more aggressively we promote smoothness via \cref{eq:prior_assumption}. 
We also treat the prior precision $\beta$ as a gamma distributed random variable $\beta \sim \mathcal{G}(c_{\beta}, d_{\beta})$ with fixed shape and rate parameters $c_{\beta}, d_{\beta} > 0$.
Consequently, we get the following \emph{SIAC prior}:
\begin{equation}\label{eq:prior}
	\pi^0(\mathbf{u}, \boldsymbol{\theta}) 
		\propto \beta^{N/2 + c_{\beta} - 1} \exp\left( -\frac{\beta}{2} \| (\mathbf{F} - \mathbf{I}) \mathbf{u} \|_2^2 - d_{\beta} \beta \right)
\end{equation} 
Finally, we remark that the rationale for choosing a gamma distribution for $\beta$ is the same as for $\alpha$, discussed in \Cref{sub:Bayesian_likelihood}.

\subsection{The posterior}
\label{sub:Bayesian_posterior}

We now return to the posterior distribution, whose density is provided by Bayes' theorem as \cref{eq:Bayes}. 
Substituting the likelihood function \cref{eq:likelihood2} and the SIAC prior density \cref{eq:prior} into \cref{eq:Bayes}, we get  
\begin{equation}\label{eq:posterior} 
\resizebox{.9\textwidth}{!}{$\displaystyle 
	\pi^{\mathbf{b}}( \mathbf{u}, \boldsymbol{\theta} ) 
		\propto \alpha^{M/2 + c_{\alpha} - 1} 
			\beta^{N/2 + c_{\beta} - 1}
			\exp\left( 
				-\frac{\alpha}{2} \| \mathbfcal{A}(\mathbf{u}) - \mathbf{b} \|_2^2 
				-\frac{\beta}{2} \| (\mathbf{F} - \mathbf{I}) \mathbf{u} \|_2^2
				- d_{\alpha} \alpha 
				- d_{\beta} \beta 
			\right). 
$}
\end{equation}   
Notably, \cref{eq:posterior} is high-dimensional and generally non-log-concave, rendering efficient inference a non-trivial task. 
We subsequently address MAP estimation and MCMC sampling in \Cref{sec:MAP,sec:MCMC}. 

We summarize the proposed Bayesian SIAC approach below in \cref{alg:BSIAC}.

\begin{algorithm}[htb]
\caption{Bayesian SIAC approach for recovering the nodal values $\mathbf{u} = [u_1,\dots,u_N]^T \in \R^N$ with $u_n = u(x_n)$, $n=1,\dots,N$, of a smooth function $u$}\label{alg:BSIAC}
\begin{algorithmic}
	\State{\textbf{Input:}} Observational data $\mathbf{b} \in \R^M$ and forward operator $\mathbfcal{A}$ in the assumed data model \cref{eq:data_model}, and grid points $\mathbf{x} = [x_1,\dots,x_N]^T \in \R^N$ at which $u$ should be recovered
    \State\textit{Step 1.} Formulate the likelihood function as in \cref{eq:likelihood2} using the observational data $\mathbf{b} \in \R^M$ and the forward operator $\mathbfcal{A}$
    \State\textit{Step 2.} Construct the SIAC matrix $\mathbf{F} \in \R^{N \times N}$ as described in \Cref{sec:SIAC-imple} 
    \State\textit{Step 3.} Use $\mathbf{F}$ to formulate the prior density as in \cref{eq:prior}
    \State\textit{Step 4.} Formulate the Bayesian SIAC posterior \cref{eq:posterior} by combining the likelihood and prior
    \State\textit{Step 5.} Fix the shape and rate parameters of the gamma hyper-priors for $\alpha$ and $\beta$, e.g., as $c_{\alpha} = c_{\beta} = 1$ and $d_{\alpha} = d_{\beta} = 10^{-3}$
    \State\textit{Step 6.} Leverage the resulting Bayesian SIAC posterior to infer the nodal values $\mathbf{u} = [u_1,\dots,u_N]^T \in \R^N$ using MAP estimation described in \Cref{sec:MAP} (see \cref{alg:MAP}) and/or MCMC sampling described in \Cref{sec:MCMC} (see \cref{alg:Gibbs})
\end{algorithmic}
\end{algorithm}

A few remarks are in order. 

\begin{remark}[General data models]
    In many existing applications of the SIAC filter, observational data are assumed to be either noisy direct measurements of the underlying signal or erroneous modal coefficients of a Galerkin approximation. 
    In the former case, the forward operator is $\mathbfcal{A}(\mathbf{u})=\mathbf{u}$, whereas in the latter it is a linear projection. 
    However, the proposed Bayesian SIAC filter framework can, in principle, be applied to arbitrary forward operators and measurement modalities. 
\end{remark}

\begin{remark}[Conflicting observational data]
    Since the Bayesian SIAC filter can be formulated for general data models, it naturally accommodates conflicting observational data---for example, multiple noisy measurements at the same $x$-location. 
    Suppose two such observations are given by $b_1 = u(x) + e_1$ and $b_2 = u(x) + e_2$. 
    Then the corresponding contribution to the likelihood \cref{eq:likelihood} (and hence to the posterior \cref{eq:posterior}) is proportional to $\exp(-\alpha[ (u(x)-b_1)^2 + (u(x)-b_2)^2 ]/2 )$, so that both measurements are incorporated simultaneously into the inference process. The prior then regularizes this information, steering the posterior toward a reconstructed value of $u(x)$ that balances agreement with the conflicting data and consistency with the assumed smoothness structure of the underlying function.
    Such scenarios cannot be accommodated within existing deterministic SIAC filtering frameworks, which typically assume a single prescribed function value at each spatial location.
\end{remark}

\begin{remark}[Relation to Gaussian regression methods]
Gaussian regression techniques, such as Kriging and Gaussian process models, are highly effective for surrogate modeling from sparse data and often outperform polynomial-based approaches in that setting, as shown, for example, in \cite{sen2015evaluation}. 
Also see \cite{gaul2014modified,gaul2015modified}, and references therein, for regression based on spline kernels with confidence intervals using Bayes’ theorem.
The objective of the present work, however, is different. 
The (Bayesian) SIAC framework is designed as a structure-exploiting reconstruction and post-processing method for discretized data, with a focus on preserving accuracy, enhancing smoothness, and quantifying uncertainty. 
From this perspective, SIAC-based methods are complementary to Gaussian regression approaches rather than direct competitors, as they leverage analytically derived kernels and known discretization structure rather than learning a global covariance model from data.
\end{remark}

%% file: 4_MAP.tex
\section{Efficient maximum a posteriori (MAP) estimation} 
\label{sec:MAP} 

We address MAP estimation for the proposed Bayesian SIAC framework. 
While we previously formulated the Bayesian SIAC model for a general forward operator $\mathbfcal{A}: \R^N \to \R^M$, we henceforth restrict the discussion to linear operators $\mathbfcal{A}(\mathbf{u}) = \mathbf{A} \mathbf{u}$ with $\mathbf{A} \in \R^{M \times N}$. 
This restriction enables efficient approximate MAP estimation by recasting it as an iteratively reweighted least-squares problem. 
Notably, the subsequent approach for efficient MAP estimation relies on well-established methods. 
We deliberately adopt existing inference tools to ensure conceptual clarity, robustness, and ease of adoption of the proposed Bayesian SIAC filter.

\subsection{The MAP estimate} 
\label{sub:MAP} 

Recall that $( \mathbf{u}^{\text{MAP}}, \alpha^{\text{MAP}}, \beta^{\text{MAP}} )$ is a MAP estimate of the posterior $\pi^{\mathbf{b}}( \mathbf{u}, \alpha, \beta )$ in \cref{eq:posterior} if it maximizes $\pi^{\mathbf{b}}( \mathbf{u}, \alpha, \beta )$. 
Equivalently---yet often computationally more robust---a MAP estimate $( \mathbf{u}^{\text{MAP}}, \alpha^{\text{MAP}}, \beta^{\text{MAP}} )$ minimizes the negative log-posterior, i.e., 
\begin{equation}\label{eq:MAP} 
	( \mathbf{u}^{\text{MAP}}, \alpha^{\text{MAP}}, \beta^{\text{MAP}} ) 
		= \argmin_{ ( \mathbf{u}, \alpha, \beta ) } \left\{ \mathcal{J}( \mathbf{u}, \alpha, \beta ) \right\},
\end{equation} 
where the objective function (also called \emph{Gibbs energy functional}) is 
\begin{equation}\label{eq:J} 
\begin{aligned} 
	\mathcal{J}( \mathbf{u}, \alpha, \beta ) 
		= & - \log \pi^{\mathbf{b}}( \mathbf{u}, \alpha, \beta ) \\ 
		= & \ \frac{\alpha}{2} \| \mathbf{A} \mathbf{u} - \mathbf{b} \|_2^2 
				+ \frac{\beta}{2} \| (\mathbf{F} - \mathbf{I}) \mathbf{u} \|_2^2
				+ d_{\alpha} \alpha 
				+ d_{\beta} \beta \\ 
		& \ - ( M/2 + c_{\alpha} - 1 ) \log \alpha
				- (N/2 + c_{\beta} - 1) \log \beta.
\end{aligned}
\end{equation} 
Here, the objective function $\mathcal{J}(\mathbf{u}, \alpha, \beta)$ is defined up to additive constants that are independent of $\mathbf{u}$, $\alpha$, and $\beta$. Since these constants do not influence the minimizer, they can be omitted without affecting the MAP estimate.

\subsection{The block-coordinate descent (BCD) approach}
\label{sub:IAS}

We use a block-coordinate descent (BCD) approach \cite{wright2015coordinate,beck2017first} to approximate the MAP estimate. 
The BCD approach computes the minimizer of the objective function $\mathcal{J}$ by alternatingly minimizing $\mathcal{J}$ w.r.t.\ (i) $\mathbf{u}$ for fixed $\alpha, \beta$, (ii) $\alpha$ for fixed $\mathbf{u}, \beta$, and (iii) $\beta$ for fixed $\mathbf{u}, \alpha$. 
Given an initial guess for the parameter vector $\mathbf{u}$, the BCD algorithm proceeds through a sequence of updates of the form 
\begin{equation}\label{eq:BCD} 
\begin{aligned}
	\mathbf{u} 
		& = \argmin_{\mathbf{u} \in \R^N} \left\{ \mathcal{J}(\mathbf{u}, \alpha, \beta) \right\}, \\ 
	\alpha 
		& = \argmin_{\alpha > 0} \left\{ \mathcal{J}(\mathbf{u}, \alpha, \beta) \right\}, \\ 
	\beta 
		& = \argmin_{\beta > 0} \left\{ \mathcal{J}(\mathbf{u}, \alpha, \beta) \right\}, 
\end{aligned}
\end{equation}
until a convergence criterion is met. 
The BCD algorithm is an attractive choice because the three subproblems \cref{eq:BCD} are significantly easier to solve than the joint minimization problem \cref{eq:MAP}. 

\begin{remark}
	Our approach is motivated by similar approaches in the context of GSBL \cite{calvetti2020sparse,glaubitz2023generalized,xiao2023sequential,glaubitz2023leveraging,lindbloom2024generalized,glaubitz2025efficient,LiGelb2025}, where conditionally Gaussian priors are combined with generalized gamma-distributed hyper-parameter vectors to promote sparsity. 
	In this context, a prevalent BCD method is the so-called iterative alternating sequential (IAS) algorithm \cite{calvetti2019hierachical,calvetti2020sparse,calvetti2020sparsity,lindbloom2024generalized,lindbloom2025priorconditioned}.
\end{remark}

\subsection{Updating the hyper-parameters $\alpha$ and $\beta$} 
\label{sub:theta_update}

We first address updating $\alpha$ given $\mathbf{u}$.  
Substituting the objective function \cref{eq:J} into the second minimization problem in \cref{eq:BCD}, and ignoring all terms that do not depend on $\alpha$, yields  
\begin{equation}\label{eq:alpha_update1} 
    \alpha = \argmin_{\alpha > 0} \left\{ 
	\frac{\alpha}{2} \| \mathbf{A} \mathbf{u} - \mathbf{b} \|_2^2 
	+ d_{\alpha} \alpha 
	- \left( M/2 + c_{\alpha} - 1 \right) \log \alpha 
    \right\}. 
\end{equation} 
We choose $M/2 + c_{\alpha} - 1 > 0$. 
Then, differentiating the objective function in \cref{eq:alpha_update1} w.r.t.\ $\alpha$, setting this derivative to zero, and solving for $\alpha$ yields 
\begin{equation}\label{eq:alpha_update}
    \alpha = \frac{M/2 + c_{\alpha} - 1}{ \| \mathbf{A} \mathbf{u} - \mathbf{b} \|_2^2/2 + d_{\alpha} }
\end{equation} 
for the $\alpha$-update. 
The same arguments as above result in 
\begin{equation}\label{eq:beta_update}
    \beta = \frac{N/2 + c_{\beta} - 1}{ \| (\mathbf{F} - \mathbf{I}) \mathbf{u} \|_2^2/2 + d_{\beta} }
\end{equation} 
for the $\beta$-update.

\subsection{Updating the parameter vector $\mathbf{u}$} 
\label{sub:x_update}

To update $\mathbf{u}$ given $\alpha,\beta$, we must solve 
the first minimization problem in \cref{eq:BCD}. 
Substituting the objective function \cref{eq:J} into it and ignoring all terms that do not depend on $\mathbf{u}$, the $\mathbf{u}$-update reduces to solving the quadratic minimization problem 
\begin{equation}\label{eq:u_update_quad} 
	\mathbf{u} 
		= \argmin_{ \mathbf{u} \in \R^N } \left\{ 
			\alpha \| \mathbf{A} \mathbf{u} - \mathbf{b} \|_2^2 
			+ \beta \| (\mathbf{F} - \mathbf{I}) \mathbf{u} \|_2^2
	\right\}, 
\end{equation} 
where we drop the constant factor 1/2.
Solving \cref{eq:u_update_quad} is equivalent to computing the least squares solution of the overdetermined linear system 
\begin{equation}\label{eq:u_update_LS}
    \begin{bmatrix} 
    		\alpha^{1/2} \mathbf{A} \\ 
		\beta^{1/2} (\mathbf{F} - \mathbf{I})   
	\end{bmatrix} 
	\mathbf{u} 
     = 
     \begin{bmatrix} 
     	\alpha^{1/2} \mathbf{b} \\ 
		\mathbf{0} 
	\end{bmatrix},
\end{equation} 
which, in turn, is equivalent to solving the regular linear system
\begin{equation}\label{eq:u_update_reg}
    \underbrace{\left( \alpha \mathbf{A}^T \mathbf{A} 
    		+ \beta (\mathbf{F}-\mathbf{I})^T (\mathbf{F}-\mathbf{I}) \right)}_{= \mathbf{C}} 
	\mathbf{u} 
    		= \alpha \mathbf{A}^T \mathbf{b}. 
\end{equation}
Notably, if the coefficient matrix $\mathbf{C}$ on the left-hand side of \cref{eq:u_update_reg} is symmetric positive-definite (SPD), each of \cref{eq:u_update_quad,eq:u_update_LS,eq:u_update_reg} shares the same unique solution. 

We summarize the BCD approach described above for efficient MAP estimation in \cref{alg:MAP}.

\begin{algorithm}[htb]
\caption{MAP estimation for the Bayesian SIAC posterior}\label{alg:MAP}
\begin{algorithmic}
	\State{\textbf{Input:}} Initial values $\alpha_0$ and $\beta_0$, and number of iterations $J_{\text{iter}}$
	\For{$j=1,\dots,J_{\text{iter}}$}
		\State\textit{Step 1.} Compute $\mathbf{u}_{j}$ by solving \cref{eq:u_update_reg} for $\alpha = \alpha_{j-1}$ and $\beta = \beta_{j-1}$
		 \State\textit{Step 2.} Compute $\alpha_{j}$ via \cref{eq:alpha_update} for $\mathbf{u} = \mathbf{u}_{j}$
		 \State\textit{Step 3.} Compute $\beta_{j}$ via \cref{eq:beta_update} for $\mathbf{u} = \mathbf{u}_{j}$
	\EndFor
	\State{\textbf{Output:}} Approximate MAP estimate $(\mathbf{u}_{J_{\text{iter}}}, \alpha_{J_{\text{iter}}}, \beta_{J_{\text{iter}}} ) \approx ( \mathbf{u}^{\text{MAP}}, \alpha^{\text{MAP}}, \beta^{\text{MAP}} )$
\end{algorithmic}
\end{algorithm} 

Finally, a few remarks are in order.

\begin{remark}[The common kernel condition]
Observe that $\mathbf{C}$ in \cref{eq:u_update_reg} is always symmetric.   
Hence, $\mathbf{C}$ is SPD if and only if it is invertible. 
This, in turn, is equivalent to $\kernel(\mathbf{C}) = 0$, which can be reformulated as the \emph{common kernel condition}
\begin{equation}\label{eq:common_kernel_C}
    \kernel(\mathbf{A}) \cap \kernel(\mathbf{F}-\mathbf{I}) = \{ \mathbf{0} \}.
\end{equation} 
Here, $\kernel(\mathbf{G}) = \{ \, \mathbf{x} \in \R^N \mid \mathbf{G} \mathbf{x} = \mathbf{0} \, \}$ is the kernel of $\mathbf{G}$.
The common kernel condition \cref{eq:common_kernel_C} is an often-made assumption in Tikhonov regularization \cite{kaipio2006statistical,tikhonov2013numerical} and  GSBL \cite{glaubitz2023generalized,glaubitz2023leveraging,lindbloom2024generalized}. 
It ensures that the combination of prior information and the given observational data yields a well-posed problem. 
\end{remark} 

\begin{remark}[Iterative methods and prior-conditioning]\label{rem:iter_methods}
Various methods can be used to solve \cref{eq:u_update_quad,eq:u_update_LS,eq:u_update_reg}. 
For sufficiently small $N$, employing direct methods at a cost of $\mathcal{O}(N^3)$ flops is reasonable. 
If $N$ is large, iterative methods such as the conjugate gradient (CG) method \cite{saad2003iterative}, the conjugate gradient for least squares (CGLS) algorithm \cite{calvetti2018bayes,calvetti2020sparse}, the gradient descent approach \cite{glaubitz2023generalized}, or hybrid projection methods \cite{chung2024computational}, should be used---potentially in combination with priorconditioning \cite{lindbloom2024generalized,lindbloom2025priorconditioned}. 
For a potentially dense forward operator $\mathbf{A}$, the above iterative methods have a typical computational cost of $\mathcal{O}(N^2)$ per iteration and memory costs of $\mathcal{O}(N^2)$.
We demonstrate the computational efficiency of the BCD approach in our later numerical experiments (see \cref{tab:signal_MAP_results}).
\end{remark} 

%% file: 5_MCMC.tex
\section{Efficient MCMC sampling} 
\label{sec:MCMC} 

We now discuss MCMC sampling, enabling Monte Carlo approximations of posterior expectations. 
Specifically, the posterior expectation of a quantity $\mathcal{H}(\mathbf{u}, \alpha, \beta)$ is approximated as 
\begin{equation}\label{eq:MC_approximation} 
	\int \mathcal{H}(\mathbf{u}, \alpha, \beta) \pi^\mathbf{b}(\mathbf{u}, \alpha, \beta) \intd (\mathbf{u}, \alpha, \beta) 
		\approx \frac{1}{J_{\text{samp}}} \sum_{j=1}^{J_{\text{samp}}} \mathcal{H}( \mathbf{u}_j, \alpha_j, \beta_j ),
\end{equation} 
with $\pi^\mathbf{b}(\mathbf{u}, \alpha, \beta)$ denoting the posterior density and each tuple $(\mathbf{u}_j, \alpha_j, \beta_j)$ being a posterior sample. 
Notably, the subsequent approach for efficient MCMC sampling relies on well-established methods. 
We deliberately adopt existing inference tools to ensure conceptual clarity, robustness, and ease of adoption of the proposed Bayesian SIAC filter. 
Specifically, we follow \cite{bardsley2012mcmc} to derive an efficient Gibbs sampler for the proposed Bayesian SIAC filter.  
To this end, we again restrict the discussion to linear operators $\mathbfcal{A}(\mathbf{u}) = \mathbf{A} \mathbf{u}$ with $\mathbf{A} \in \R^{M \times N}$.

We note from the joint posterior $\pi^{\mathbf{b}}( \mathbf{u}, \alpha, \beta )$ in \cref{eq:posterior} that the full conditional densities are 
\begin{equation}\label{eq:cond_densities}
\begin{aligned}
	\pi^{\mathbf{b}}( \mathbf{u} | \alpha, \beta ) 
		& \propto \exp\left( 
				-\frac{\alpha}{2} \| \mathbf{A} \mathbf{u} - \mathbf{b} \|_2^2 
				-\frac{\beta}{2} \| (\mathbf{F} - \mathbf{I}) \mathbf{u} \|_2^2
			\right), \\ 
	\pi^{\mathbf{b}}( \alpha | \mathbf{u}, \beta ) 
		& \propto \alpha^{M/2 + c_{\alpha} - 1} 
			\exp\left( 
				-\frac{\alpha}{2} \| \mathbf{A} \mathbf{u} - \mathbf{b} \|_2^2 
				- d_{\alpha} \alpha 
			\right), \\
	\pi^{\mathbf{b}}( \beta | \mathbf{u}, \alpha ) 
		& \propto \beta^{N/2 + c_{\beta} - 1}
			\exp\left( 
				-\frac{\beta}{2} \| (\mathbf{F} - \mathbf{I}) \mathbf{u} \|_2^2 
				- d_{\beta} \beta 
			\right). 
\end{aligned}
\end{equation}
The above densities imply the following full conditional distributions: 
\begin{equation}\label{eq:cond_distributions}
\begin{aligned}
	\mathbf{u} | \alpha, \beta, \mathbf{b} 
		& \sim \mathcal{N}( \mathbf{C}^{-1} \alpha \mathbf{A}^T \mathbf{b}, \mathbf{C}^{-1} ), \\ 
	\alpha | \mathbf{u}, \beta, \mathbf{b} 
		& \sim \mathcal{G}( M/2 + c_{\alpha}, \| \mathbf{A} \mathbf{u} - \mathbf{b} \|_2^2/2 + d_{\alpha} ), \\
	\beta | \mathbf{u}, \alpha, \mathbf{b} 
		& \sim \mathcal{G}( N/2 + c_{\beta}, \| (\mathbf{F} - \mathbf{I}) \mathbf{u} \|_2^2/2 + d_{\beta} ), 
\end{aligned}
\end{equation}
where $\mathcal{N}$ and $\mathcal{G}$ again denote Gaussian and gamma distributions. 
Furthermore, the precision matrix in \cref{eq:cond_distributions} is $\mathbf{C} = \alpha \mathbf{A}^T \mathbf{A} + \beta (\mathbf{F}-\mathbf{I})^T (\mathbf{F}-\mathbf{I})$---the same as in \cref{eq:u_update_reg}. 
Building upon \cref{eq:cond_distributions}, we can use a Gibbs sampler that sequentially samples from the full conditional distributions to generate samples from the joint posterior $\pi^b(\mathbf{u}, \alpha, \beta)$. 
We summarize the proposed Gibbs sampler below in \cref{alg:Gibbs}.

\begin{algorithm}[htb]
\caption{Gibbs sampling for the Bayesian SIAC posterior}\label{alg:Gibbs}
\begin{algorithmic}
	\State{\textbf{Input:}} Initial values $\alpha_1$ and $\beta_1$, and total number of samples $J_{\text{samp}}$ 
	\For{$j=1,\dots,J_{\text{samp}}$}
		\State\textit{Step 1.} Compute $\mathbf{u}_{j} \sim \mathcal{N}( \mathbf{C}_j^{-1} \alpha_j \mathbf{A}^T \mathbf{b}, \mathbf{C}_j^{-1} )$, $\mathbf{C}_j = \alpha_j \mathbf{A}^T \mathbf{A} + \beta_j (\mathbf{F} - \mathbf{I})^T (\mathbf{F} - \mathbf{I})$
		 \State\textit{Step 2.} Compute $\alpha_{j+1} \sim \mathcal{G}( M/2 + c_{\alpha}, \| \mathbf{A} \mathbf{u}_j - \mathbf{b} \|_2^2/2 + d_{\alpha} )$
		 \State\textit{Step 3.} Compute $\beta_{j+1} \sim \mathcal{G}( N/2 + c_{\beta}, \| (\mathbf{F} - \mathbf{I}) \mathbf{u}_j \|_2^2/2 + d_{\beta} )$
	\EndFor
	\State{\textbf{Output:}} Samples $\{ (\mathbf{u}_j, \alpha_j, \beta_j ) \}_{j=1}^{J_{\text{samp}}}$
\end{algorithmic}
\end{algorithm} 

We note that generating the scalar random draws for $\alpha$ and $\beta$ in Steps 1 and 2 of \cref{alg:Gibbs} is efficient and easy using readily available software. 
In our implementation, we use the Julia programming language \cite{bezanson2017julia}. 
The computational bottleneck of \cref{alg:Gibbs} is Step 1, which demands more attention. 
For each iteration of $\mathbf{u}_j$ in Step 1 of \cref{alg:Gibbs}, we need to solve the linear system 
\begin{equation}\label{eq:Gibbs_u} 
	\mathbf{C}_j \mathbf{u}_{j} = \alpha_j \mathbf{A}^T \mathbf{b} + \mathbf{w} 
	\quad \text{with} \quad 
	\mathbf{w} \sim \mathcal{N}( \mathbf{0}, \mathbf{C}_j ),
\end{equation} 
where $\mathbf{C}_j = \alpha_j \mathbf{A}^T \mathbf{A} + \beta_j (\mathbf{F} - \mathbf{I})^T (\mathbf{F} - \mathbf{I})$. 
We can efficiently solve \cref{eq:Gibbs_u} using the methods discussed in \cref{rem:iter_methods}. 
Furthermore, the samples $\mathbf{w}$ needed to set up the right-hand side of \cref{eq:Gibbs_u} can be computed efficiently as $\mathbf{w} = \sqrt{\alpha_j} \mathbf{A}^T \mathbf{v}_1 + \sqrt{\beta_j} (\mathbf{F} - \mathbf{I})^T \mathbf{v}_2$ with $\mathbf{v}_1 \sim \mathcal{N}(\mathbf{0},\mathbf{I}_M)$ and $\mathbf{v}_2 \sim \mathcal{N}(\mathbf{0},\mathbf{I}_N)$. 
Here, $\mathbf{I}_M$ and $\mathbf{I}_N$ denote the identity matrices of dimension $M \times M$ and $N \times N$, respectively. 

%% file: 6_numerics.tex
 \section{Computational examples} 
\label{sec:numerics} 

We now demonstrate the proposed Bayesian SIAC filter for a series of computational experiments. 
The Julia \cite{bezanson2017julia} code to reproduce our computational experiments is openly available at \cite{glaubitz2025BayesianRepro}.

\subsection{Signal denoising} 
\label{sub:numerics_signal}

\begin{figure}[tb]
	\centering 
	\begin{subfigure}[b]{0.49\textwidth}
		\includegraphics[width=\textwidth]{%
			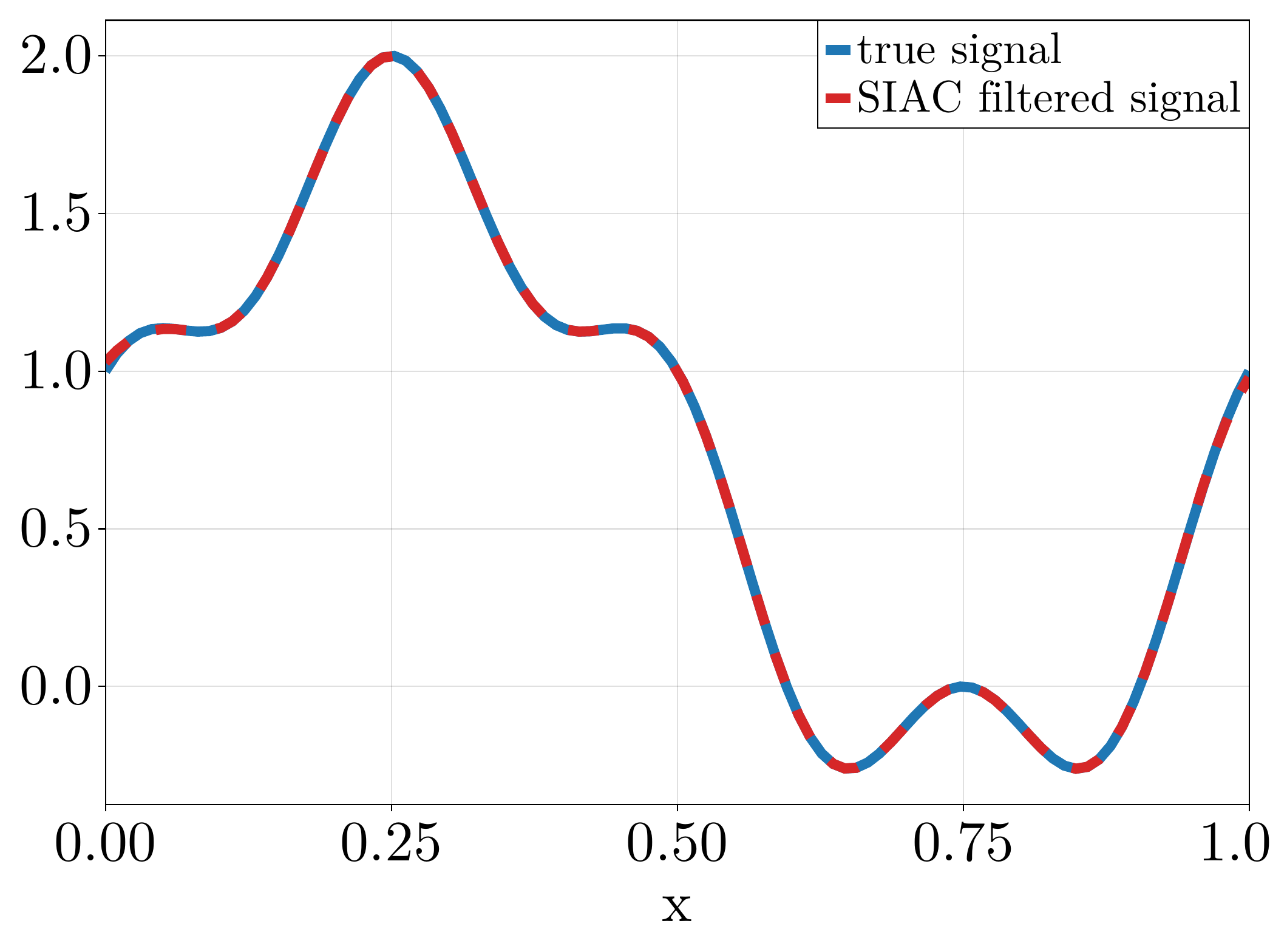}
    \caption{True signal $u$ and SIAC-filtered signal $\mathcal{F}[u]$}
    \label{fig:true_signal}
    \end{subfigure}%
	~
	\begin{subfigure}[b]{0.49\textwidth}
		\includegraphics[width=\textwidth]{%
			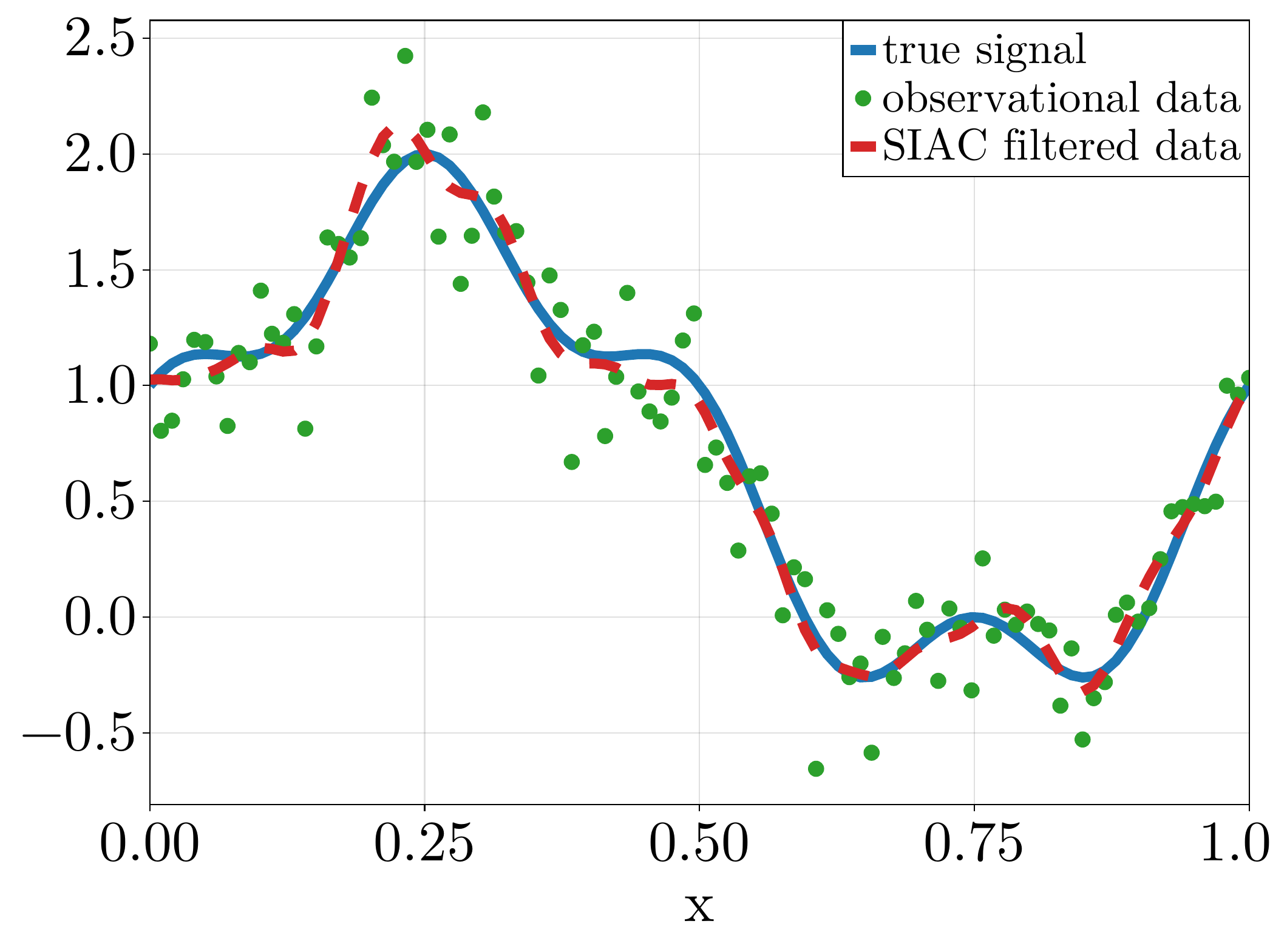}
    \caption{Data $\mathbf{b}$ and SIAC-filtered data $\mathbf{F} \mathbf{b}$}
    \label{fig:data}
    \end{subfigure}
    \caption{ 
	   Set-up for the one-dimensional denoising problem 
	}
    \label{fig:signal_data}
\end{figure}

We consider a one-dimensional denoising problem where the goal is to estimate the nodal values $\mathbf{u} = [ u(x_1), \dots, u(x_N) ]^T$ of a smooth signal $u: [0,1] \to \mathbb{R}$ at $N=100$ equidistant grid points from noisy direct observations. 
This setup leads to the data-generating model
\begin{equation}\label{eq:signal_denoising}
	\mathbf{b} = \mathbf{u} + \mathbf{e},
\end{equation}
where the forward operator is $\mathbf{A} = \mathbf{I}$. 
Furthermore, we assume i.i.d.\ Gaussian noise $\mathbf{e} \sim \mathcal{N}(\mathbf{0}, \sigma^2 \mathbf{I} )$ with noise variance $\sigma^2 = 5 \cdot 10^{-2}$. 
The true underlying signal is $u(x) = \sin(2\pi x)+0.5\cos(4\pi x)^2 + 0.5$. 
\Cref{fig:signal_data} illustrates the problem set-up. 
Specifically, \Cref{fig:true_signal} shows the true underlying signal $u$ and its SIAC-filtered version $\mathcal{F}[u]$. 
We use the B-Spline SIAC filter $K^{(7,4)}$, assuming $k=3$ as described in \cref{sec:SIAC-imple}.
We observe that the two are approximately the same, verifying that assumption \cref{eq:prior_assumption} is satisfied and thus motivating the application of the proposed SIAC prior in \Cref{sub:Bayesian_prior}. 
At the same time, \Cref{fig:data} illustrates the noisy observational data $\mathbf{b}$ and its SIAC-filtered version $\mathbf{F} \mathbf{b}$. 
The SIAC filter largely denoises the data, although small differences remain between the filtered data and the true signal.

\begin{figure}[tb]
	\centering 
	\begin{subfigure}[b]{0.49\textwidth}
		\includegraphics[width=\textwidth]{%
			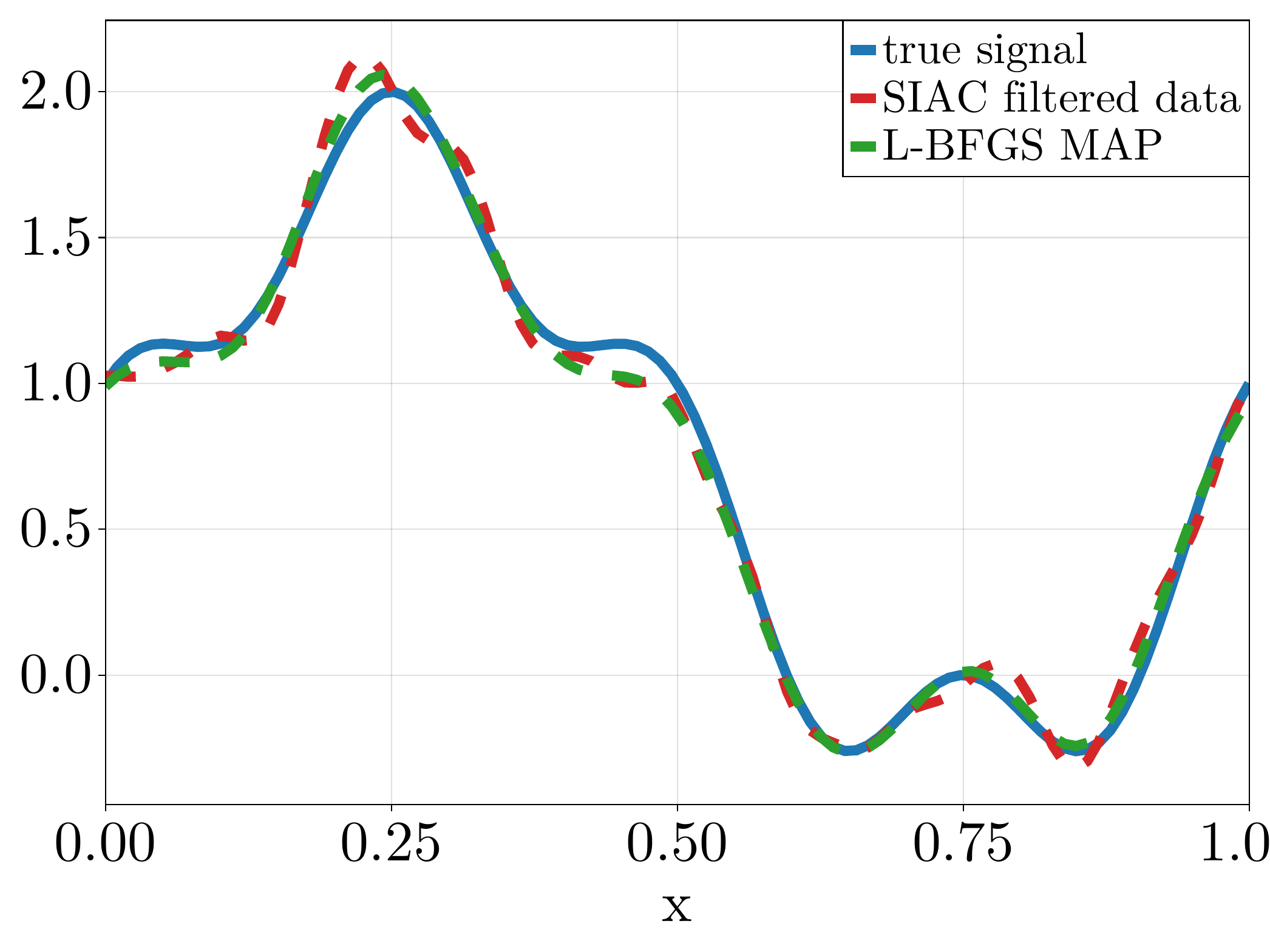} 
    \caption{MAP estimate using L-BFGS}
    \label{fig:signal_MAP_opt}
    \end{subfigure}%
	~
	\begin{subfigure}[b]{0.49\textwidth}
		\includegraphics[width=\textwidth]{%
			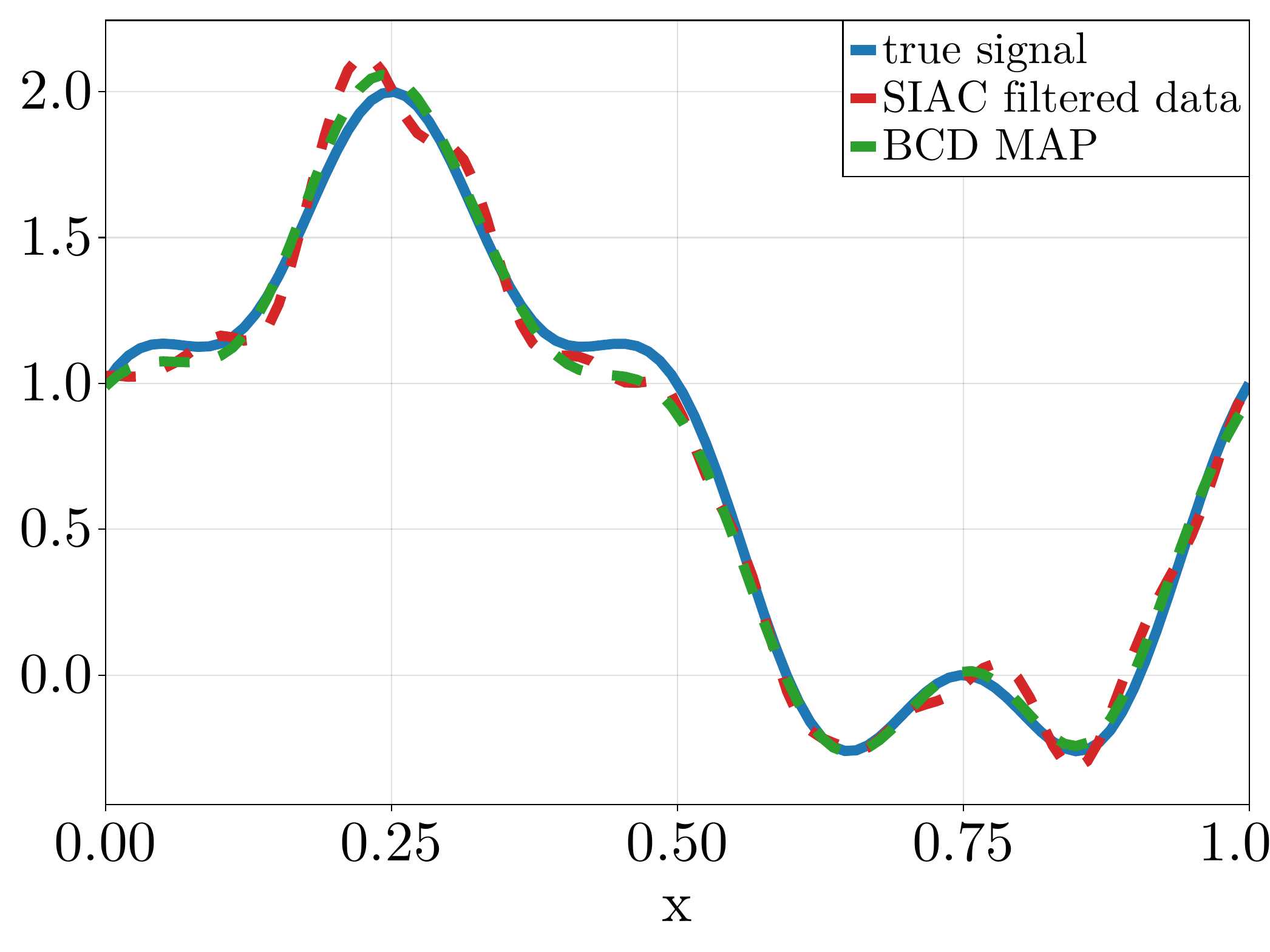} 
    \caption{MAP estimate using proposed BCD}
    \label{fig:signal_MAP_BCD}
    \end{subfigure}
    \caption{ 
	Comparison of the SIAC-filtered data $\mathbf{F} \mathbf{b}$ with the MAP estimate of the Bayesian SIAC filter $\mathbf{u}^{\rm MAP}$. 
    The MAP estimates were obtained using the existing L-BFGS (left) and proposed BCD (right) methods. 
	}
    \label{fig:signal_MAP_estimate}
\end{figure}

\Cref{fig:signal_MAP_estimate} illustrates the MAP estimates of the proposed Bayesian SIAC model. 
We determined these using (a) the `off-the-shelf' limited-memory Broyden--Fletcher--Goldfarb--Shanno (L-BFGS) algorithm \cite{nocedal2006numerical} and (b) the BCD approach proposed in \Cref{sec:MAP}. 
Both MAP estimates approximate the true signal roughly as accurately as the deterministic SIAC filter. 
At the same time, we do not observe visible accuracy differences between the MAP estimates (a) and (b).

\begin{table}[tb]
    \centering
    \begin{tabular}{c | c c | c c c}
        \toprule
        $N$ & \multicolumn{2}{c}{time [s]} & \multicolumn{3}{c}{relative $\ell^2$-error} \\ 
        \midrule
         & L-BFGS & BCD & L-BFGS & BCD & SIAC filtered data\\
        \midrule
        100 & 1.0e-01 & 1.7e-02 & 1.3e-01 & 4.8e-02 & 7.2e-02\\
        200 & 5.4e-01 & 8.4e-03 & 1.4e-01 & 6.4e-02 & 8.7e-02\\
        400 & 1.0e+01 & 1.4e-02 & 9.5e-02 & 5.5e-02 & 9.3e-02\\
        800 & 1.7e+01 & 5.4e-02 & 8.6e-02 & 5.6e-02 & 9.0e-02\\
        \bottomrule
    \end{tabular}
    \caption{
    Run times (in seconds) and relative $\ell^2$-errors of the Bayesian SIAC filter's MAP estimate---using the existing L-BFGS and proposed BCD algorithms---and the deterministic SIAC filter for different numbers of equidistant grid points $N$. 
    }
    \label{tab:signal_MAP_results}
\end{table}

In \Cref{tab:signal_MAP_results}, we quantify the above observations and further compare the existing L-BFGS method and the proposed BCD algorithm for approximating the MAP estimate of the Bayesian SIAC filter for different numbers $N$ of equidistant grid points, reporting execution times and relative $\ell^2$-errors. 
Furthermore, we compare the relative $\ell^2$-errors of the Bayesian SIAC filter's MAP estimates with the ones for the existing deterministic SIAC filter and observe that the BCD approach yields the smallest errors in all cases. 
Both MAP estimators used a maximum number of $10^3$ iterations. 
In the BCD algorithm, we also terminate if the relative or absolute change in $\mathbf{u}$ from one iteration to the next is below $10^{-4}$ or $10^{-8}$, respectively. 
For L-BFGS, we use the default setting of stopping when the absolute gradient in the infinity norm drops below $10^{-8}$. 
We did not fine-tune these stopping criteria for this particular problem. 
Yet, we observe that the proposed BCD algorithm is significantly faster and slightly more accurate than L-BFGS. 
This improvement can be attributed to the BCD algorithm’s ability to exploit the specific structure of the Bayesian SIAC filter's posterior, underscoring its advantage in efficient MAP estimation over the off-the-shelf optimizer.

We next illustrate how the proposed Bayesian SIAC filter enables UQ. 
In practice, uncertainty arises from both input data and the underlying data-generating model, and this uncertainty propagates through the reconstruction to affect subsequent predictions and decision-making. 
Quantifying such uncertainty is often critical in real-world applications. 
Here, we demonstrate how the Bayesian SIAC filter addresses this need by enabling UQ in the reconstructed signal.
Specifically, we sample from the Bayesian SIAC filter's posterior \cref{eq:posterior} using two methods:
\begin{itemize}
    \item the Gibbs sampler derived in \Cref{sec:MCMC}, and 
    \item the standard Adaptive Metropolis (AM) algorithm \cite{haario2001adaptive,andrieu2008tutorial,atchade2010limit}, implemented as in \cite{AdaptiveMCMC}, with a target mean acceptance rate of 23.4\%.
\end{itemize}
We include the AM algorithm as a reference point for the proposed Gibbs sampler. 
For both methods, we generated four independent MCMC chains in parallel (via multi-threading). 
All chains were initialized with random draws from the prior. 
Furthermore, we don't use thinning, but allow for a burn-in phase corresponding to 10\% of the number of samples.

\begin{figure}[tb]
	\centering 
	\begin{subfigure}[b]{0.49\textwidth}
		\includegraphics[width=\textwidth]{%
			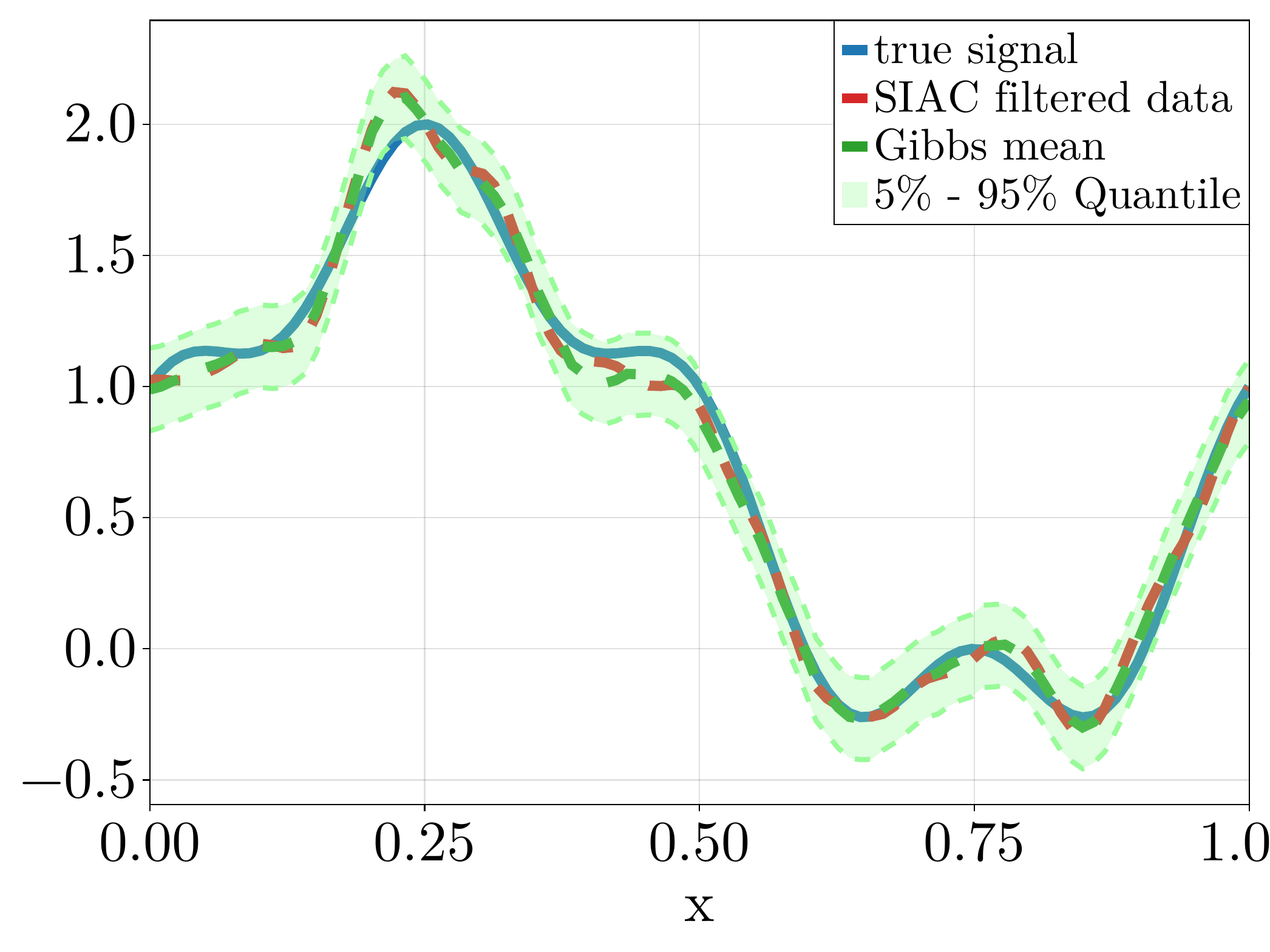}
    \caption{Gibbs sampling}
    \label{fig:signal_UQ_Gibbs}
    \end{subfigure}%
	~
	\begin{subfigure}[b]{0.49\textwidth}
		\includegraphics[width=\textwidth]{%
			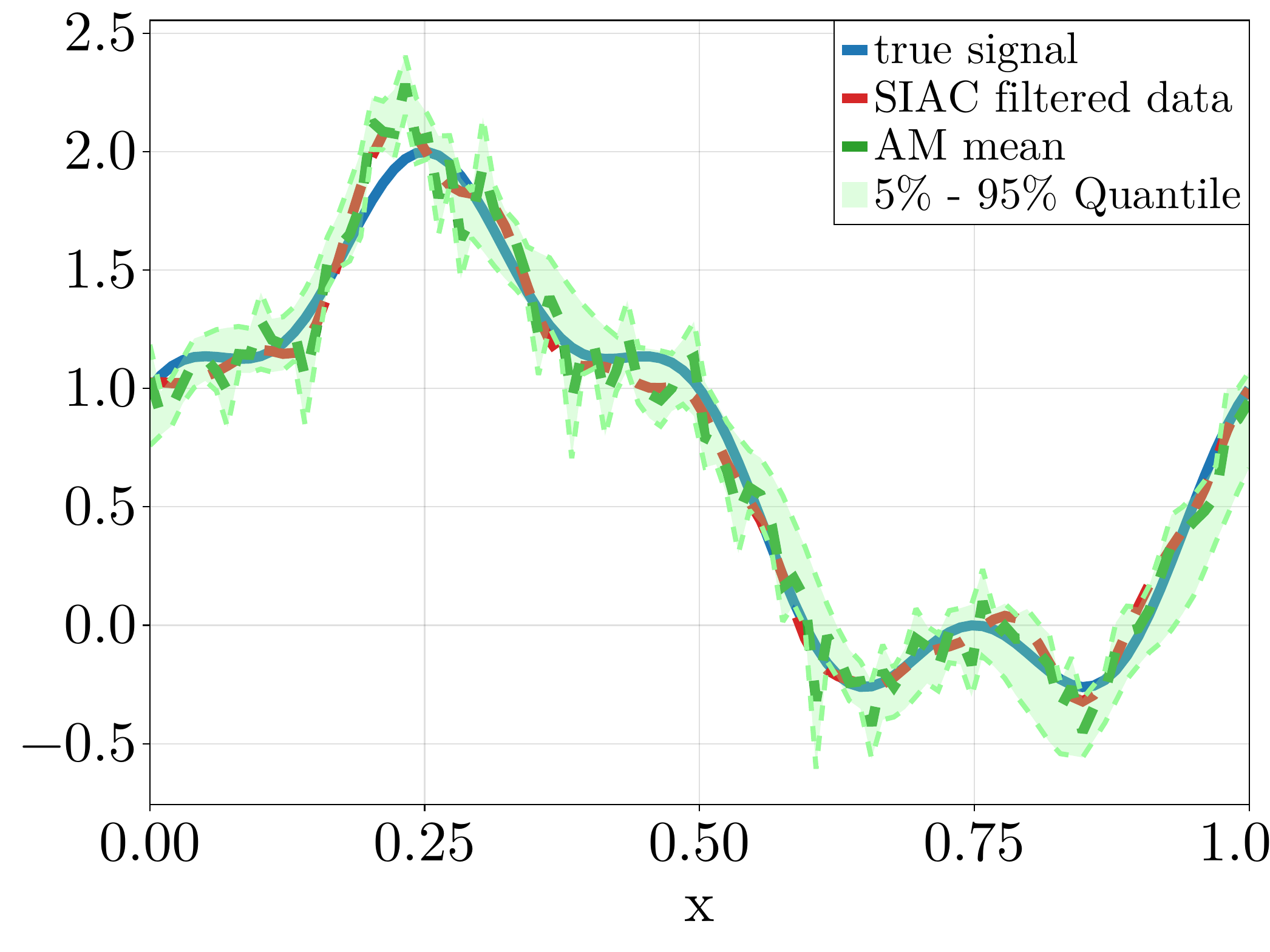}
    \caption{AM sampling}
    \label{fig:signal_UQ_AM}
    \end{subfigure}
    \caption{
        True signal $u$, SIAC-filtered noisy data $\mathbf{F} \mathbf{b}$, and mean and $90\%$ quantile ranges of the $u$--samples of the Bayesian SIAC filter. 
        We used the Gibbs sampler (\cref{fig:signal_UQ_Gibbs}) proposed in \Cref{sec:MCMC} and the standard AM sampler (\cref{fig:signal_UQ_AM}) to generate $10^5$ samples, respectively. 
	}
    \label{fig:signal_UQ}
\end{figure}

\cref{fig:signal_UQ} shows the mean and 90\% quantile range of the $u$-samples produced by the proposed Gibbs sampler (\cref{fig:signal_UQ_Gibbs}) and the standard AM sampler (\cref{fig:signal_UQ_AM}). 
We respectively generated $10^5$ samples.
Notably, the mean estimate obtained from the Gibbs sampler is comparable to both the MAP estimate and the reconstruction produced by the existing deterministic SIAC filter (see \cref{fig:signal_data,fig:signal_MAP_estimate}).
Yet, a key advantage of sampling from the Bayesian SIAC filter's posterior is the ability to compute additional statistics---such as the illustrated quantile ranges---which enables us to move beyond deterministic point estimates and quantify the uncertainty in the reconstructed function. 
By contrast, in \cref{fig:signal_UQ_AM}, the standard AM sampler produces a less accurate, oscillatory profile, suggesting that it explores the Bayesian SIAC filter's posterior less effectively than the proposed Gibbs approach.

\begin{figure}[tb] 
    \centering
    \begin{subfigure}[b]{0.32\textwidth}
		\includegraphics[width=\textwidth]{%
			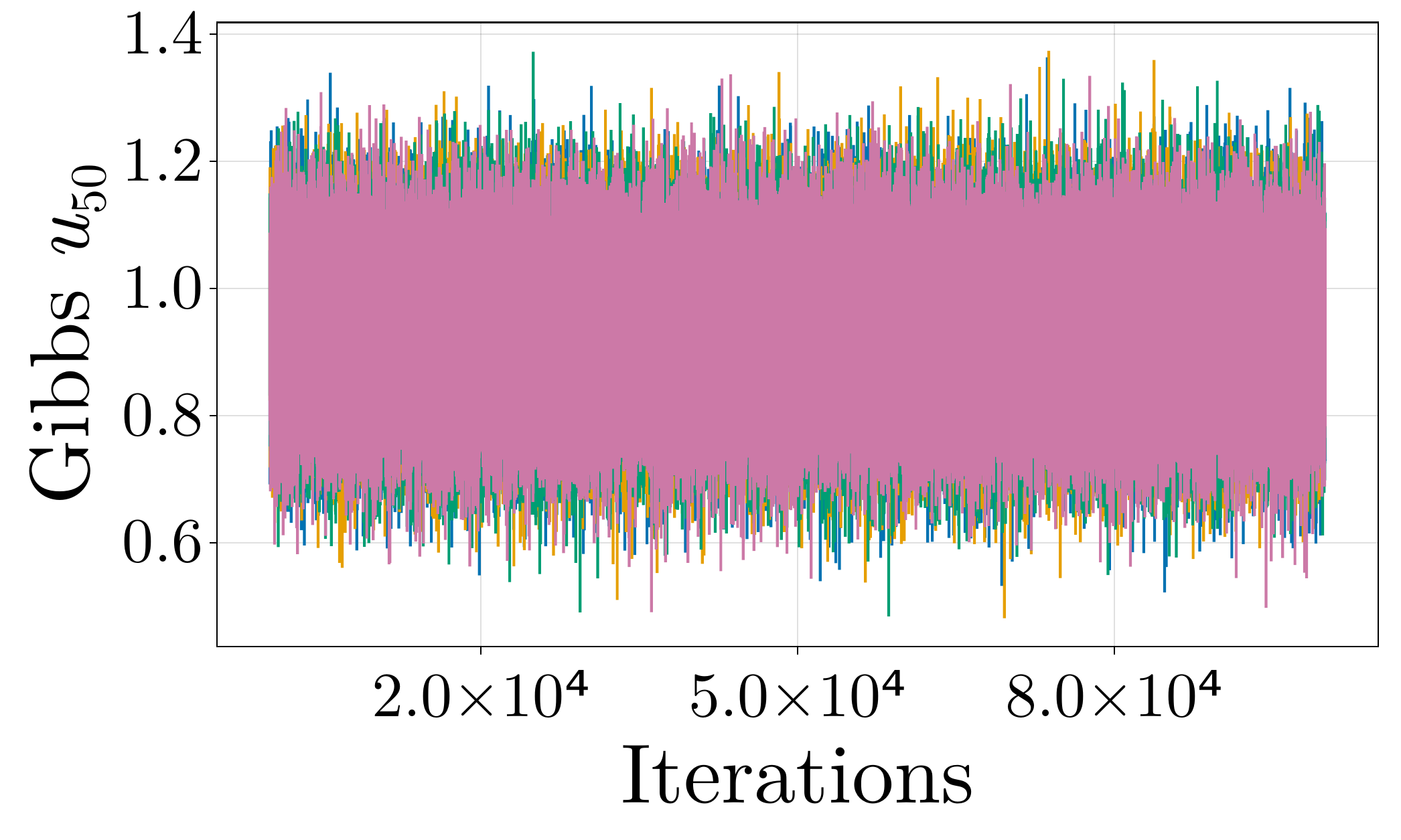} 
        \caption{Gibbs, $u_{50}$}
        \label{fig:signal_traces_Gibbs_u}
    \end{subfigure}
    \begin{subfigure}[b]{0.32\textwidth}
		\includegraphics[width=\textwidth]{%
			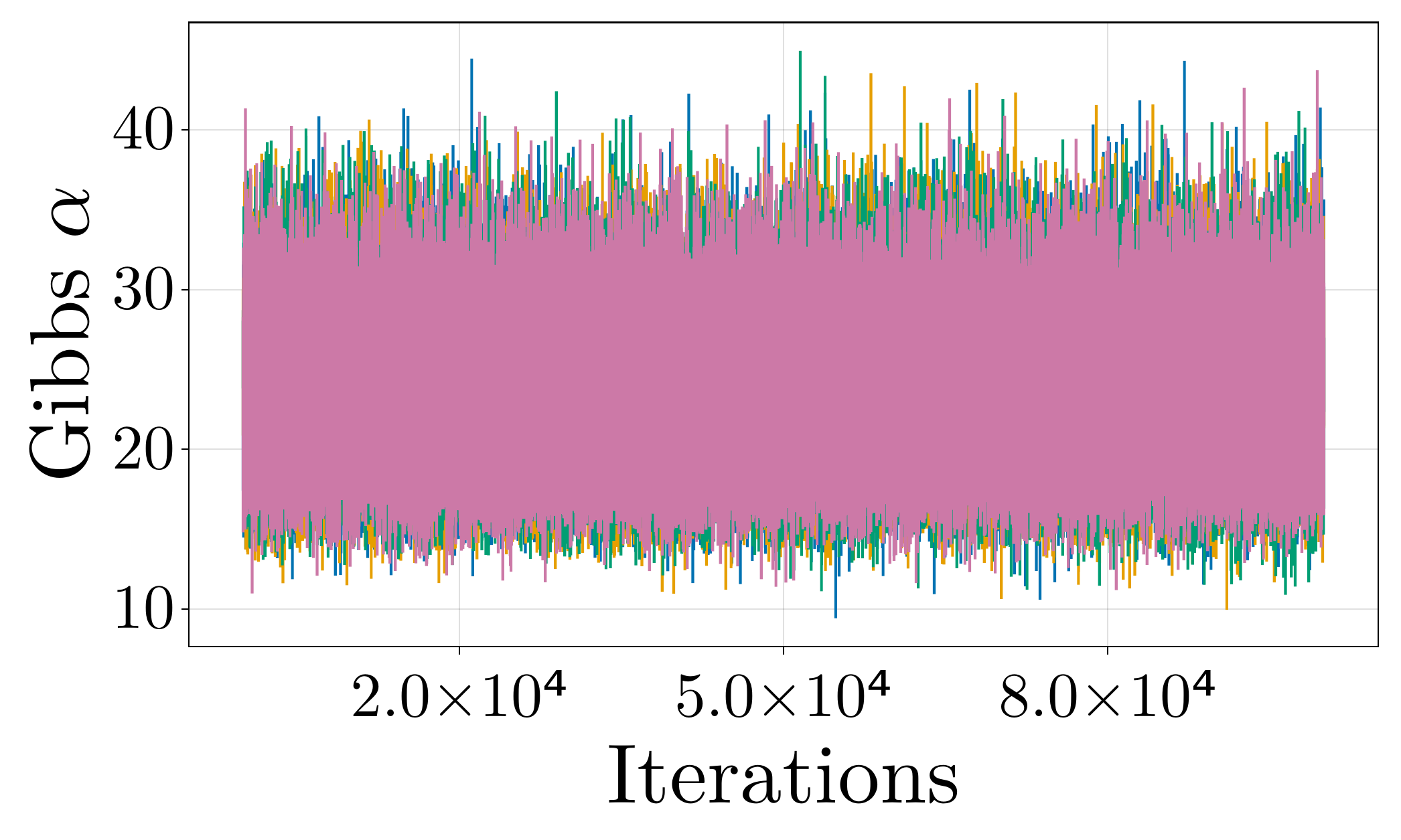} 
        \caption{Gibbs, $\alpha$}
        \label{fig:signal_traces_Gibbs_alpha}
    \end{subfigure}
    \begin{subfigure}[b]{0.32\textwidth}
		\includegraphics[width=\textwidth]{%
			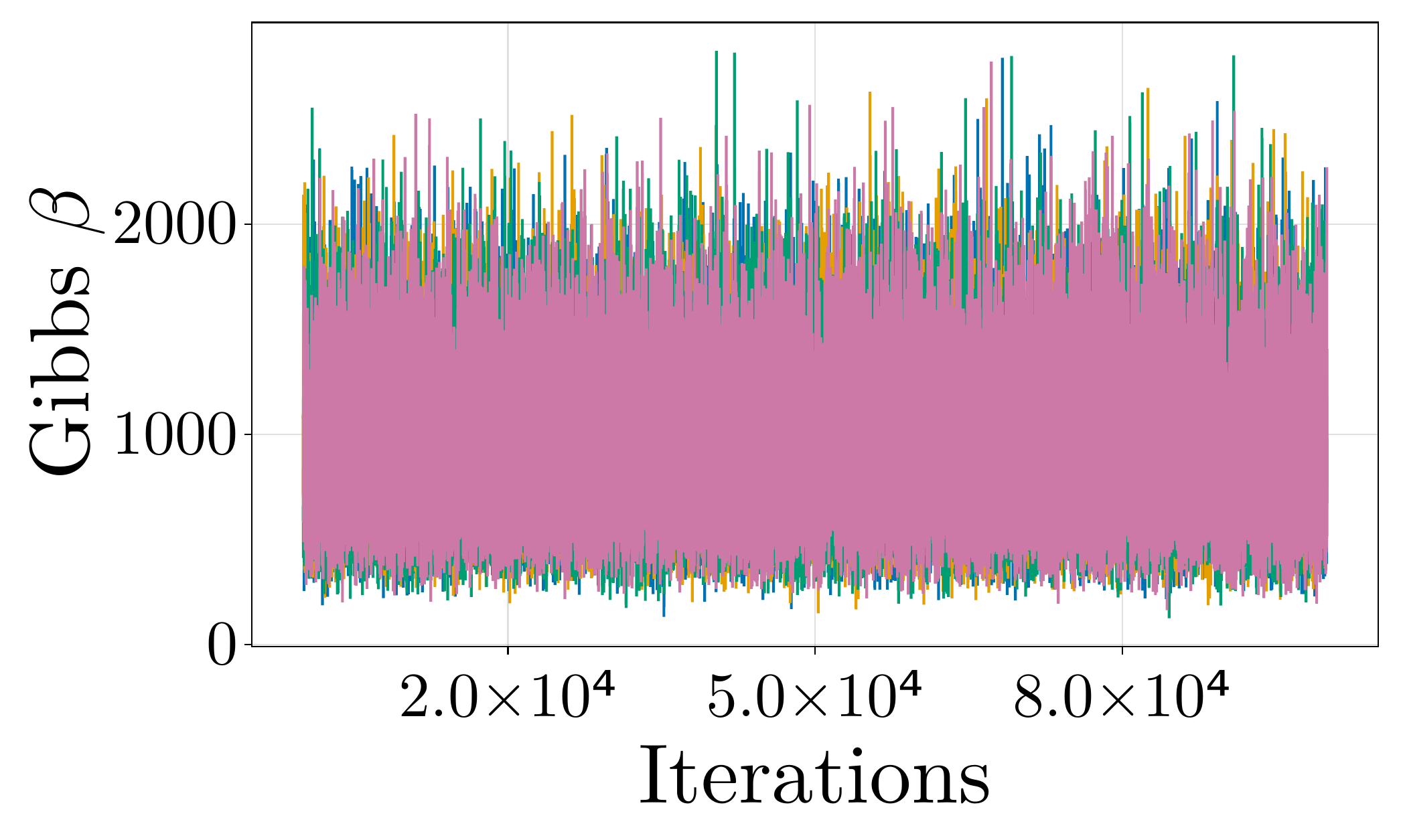} 
        \caption{Gibbs, $\beta$}
        \label{fig:signal_traces_Gibbs_beta}
    \end{subfigure}
    \\ 
    \begin{subfigure}[b]{0.32\textwidth}
		\includegraphics[width=\textwidth]{%
			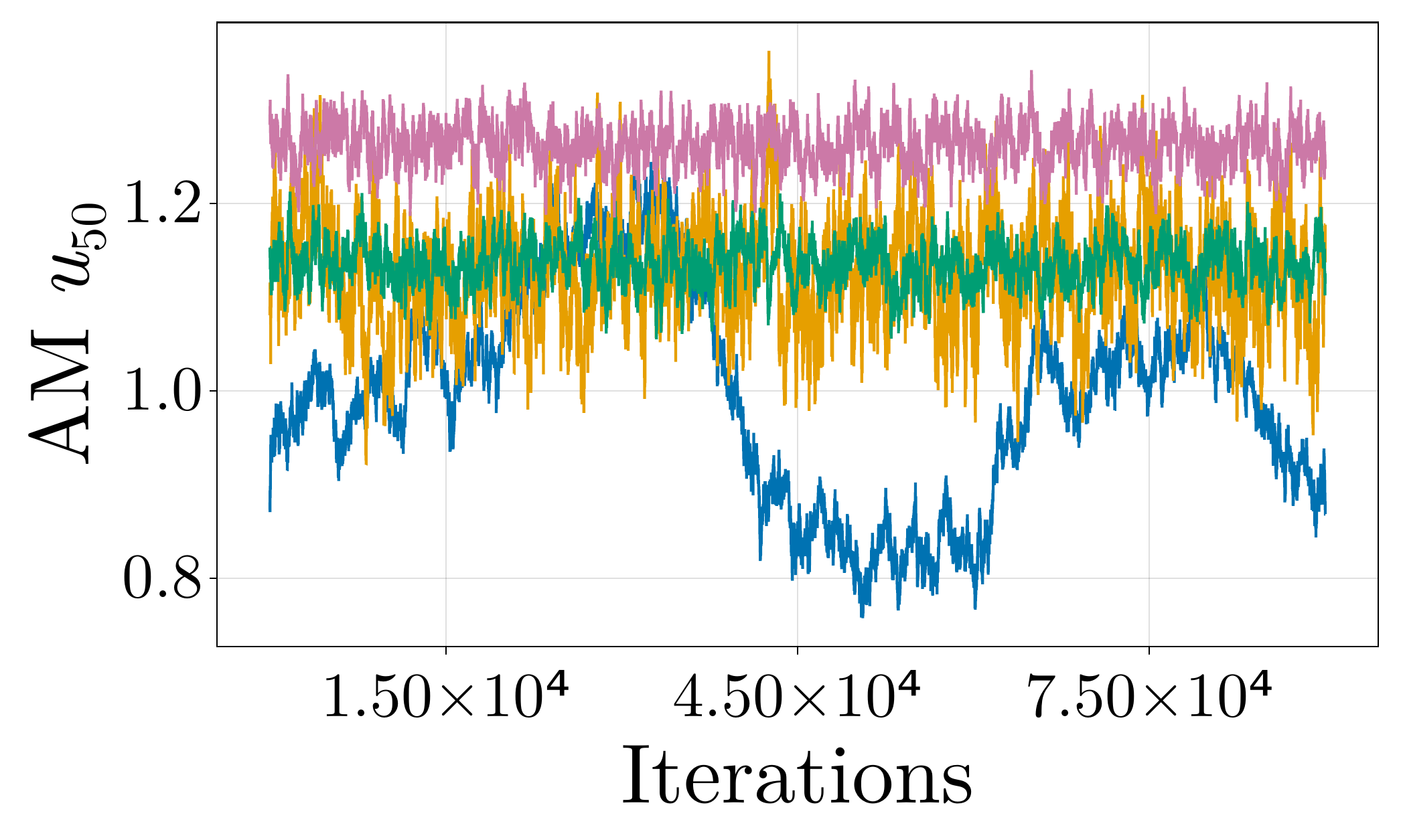} 
        \caption{AM, $u_{50}$}
        \label{fig:signal_traces_AM_u}
    \end{subfigure}
    \begin{subfigure}[b]{0.32\textwidth}
		\includegraphics[width=\textwidth]{%
			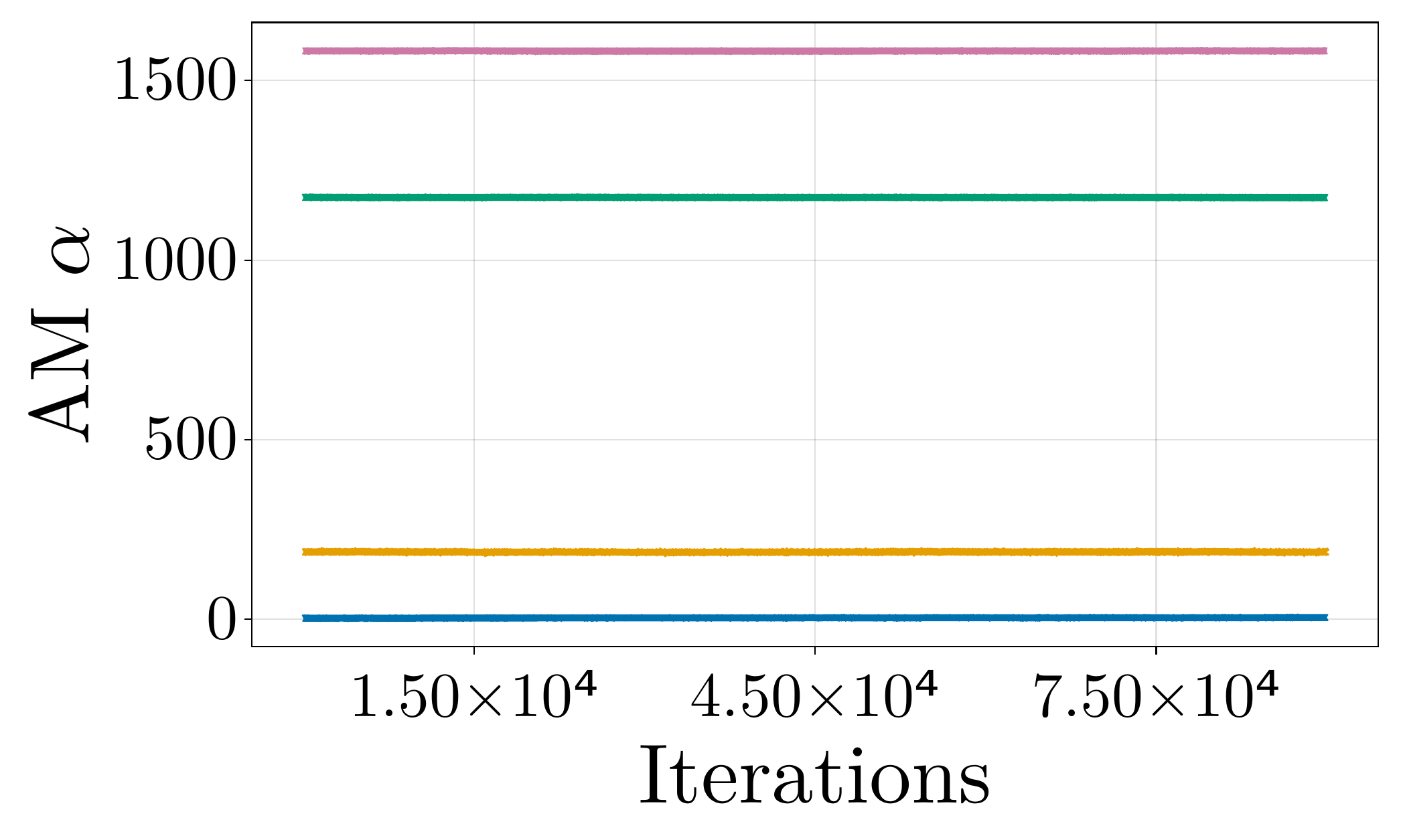} 
        \caption{AM, $\alpha$}
        \label{fig:signal_traces_AM_alpha}
    \end{subfigure}
    \begin{subfigure}[b]{0.32\textwidth}
		\includegraphics[width=\textwidth]{%
			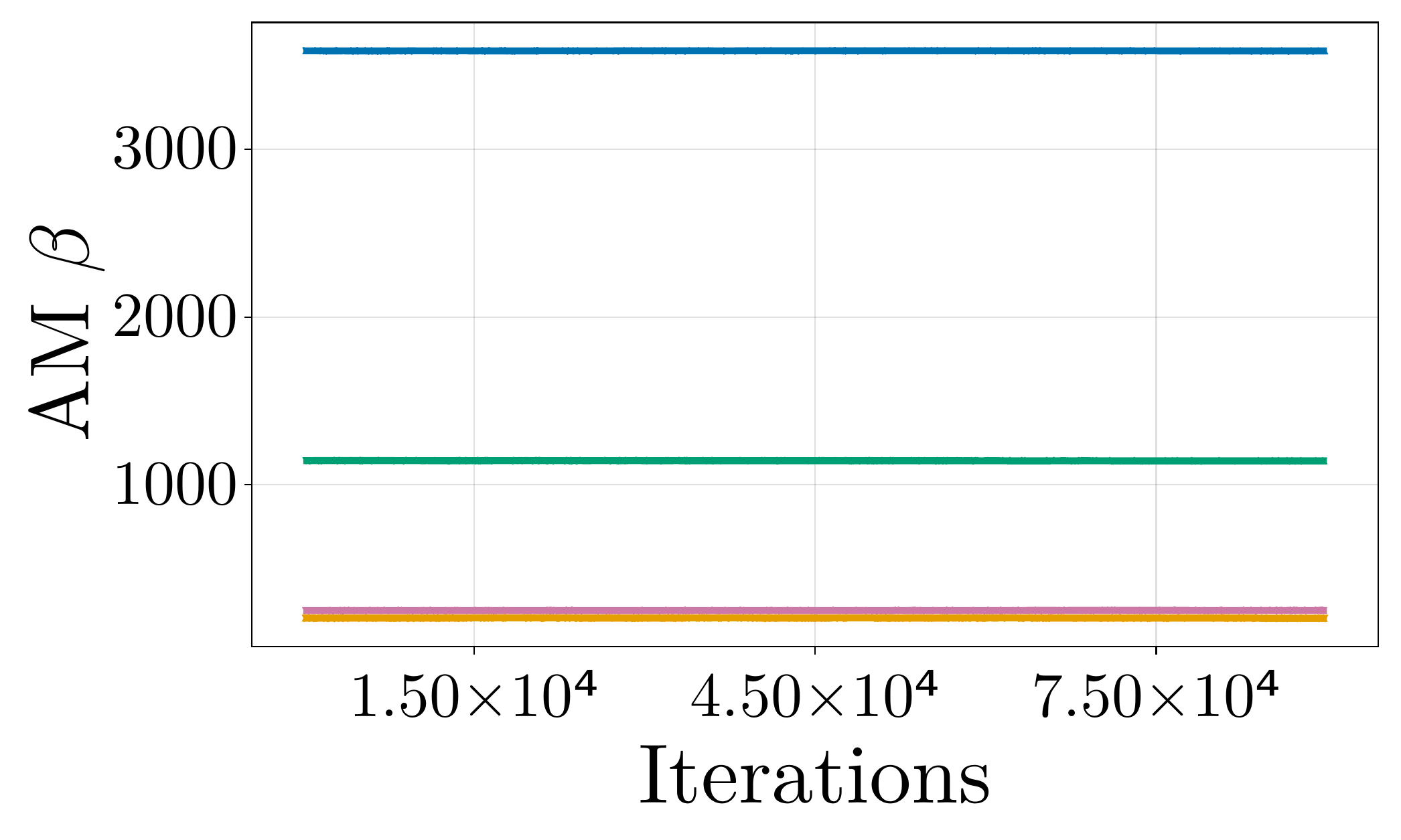} 
        \caption{AM, $\beta$}
        \label{fig:signal_traces_AM_beta}
    \end{subfigure}
    \caption{ 
    Traces for the $u_{50}$- (left), $\alpha$- (middle), and $\beta$-samples (right) of the Bayesian SIAC filter.
    We used the Gibbs sampler from \Cref{sec:MCMC} (top) and the AM algorithm (bottom).
	}
    \label{fig:signal_traces}
\end{figure}

To further compare the two samplers, \cref{fig:signal_traces} displays the trace plots for the $u_{50}$-, $\alpha$-, and $\beta$-samples of the Bayesian SIAC filter. 
The results in \cref{fig:signal_traces} underscore that the proposed Gibbs sampler explores the Bayesian SIAC filter's posterior more effectively than the AM method. 
Specifically, in \cref{fig:signal_traces_AM_alpha,fig:signal_traces_AM_beta}, we observe the AM chains to be stuck in different modes, exhibiting poor mixing. 
By contrast, in \cref{fig:signal_traces_Gibbs_alpha,fig:signal_traces_Gibbs_beta}, each Gibbs chain appears to traverse the posterior space efficiently.

\begin{table}[tb]
    \centering
    \begin{tabular}{c c r r r r r}
        \toprule
        & \#samples & time & MPSRF-1 & rel.\ MPSRF & ESS & rel.\ ESS \\
        \midrule 
        Gibbs & $10^4$ & 8.5e+0 & 8.4e-3 & 7.2e-2 & 3.8e+4 & 4.6e+4 \\
        Gibbs & $10^5$ & 3.7e+1 & 5.3e-4 & 2.0e-2 & 3.8e+5 & 1.0e+4 \\
        AM & $10^5$ & 1.6e+0 & 7.8e+5 & 1.3e+6 & 8.6e+2 & 5.3e+2 \\
        AM & $10^6$ & 1.3e+1 & 5.7e+5 & 7.6e+5 & 8.4e+3 & 6.3e+2 \\
        \bottomrule
    \end{tabular}
    \caption{
    The wall-clock runtime in seconds, the number of samples, the MPSRF minus one, relative MPSRF (MPSRF-1, times the runtime), the mean ESS, and the relative mean ESS per second of the Bayesian SIAC model using the Gibbs sampling and the AM sampling. 
    } 
    \label{tab:signal}
\end{table}

We observe comparable behavior across different sample sizes for both sampling methods. 
To demonstrate this, \cref{tab:signal} summarizes the number of samples, the computational wall-clock time in seconds (“time”), the multivariate potential scale reduction factor (MPSRF) minus one, relative MPSRF (MPSRF-1, times the runtime), the mean effective sample size (ESS), and the relative mean ESS per second for the proposed Gibbs sampler and the standard AM sampler. 
The ESS of a Markov chain is defined as the number of independent samples of the posterior that are needed to estimate $\mathbb{E}[ \mathcal{G} ]$ (for some quantity of interest $\mathcal{G}(\mathbf{z})$) with the same statistical accuracy as an estimate from the Markov chain.
It is a measure of how much information is contained in the MCMC chain and should be as large as possible; see \cite{gelman2003bayesian,wolff2004monte}. 
We observe in \cref{tab:signal} that the mean ESS and the relative mean ESS of the proposed Gibbs sampler are around two orders of magnitude larger than for the AM algorithm, indicating that the Gibbs sampler explores the Bayesian SIAC filter's posterior more efficiently.  
At the same time, the MPSRF is commonly used to assess the convergence of multiple MCMC chains \cite{brooks1998general}.
One typically has $\operatorname{MPSRF} \geq 1$ (assuming the chains have overdispersed starting points that cause the inter-chain variance to be larger than the within-chain variance). 
When the MPSRF approaches $1$, the variance within each sequence approaches the variance across sequences, thus indicating that each individual chain has converged to the target distribution. 
The literature contains several recommendations for the values of MPSRF that indicate convergence. 
For example, \cite{gelman1992inference} suggests the commonly used value $\operatorname{MPSRF} < 1.1$ while \cite{vehtari2021rank} argues for the more conservative threshold $\operatorname{MPSRF} < 1.01$.
We observe in \cref{tab:signal} that the MPSRF and the relative MPSRF of the proposed Gibbs sampler are around six and five orders of magnitude smaller than for the AM algorithm, indicating that the Gibbs sampler converges significantly faster. 
Specifically, even when we compare Gibbs for $10^4$ samples with AM for $10^6$ samples---in which case Gibbs takes only 65\% of AM's computational time---we observe the ``MPSRF-1'' and ``ESS'' values to be eight and one order of magnitude smaller and larger for Gibbs compared to AM, respectively, again indicating that the proposed Gibbs sampler explores the Bayesian SIAC filter's posterior more effectively. 

\begin{remark}
We briefly discuss the computational cost and scalability of the Bayesian SIAC filter.
For MAP estimation, the dominant computational cost arises from repeatedly solving the linear system \cref{eq:u_update_reg} of size $N \times N$, for which the computational and memory costs have been discussed in \cref{rem:iter_methods}. 
As shown in \cref{tab:signal_MAP_results}, the proposed BCD approach exhibits favorable scaling behavior with respect to $N$ for the tested problem sizes and provides substantial runtime improvements over a generic L-BFGS optimizer.
For posterior sampling, the computational cost is again dominated by solving similar linear systems at each iteration of the Gibbs sampler.
While the Gibbs sampler has higher per-sample costs, its improved mixing behavior yields significantly smaller MPSRFs and larger ESS per unit time than AM sampling, as demonstrated in \cref{tab:signal}.
\end{remark}

\subsection{Robustness with respect to the kernel}

We next compare the robustness of the deterministic and proposed Bayesian SIAC filters $\mathcal{F}$ with respect to the kernel $K$.
To this end, we consider the same set-up as above in \Cref{sub:numerics_signal} with different B-spline SIAC filters using kernels $K = K^{(2k+1,k+1)}$ of degree $k$ as described in detail in \cref{sec:SIAC-imple}. 

\begin{figure}[tb]
	\centering 
	\begin{subfigure}[b]{0.49\textwidth}
		\includegraphics[width=\textwidth]{%
			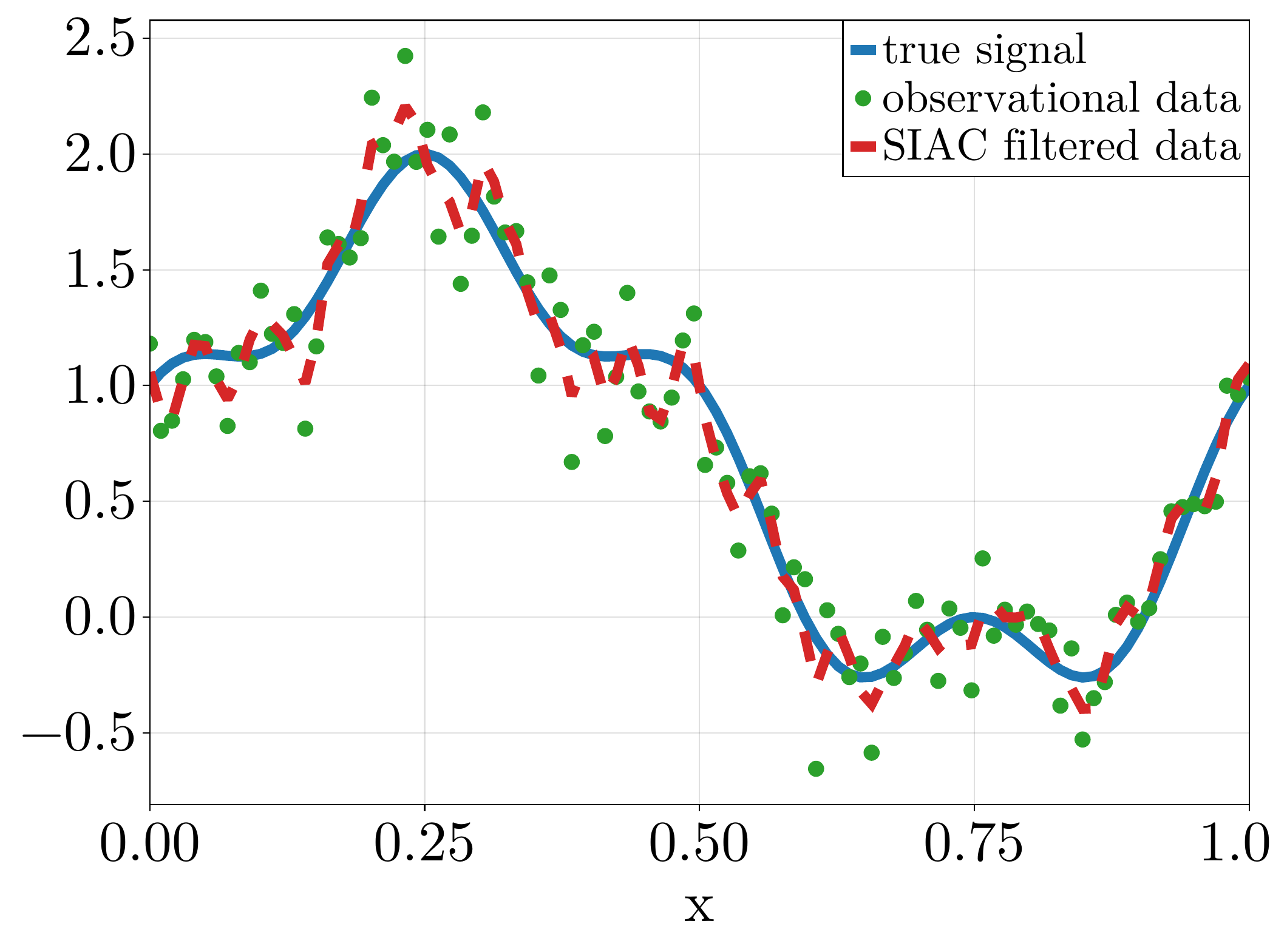}
    \caption{$k=1$}
    \label{fig:differentKernels_data_k1}
    \end{subfigure}%
	~
	\begin{subfigure}[b]{0.49\textwidth}
		\includegraphics[width=\textwidth]{%
			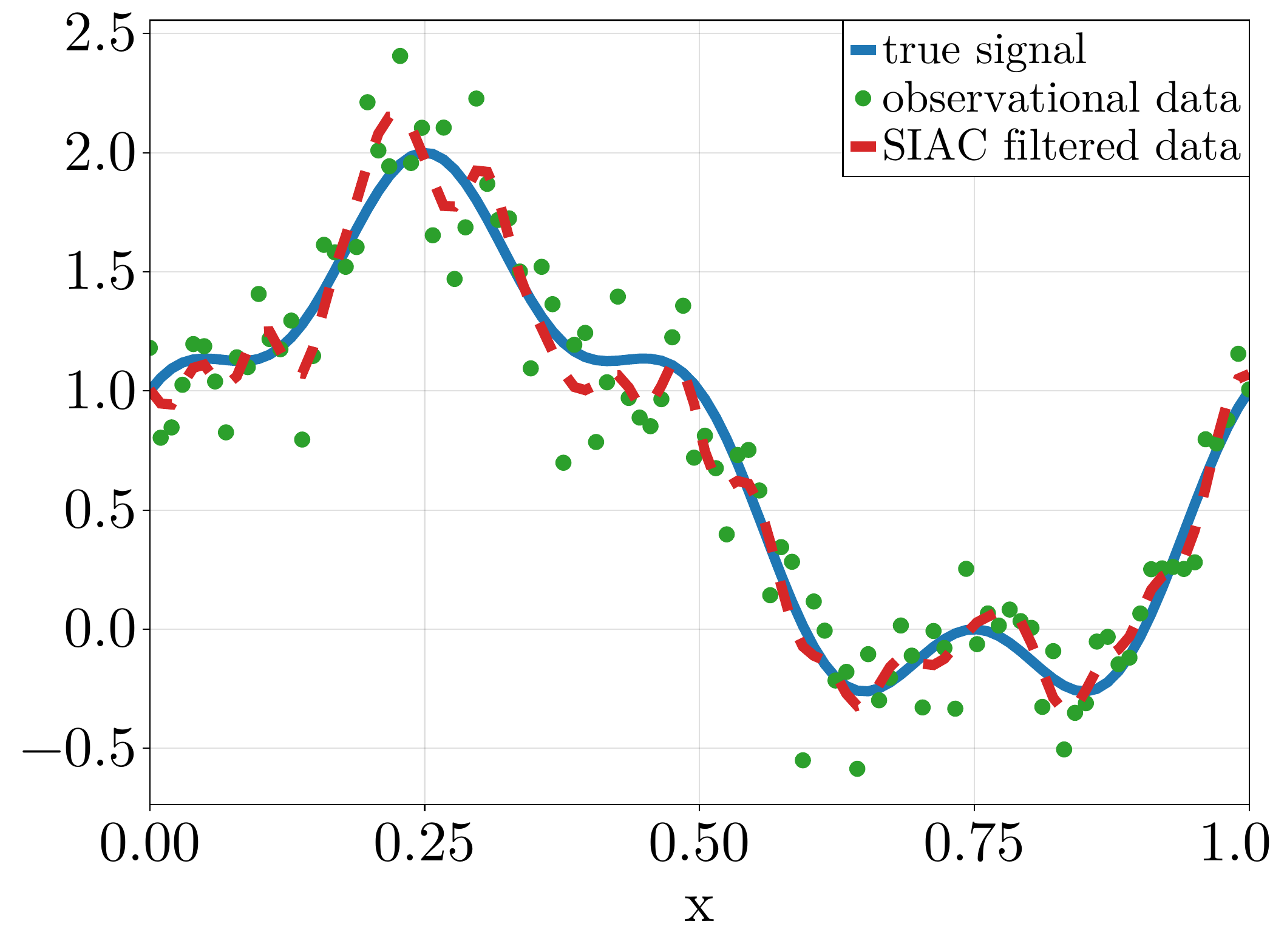}
    \caption{$k=2$}
    \label{fig:differentKernels_data_k2}
    \end{subfigure}%
    \\ 
    \begin{subfigure}[b]{0.49\textwidth}
		\includegraphics[width=\textwidth]{%
			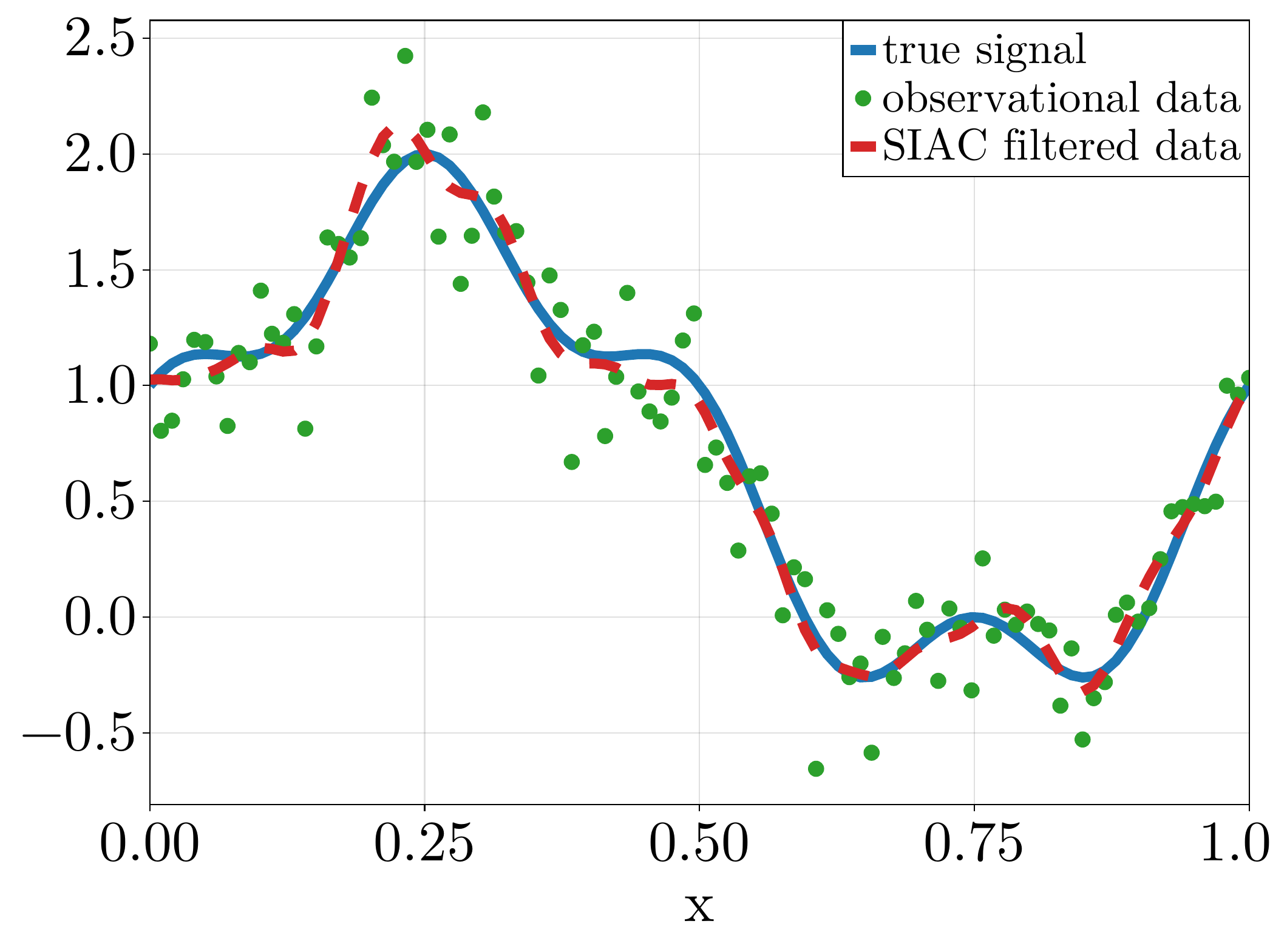}
    \caption{$k=3$}
    \label{fig:differentKernels_data_k3}
    \end{subfigure}%
    ~ 
    \begin{subfigure}[b]{0.49\textwidth}
		\includegraphics[width=\textwidth]{%
			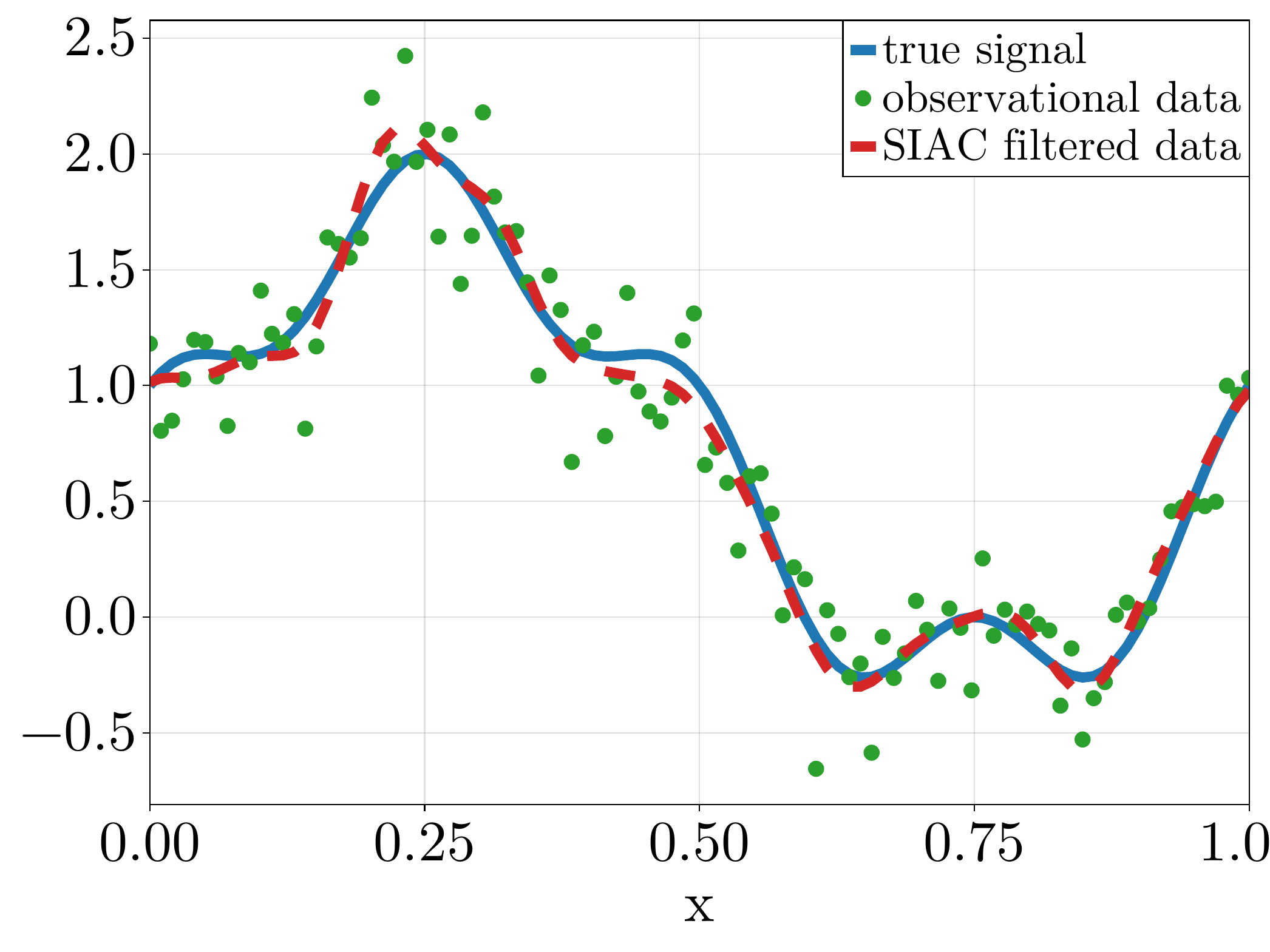}
    \caption{$k=4$}
    \label{fig:differentKernels_data_k4}
    \end{subfigure}%
    \caption{
        The true signal $u$, the observational data $\mathbf{b}$, and the SIAC-filtered data $\mathbf{F} \mathbf{b}$ for SIAC filters $\mathcal{F}$ with different kernels $K^{(2k+1,k+1)}$, where $k$ is the degree
	}
    \label{fig:differentKernels_data}
\end{figure}

\begin{figure}[tb]
	\centering 
	\begin{subfigure}[b]{0.49\textwidth}
		\includegraphics[width=\textwidth]{%
			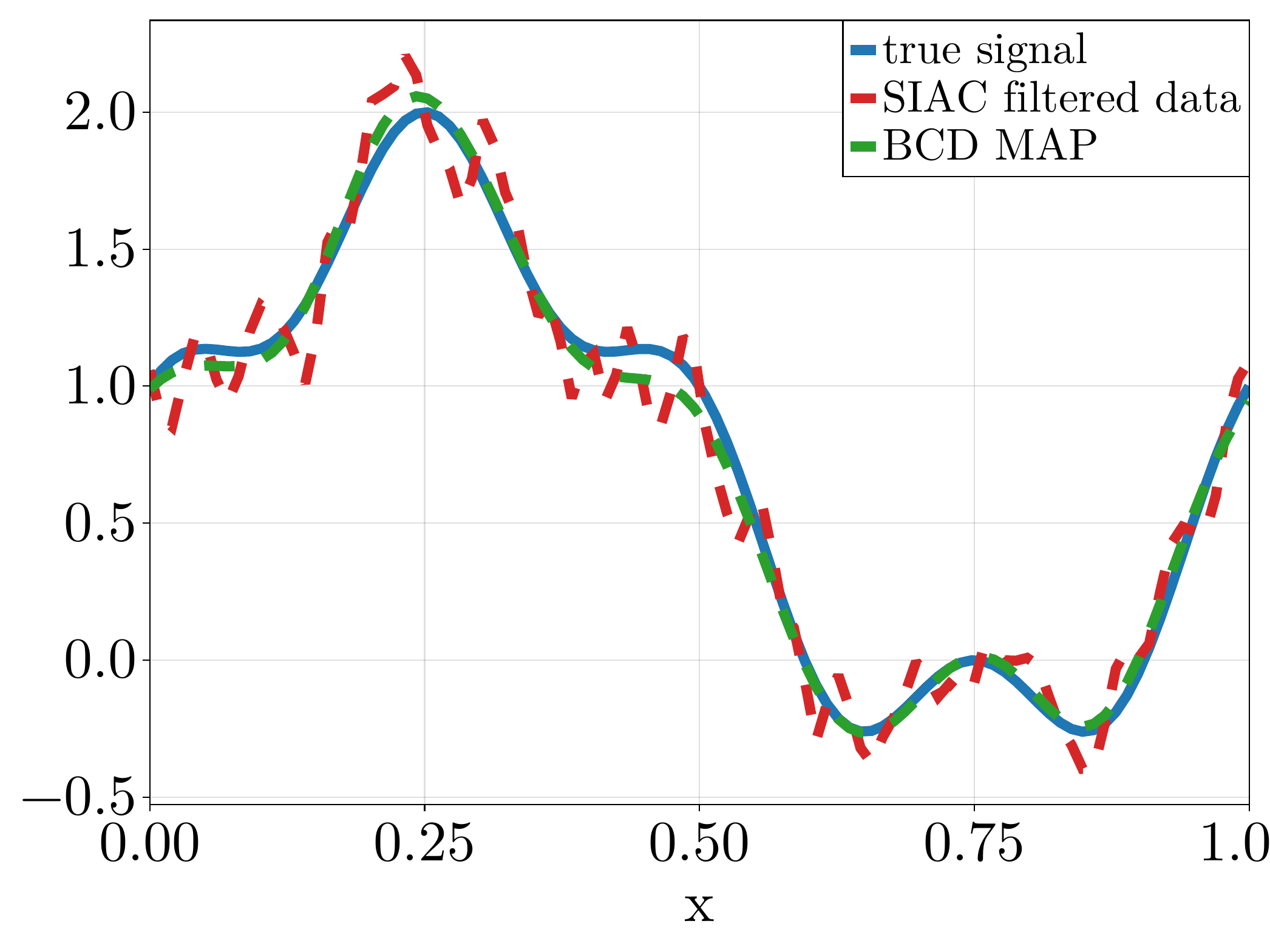}
    \caption{$k=1$}
    \label{fig:differentKernels_MAP_k1}
    \end{subfigure}%
	~
	\begin{subfigure}[b]{0.49\textwidth}
		\includegraphics[width=\textwidth]{%
			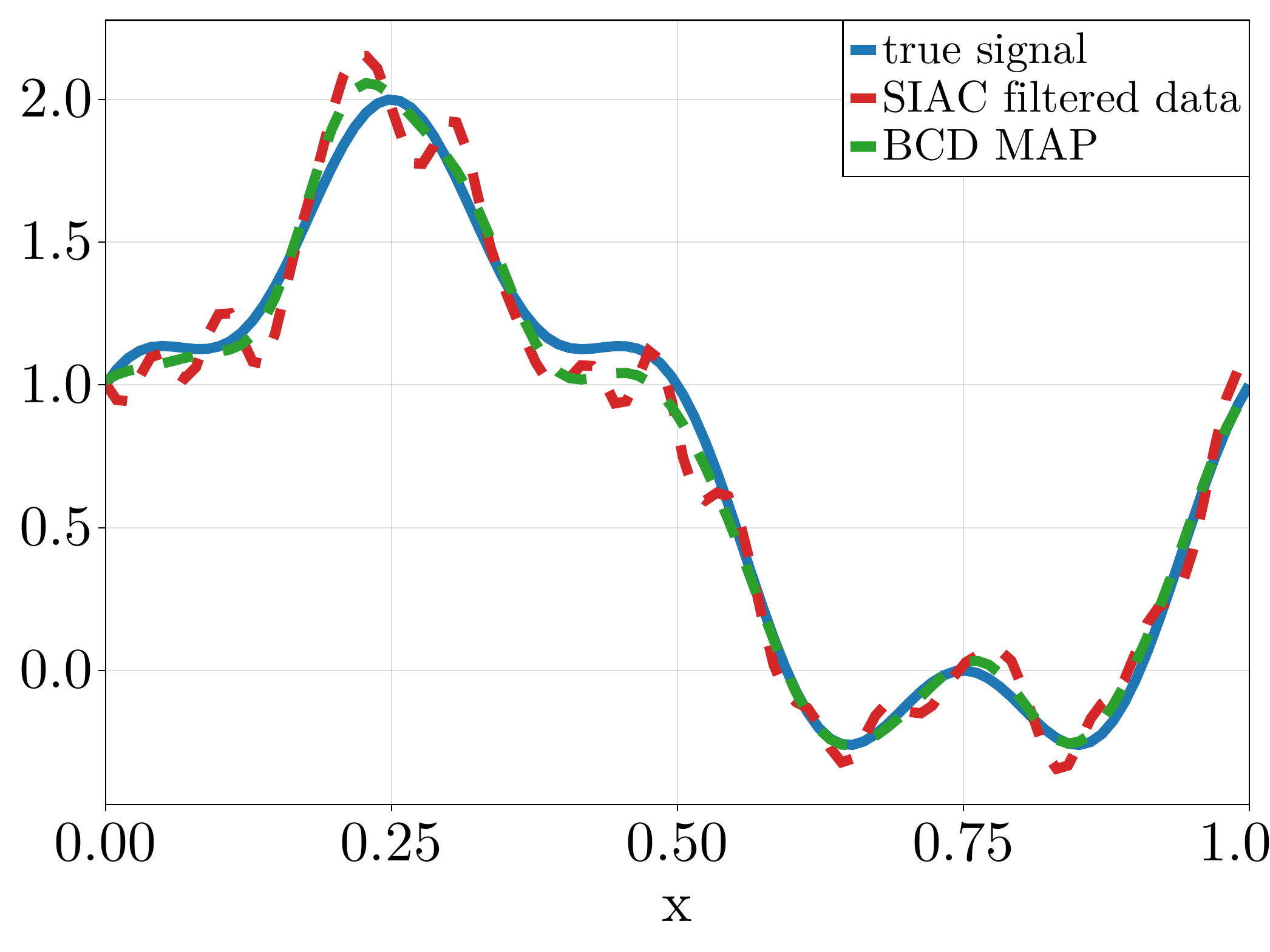}
    \caption{$k=2$}
    \label{fig:differentKernels_MAP_k2}
    \end{subfigure}%
    \\ 
    \begin{subfigure}[b]{0.49\textwidth}
		\includegraphics[width=\textwidth]{%
			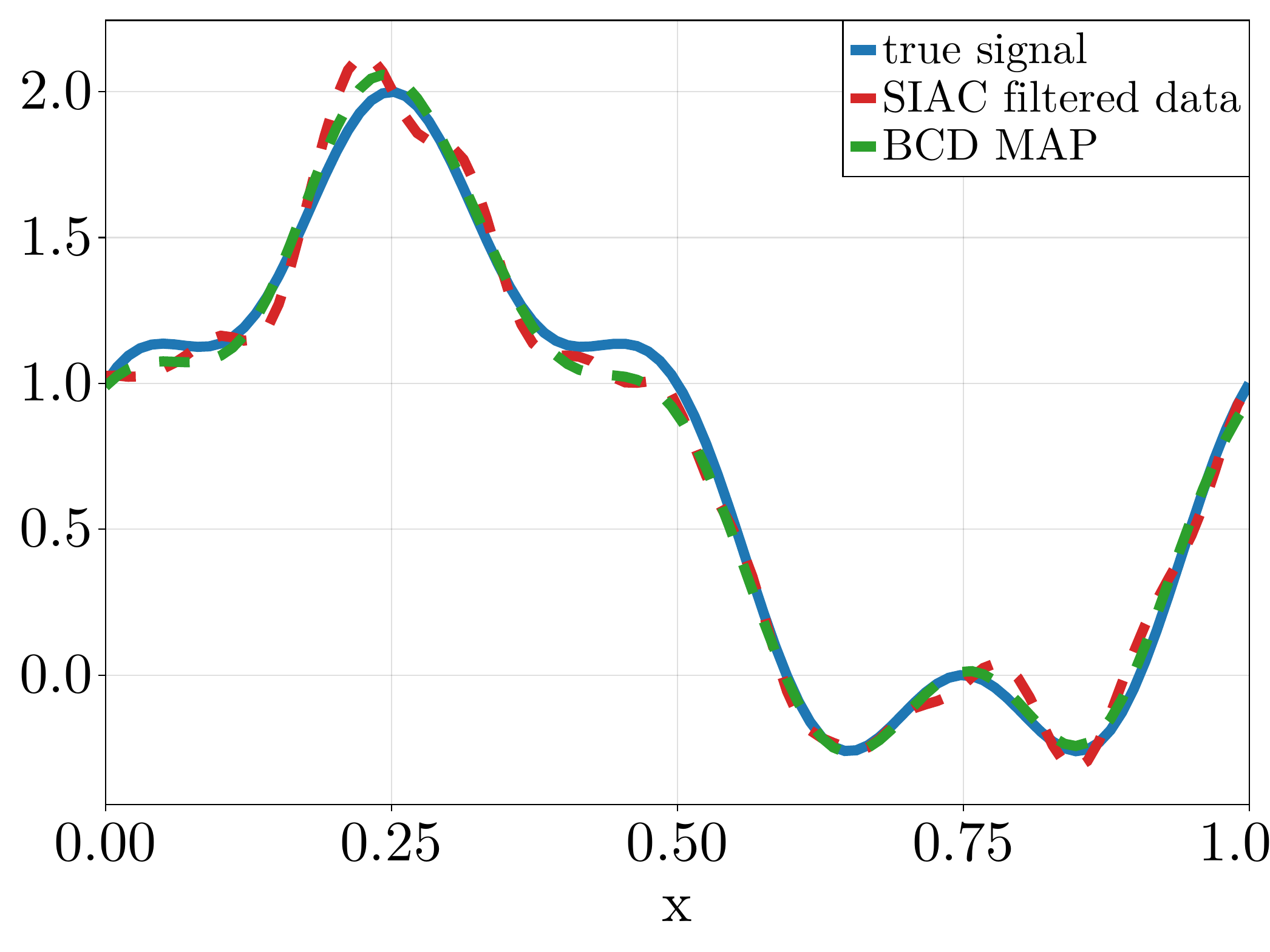}
    \caption{$k=3$}
    \label{fig:differentKernels_MAP_k3}
    \end{subfigure}%
    ~ 
    \begin{subfigure}[b]{0.49\textwidth}
		\includegraphics[width=\textwidth]{%
			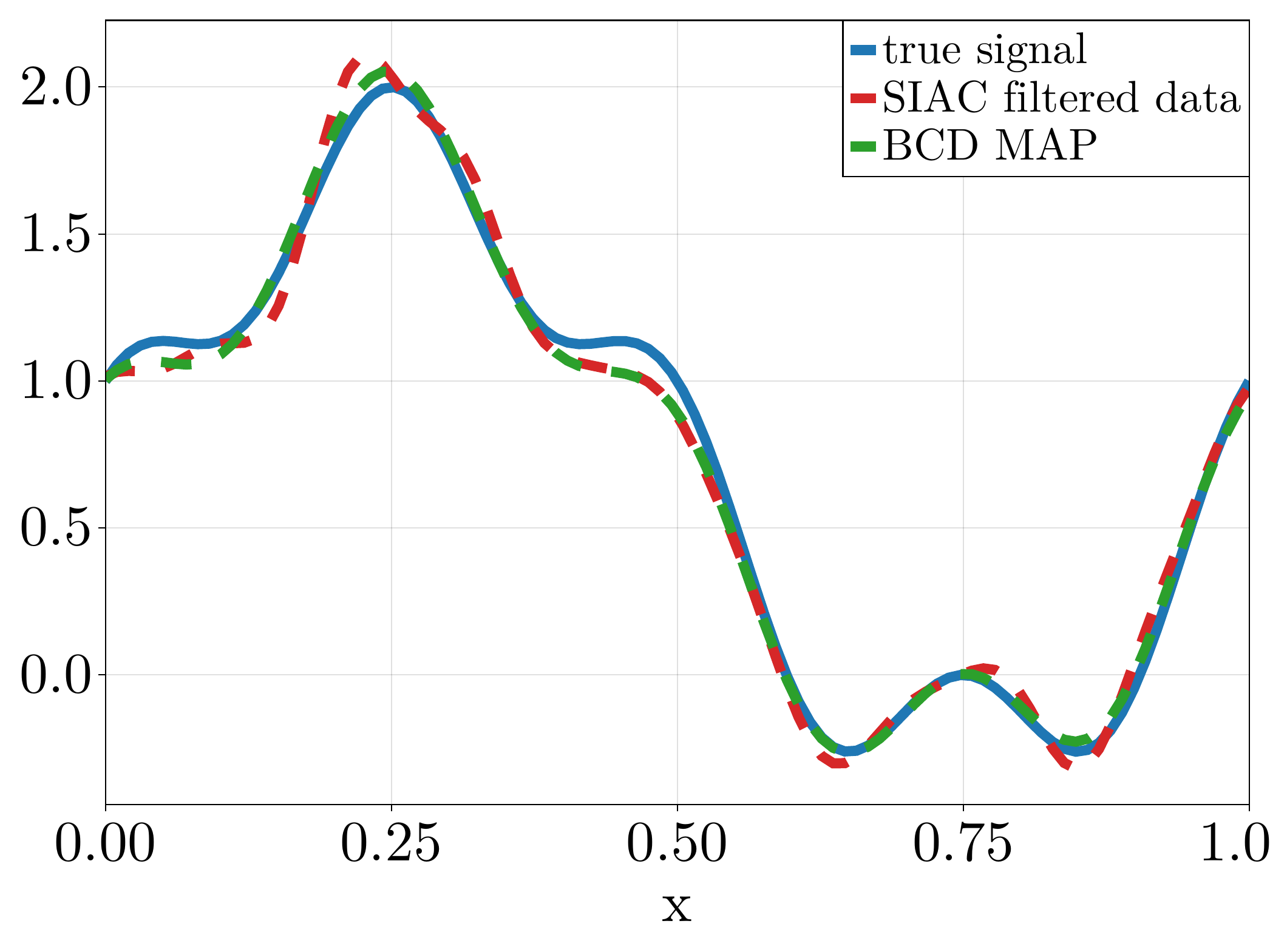}
    \caption{$k=4$}
    \label{fig:differentKernels_MAP_k4}
    \end{subfigure}%
    \caption{
        The true signal $u$, the SIAC-filtered data $\mathbf{F} \mathbf{b}$, and the MAP estimate of the Bayesian SIAC filter for SIAC filters $\mathcal{F}$ with different kernels $K^{(2k+1,k+1)}$, where $k$ is the degree
	}
    \label{fig:differentKernels_MAP}
\end{figure}

\begin{figure}[tb]
	\centering 
	\begin{subfigure}[b]{0.49\textwidth}
		\includegraphics[width=\textwidth]{%
			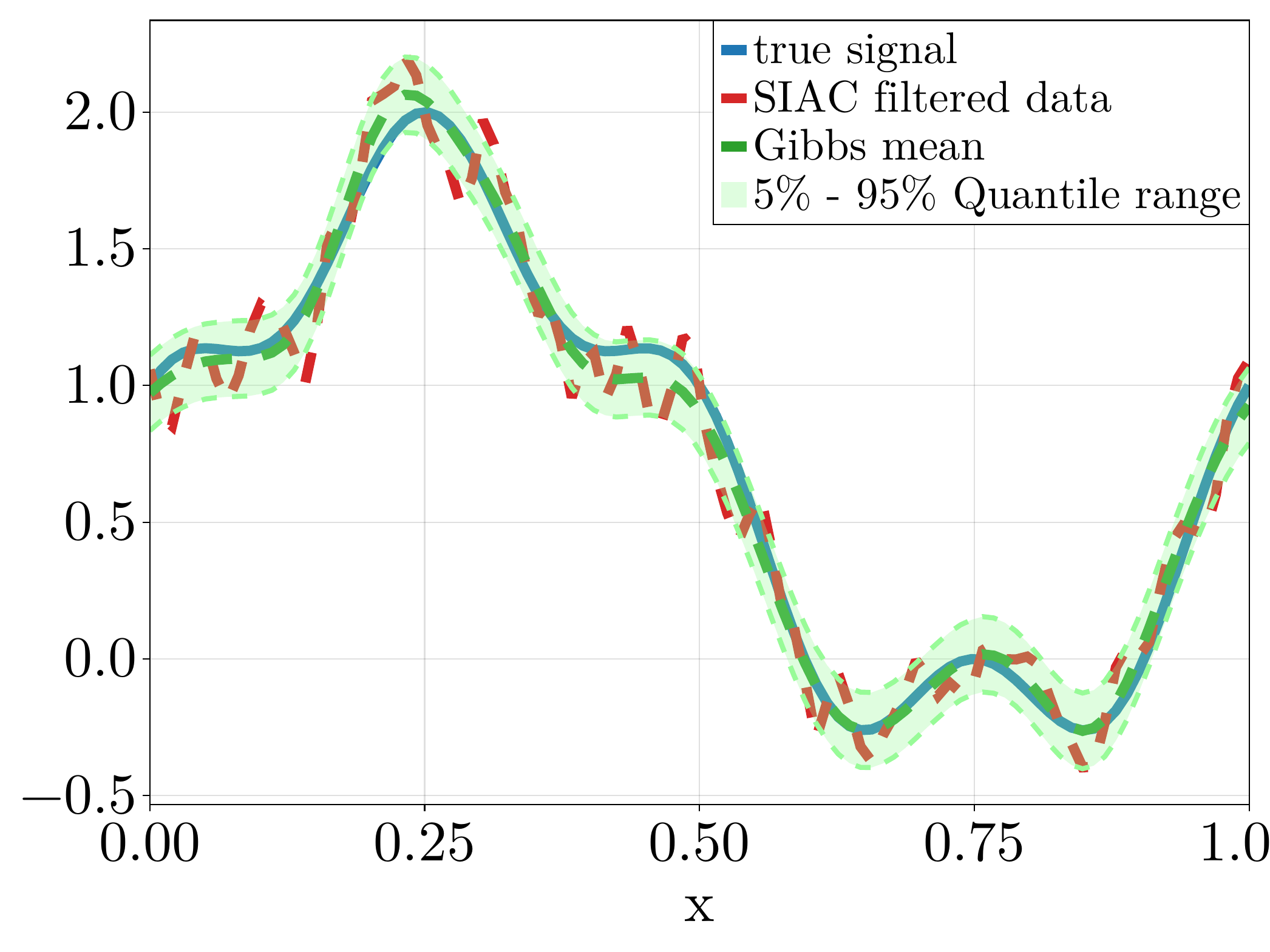}
    \caption{$k=1$}
    \label{fig:differentKernels_Gibbs_k1}
    \end{subfigure}%
	~
	\begin{subfigure}[b]{0.49\textwidth}
		\includegraphics[width=\textwidth]{%
			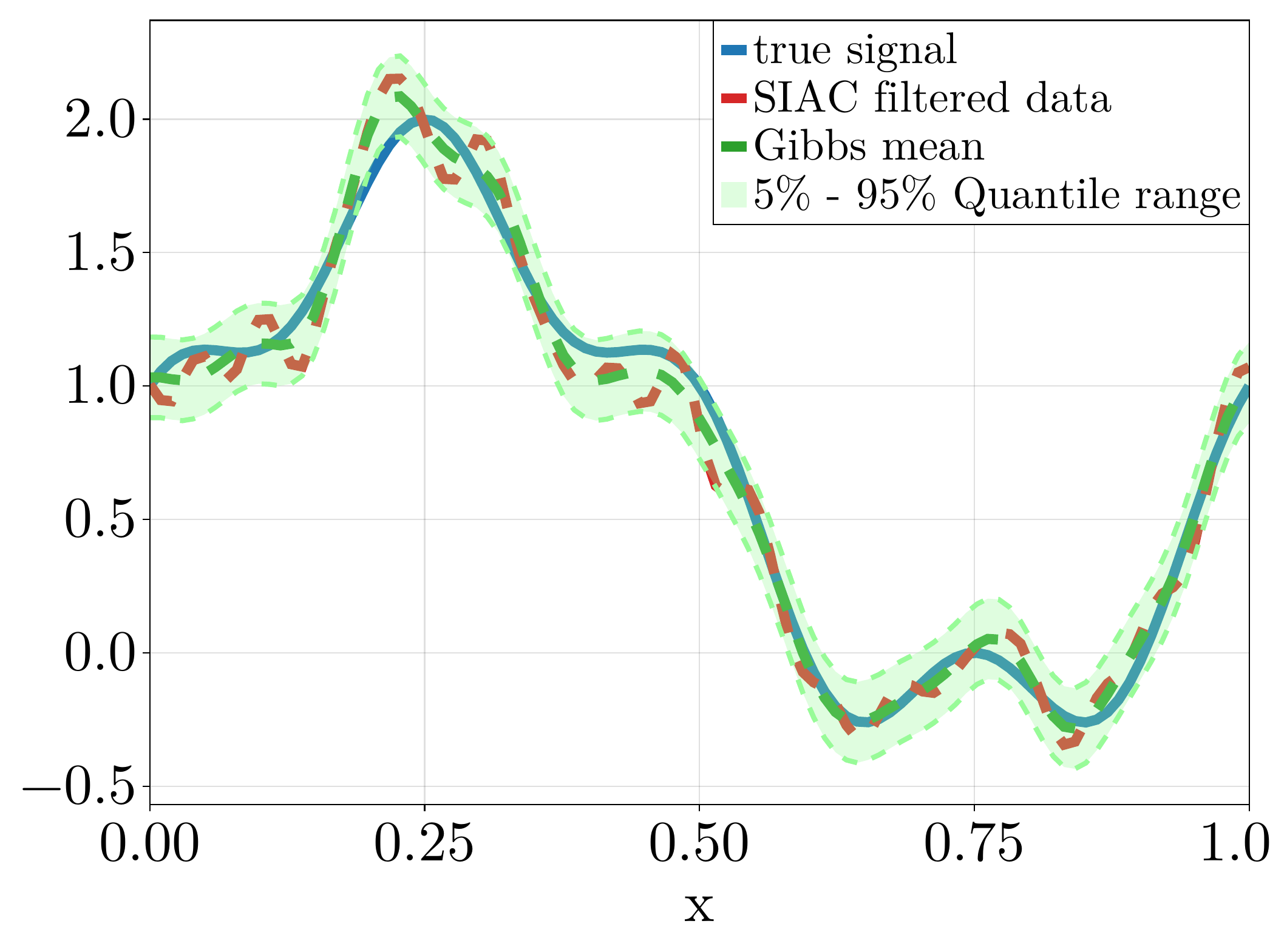}
    \caption{$k=2$}
    \label{fig:differentKernels_Gibbs_k2}
    \end{subfigure}%
    \\ 
    \begin{subfigure}[b]{0.49\textwidth}
		\includegraphics[width=\textwidth]{%
			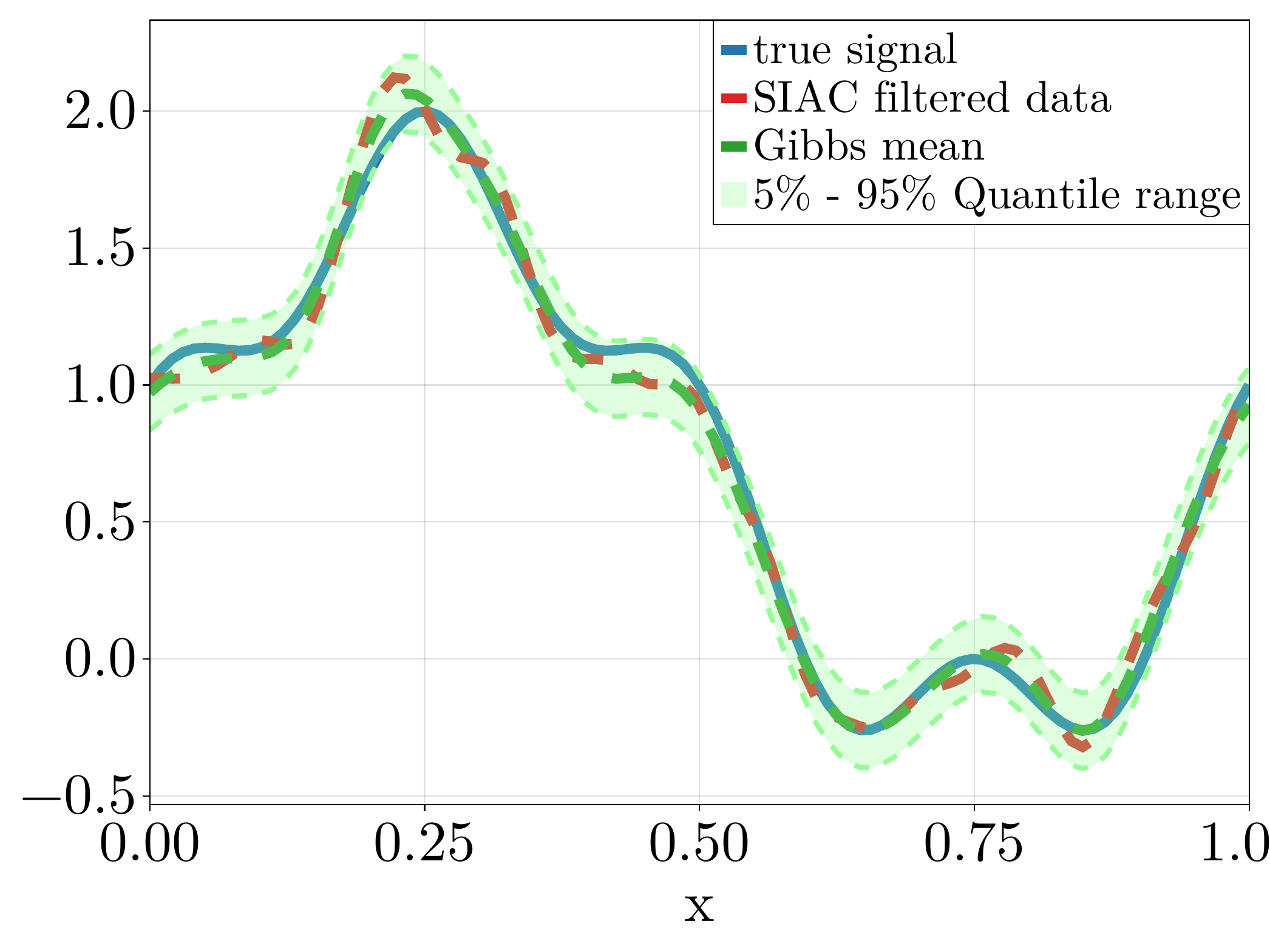}
    \caption{$k=3$}
    \label{fig:differentKernels_Gibbs_k3}
    \end{subfigure}%
    ~ 
    \begin{subfigure}[b]{0.49\textwidth}
		\includegraphics[width=\textwidth]{%
			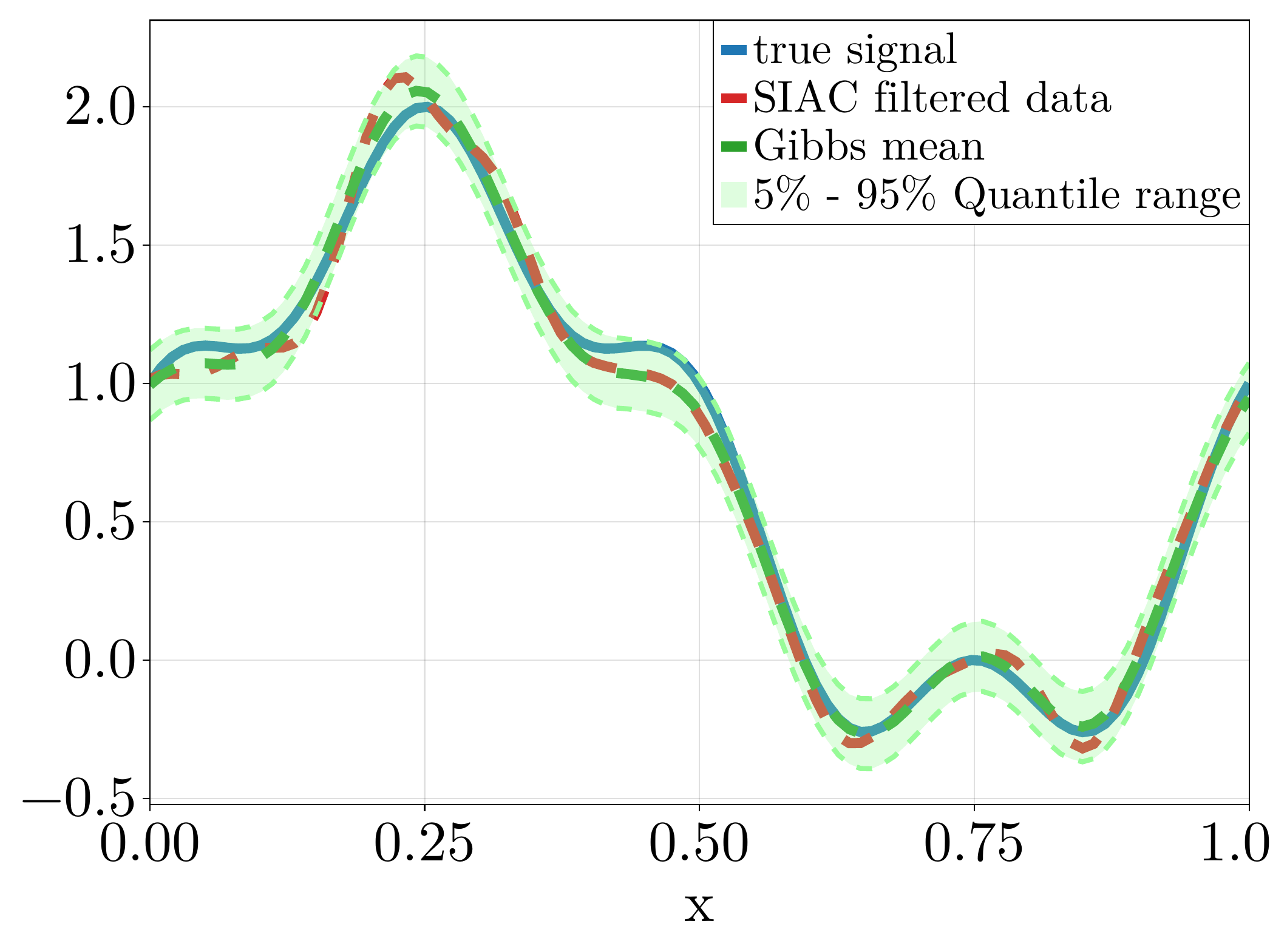}
    \caption{$k=4$}
    \label{fig:differentKernels_Gibbs_k4}
    \end{subfigure}%
    \caption{
        The true signal $u$, the SIAC-filtered data $\mathbf{F} \mathbf{b}$, and the sample mean of the Bayesian SIAC filter for SIAC filters $\mathcal{F}$ with different kernels $K^{(2k+1,k+1)}$, where $k$ is the degree
	}
    \label{fig:differentKernels_Gibbs}
\end{figure}

\cref{fig:differentKernels_data} illustrates the true signal $u$, the observational data $\mathbf{b}$, and the SIAC-filtered data $\mathbf{F} \mathbf{b}$ for different kernels with $k=1,2,3,4$. 
We observe that the choice of kernel in the SIAC filter visibly affects the quality of the deterministic SIAC filter, particularly for small values of the degree $k$. 
At the same time, \cref{fig:differentKernels_MAP} visualizes the true signal $u$, the SIAC-filtered data $\mathbf{F} \mathbf{b}$, and the MAP estimate of the Bayesian SIAC filter for the same kernels $K^{(2k+1,k+1)}$ with degrees $k=1,2,3,4$. 
\cref{fig:differentKernels_MAP} demonstrates that, while the MAP estimates of the Bayesian SIAC filter are also influenced by the choice of the kernel, the reconstructions of the Bayesian SIAC filter are both more robust with respect to the choice of the kernel and more accurate in all cases. 
\cref{fig:differentKernels_Gibbs} illustrates the same advantage of the Bayesian SIAC filter over the deterministic one for the sample mean of the Bayesian SIAC filter. 
The sample mean has been computed from $10^5$ samples generated by the Gibbs sampler, as detailed in \Cref{sub:numerics_signal}.
Furthermore, we can observe an increase in uncertainty in \cref{fig:differentKernels_Gibbs} for the kernels with smaller degrees $k$---for instance, compare \cref{fig:differentKernels_Gibbs_k1} and \cref{fig:differentKernels_Gibbs_k4}. 
This highlights an additional advantage of the Bayesian SIAC filter: the quantified uncertainty in the reconstructed signal can indicate misspecified (i.e., poorly chosen) kernels. 

\begin{table}[tb]
    \centering
    \resizebox{\textwidth}{!}{
    \begin{tabular}{l | c c c | c c c | c c c}
        \toprule
         & \multicolumn{3}{c}{relative $\ell^1$-errors} 
            & \multicolumn{3}{c}{relative $\ell^2$-errors} 
            & \multicolumn{3}{c}{relative $\ell^{\infty}$-errors} \\ 
        \midrule
        $k$ & det. & MAP & mean
            & det. & MAP & mean
            & det. & MAP & mean \\
        \midrule
        1 & 1.3e-01 & \textbf{8.6e-02} & 9.0e-02 
            & 1.2e-01 & \textbf{8.5e-02} & 8.9e-02 
            & 1.4e-01 & \textbf{1.2e-01} & \textbf{1.2e-01} \\
        2 & 1.1e-01 & \textbf{5.1e-02} & 6.5e-02 
            & 1.0e-01 & \textbf{5.4e-02} & 6.7e-02 
            & 1.3e-01 & \textbf{7.2e-02} & 9.8e-02 \\
        3 & 7.1e-02 & \textbf{4.7e-02} & 5.1e-02 
            & 7.2e-02 & \textbf{4.8e-02} & 5.2e-02 
            & 1.0e-01 & \textbf{5.8e-02} & 6.0e-02 \\
        4 & 6.9e-02 & \textbf{4.6e-02} & 4.7e-02 
            & 6.8e-02 & \textbf{5.8e-02} & 5.9e-02 
            & 9.3e-02 & \textbf{5.8e-02} & 5.9e-02 \\
        5 & 6.5e-02 & \textbf{4.5e-02} & \textbf{4.5e-02} 
            & 6.1e-02 & \textbf{4.7e-02} & \textbf{4.7e-02} 
            & 6.2e-02 & \textbf{5.6e-02} & 5.8e-02 \\
        \bottomrule
    \end{tabular}
    }
    \caption{
    Relative $\ell^1$-, $\ell^2$-, and $\ell^{\infty}$-errors of the deterministic SIAC filter (``det.'') and the MAP estimate and sample mean of the Bayesian SIAC filter (``MAP'' and ``mean'') for SIAC filters $\mathcal{F}$ with different kernels $K^{(2k+1,k+1)}$, where $k$ is the degree 
    }
    \label{tab:differentKernels_errors}
\end{table}

Finally, \cref{tab:differentKernels_errors} provides the relative $\ell^1$-, $\ell^2$-, and $\ell^{\infty}$-errors of the deterministic SIAC filter and the MAP estimate and sample mean of the Bayesian SIAC filter for different kernels $K^{(2k+1,k+1)}$ with degrees $k=1,2,3,4,5$.
We observe that the MAP estimate of the Bayesian SIAC filter yields the smallest errors for all norms and degrees $k$. 
Furthermore, the sample mean of the Bayesian SIAC filter is only slightly less accurate. 
On the other hand, the deterministic SIAC filter yields notably larger errors---particularly for $k=1$ and $k=2$.
We also observe that the difference between the errors for smaller and larger degrees $k$ is less pronounced for the Bayesian SIAC filter than for the deterministic one. 
This observation reaffirms that the Bayesian SIAC filter is more robust than the deterministic SIAC filter with respect to the kernel.

\subsection{Post-processing for DG numerical solutions} 
\label{sub:numerics_DG}

As mentioned earlier, the deterministic SIAC filter was initially introduced as a post-processing technique for numerical PDE solutions \cite{bramble1977higher,cockburn2003enhanced}, and it has most prominently been developed in the context of DG methods and hyperbolic conservation laws \cite{ryan2005extension,curtis2008postprocessing,mirzaee2010quantification,mirzaee2012efficient,ji2014superconvergent}.
We thus next consider the post-processing of numerical DG solutions of the linear advection equation 
\begin{equation}\label{eq:linear_adv}
    \partial_t u(x,t) + \partial_x u(x,t) = 0, \quad 
    x \in (0,1), \ t \in (0,T),
\end{equation}
with periodic boundary conditions and smooth initial condition $u(x,0) = \sin(\pi x)/2 + 1$.
We discretize \cref{eq:linear_adv} in space using a nodal DG method with $J$ uniform cells and $k+1$ Gauss--Legendre points per element---resulting in a piecewise polynomial approximation with degree $k$. 
Furthermore, we use a full-upwind numerical flux to introduce coupling between neighboring cells. 
We refer to \cite{hesthaven2008nodal} for more details on nodal DG methods. 
To integrate the numerical solution in time, we use the ninth-order Verner Runge--Kutta scheme \cite{verner2010numerically}. 

\begin{figure}[tb]
	\centering 
	\begin{subfigure}[b]{0.49\textwidth}
		\includegraphics[width=\textwidth]{%
			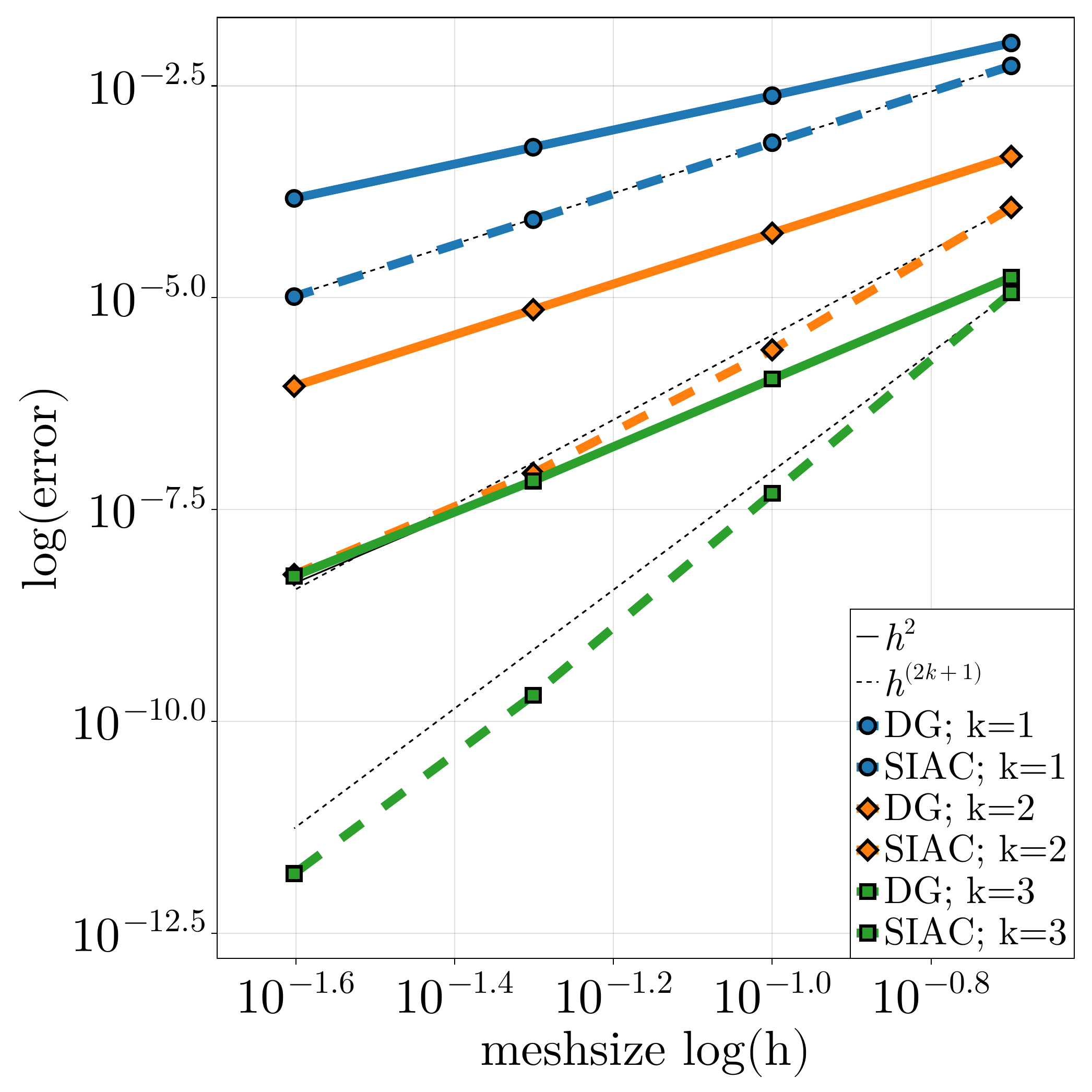} 
        \caption{Deterministic SIAC filter}
        \label{fig:DG_convergence_det}
    \end{subfigure}%
	\begin{subfigure}[b]{0.49\textwidth}
		\includegraphics[width=\textwidth]{%
			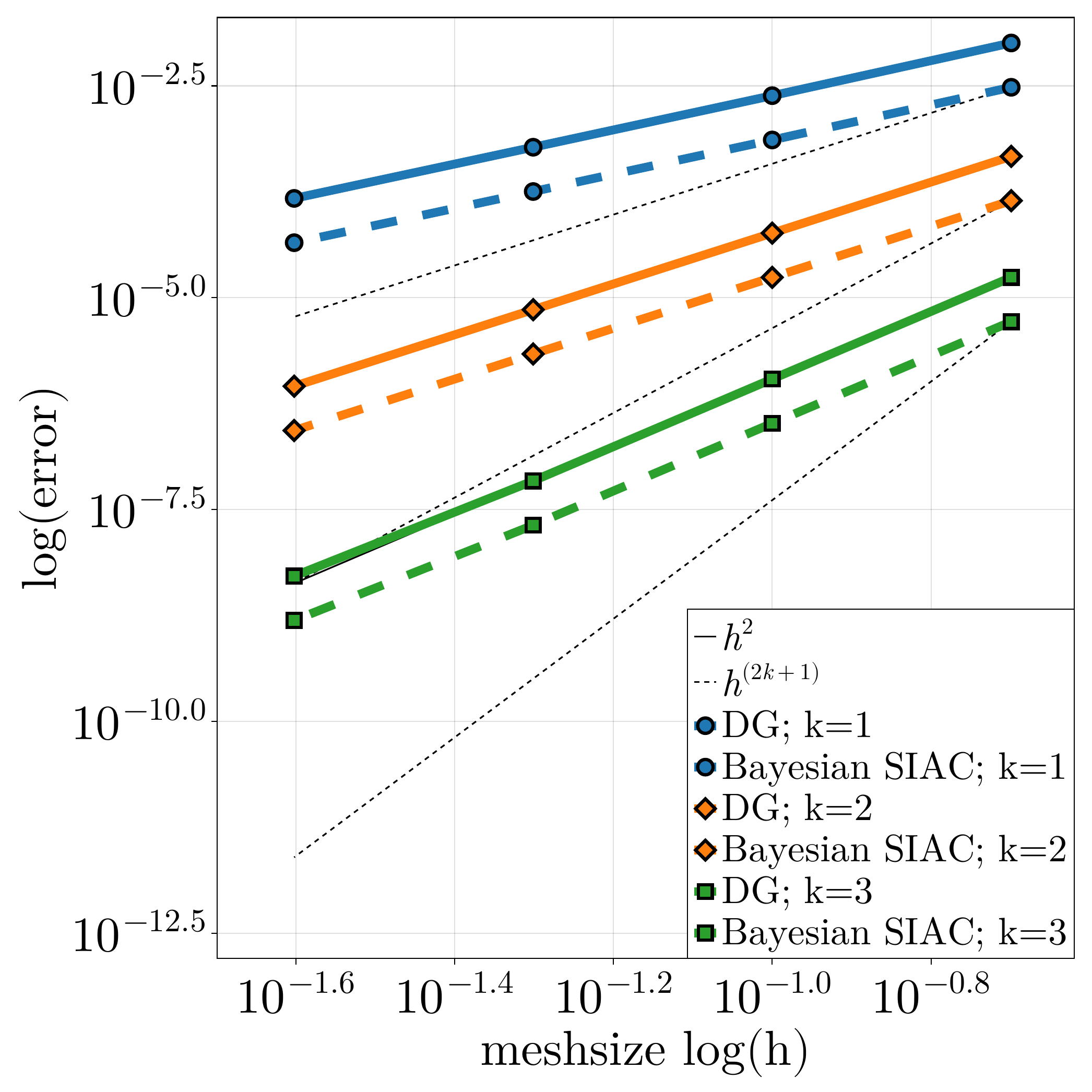} 
        \caption{Bayesian SIAC filter}
        \label{fig:DG_convergence_Bay}
    \end{subfigure}
    \caption{
        $L^2$-errors of the numerical DG solutions and their post-processed counterparts for the linear advection equation \cref{eq:linear_adv} at time $t = 1$, using polynomial degrees $k = 1, 2, 3$. 
        We compare the post-processed solutions using the deterministic SIAC filter (\cref{fig:DG_convergence_det}) with the MAP estimate of the Bayesian SIAC filter (\cref{fig:DG_convergence_Bay}), efficiently computed using the BCD algorithm. 
    }
    \label{fig:DG_convergence}
\end{figure}

In \cref{fig:DG_convergence}, we compare the $L^2$-errors of the standard deterministic SIAC filter and the Bayesian SIAC filter, applied to well-resolved DG approximations of a smooth solution of the linear advection equation \cref{eq:linear_adv} at $t=1$ for polynomial degrees $k=1,2,3$. 
For the Bayesian SIAC filter, we report the MAP estimates obtained from the efficient BCD algorithm. 
The solid black lines in \cref{fig:DG_convergence} indicate the expected DG convergence rate of $k+1$, while the dashed black lines show the superconvergent rate of $2k+1$ achieved by the deterministic SIAC filter. 
The solid colored lines represent the observed convergence of the unfiltered DG solutions, and the dashed colored lines correspond to the post-processed results from the deterministic SIAC filter (\cref{fig:DG_convergence_det}) and the MAP estimates from the Bayesian SIAC filter (\cref{fig:DG_convergence_Bay}).
The results demonstrate that the deterministic SIAC filter achieves superconvergence, in agreement with previous studies \cite{docampo2020enhancing}. 
In contrast, the Bayesian SIAC filter reduces the error of the DG solutions but does not improve the asymptotic convergence rate. 
We observed similar behavior for the posterior mean estimates of the Bayesian SIAC filter (not shown for brevity). 
These findings indicate that, while the Bayesian SIAC framework provides accurate estimates and enables UQ as well as general data models, it does not replicate the superconvergent behavior of the deterministic SIAC filter. 
Consequently, for post-processing highly resolved DG approximations of smooth solutions of hyperbolic conservation laws, the classical deterministic SIAC filter remains the method of choice.
It remains to be addressed in future works if the same is true for under-resolved and discontinuous solutions.

\begin{remark}
Deterministic SIAC filters are amenable to rigorous analysis and known to exhibit superconvergent behavior under suitable smoothness and scaling assumptions \cite{suarez2017regularization}. 
In the Bayesian SIAC formulation considered here, the probabilistic treatment fundamentally alters the mechanisms, such as exact polynomial reproduction and moment cancellation, that underlie superconvergence in the deterministic setting. 
As a result, superconvergent behavior is generally not retained in the present framework.
One possible avenue toward recovering superconvergence within a Bayesian setting would be to impose additional structural constraints, for example, by anchoring the inference around the deterministic SIAC reconstruction (e.g., via a Laplace approximation) or by restricting the probabilistic treatment to subspaces that do not interfere with the superconvergent modes. 
Similarly, while noise could in principle be mitigated in a deterministic SIAC setting through additional SIAC filtering, such successive filtering modifies the effective reconstruction operator and obscures whether the analytical mechanisms responsible for superconvergence can be preserved. 
The Bayesian SIAC approach instead incorporates uncertainty and noise directly at the model level, yielding a different accuracy–cost trade-off that prioritizes robustness and uncertainty quantification over strict superconvergence.
\end{remark}

\begin{remark}
    Approximate solutions found using Galerkin schemes are known to exhibit oscillatory error components around the true solution. 
    This phenomenon is closely related to the classical interpretation of DG errors in negative-order norms, where oscillatory components cancel and lead to superconvergent behavior; see \cite{cockburn2003enhanced, wahlbin2006superconvergence,adjerid2007discontinuous}.
    These oscillatory error components, which are present even in linear equations, may be interpreted as small-scale perturbations analogous to noise. 
    This point of view motivates applying the Bayesian SIAC filter to linear equations. 
    For nonlinear hyperbolic equations, SIAC filtering can improve the convergence rate in smooth regions away from discontinuities, similarly to the linear case, as demonstrated in \cite{cockburn2003enhanced,li2015siac}. 
    In the proximity of discontinuities, a different type of oscillatory behavior occurs  in the form of the Gibbs phenomenon, degrading the accuracy and stability of the approximation. 
    SIAC filters have been shown to be effective tools for shock regularization in this context \cite{bohm2019multi}. 
    Extending the Bayesian SIAC framework to nonlinear problems with discontinuities is beyond the scope of the present work and is left for future investigation.
\end{remark}

\subsection{Image reconstruction from indirect data}
\label{sub:numerics_image}

\begin{figure}[tb]
	\centering
    \begin{subfigure}[b]{0.49\textwidth}
		\includegraphics[width=\textwidth]{%
            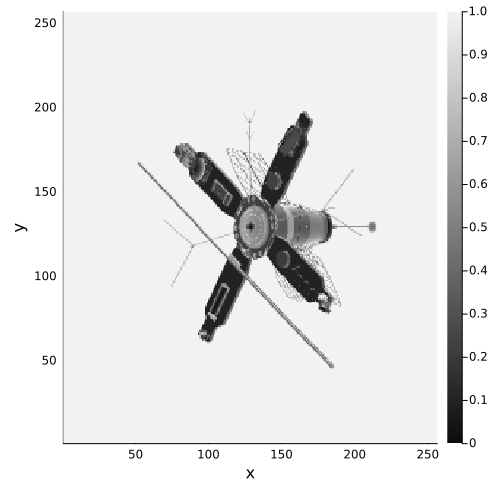} 
    	\caption{Reference image}
    	\label{fig:deconvolution_2d_ref}
    \end{subfigure}%
    ~
    \begin{subfigure}[b]{0.49\textwidth}
		\includegraphics[width=\textwidth]{%
            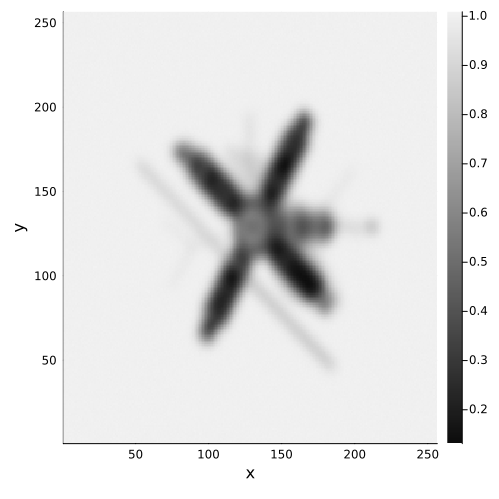} 
    	\caption{Noisy blurred image}
    	\label{fig:deconvolution_2d_blurred}
    \end{subfigure}%
	\\ 
    \begin{subfigure}[b]{0.49\textwidth}
		\includegraphics[width=\textwidth]{%
        		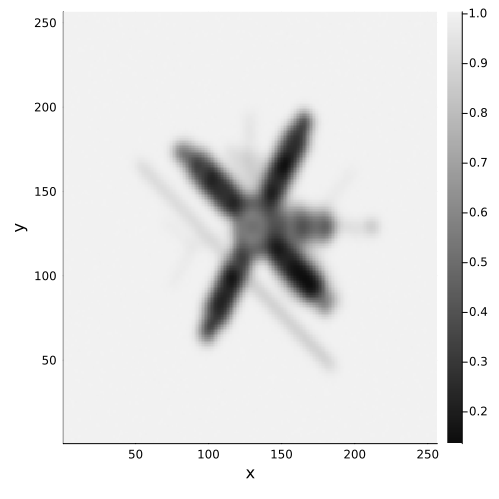} 
    	\caption{Deterministic SIAC filter}
    	\label{fig:deconvolution_2d_detSIACY}
    \end{subfigure}%
    ~
    \begin{subfigure}[b]{0.49\textwidth}
		\includegraphics[width=\textwidth]{%
        		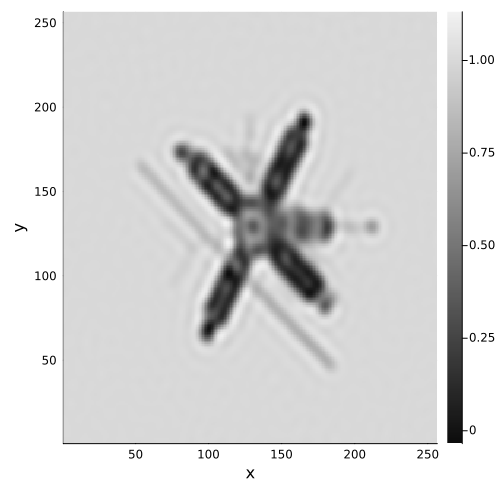} 
    	\caption{Bayesian SIAC filter (MAP estimate)}
    	\label{fig:deconvolution_2d_MAP}
    \end{subfigure}%
    \caption{ 
  	The reference image, the corresponding noisy blurred image, the SIAC-filtered reference image, and the SIAC-filtered noisy blurred image
    }
    \label{fig:deconvolution_2d}
\end{figure}

The deterministic SIAC filter is designed for use with noisy direct observational data, removing spurious oscillations while preserving important higher-frequency features. 
However, the deterministic SIAC filter cannot directly be applied to indirect observational data, such as blurred, Fourier, or Radon data. 
In contrast, the proposed Bayesian SIAC filter can accommodate any observational data straightforwardly, making it broadly applicable. 
To demonstrate this feature, we apply the Bayesian SIAC filter approach to recover the $256 \times 256$ reference image shown in \cref{fig:deconvolution_2d_ref} from its noisy and blurred counterpart $B$ depicted in \cref{fig:deconvolution_2d_blurred}. 
\cref{fig:deconvolution_2d_detSIACY} shows the result of directly applying the deterministic SIAC filter to the noisy, blurred observation. 
We have used the B-spline filter $K^{(3,2)}$ with $k=1$ as described in \cref{sec:SIAC-imple}.
As expected, the output fails to improve the reconstruction, demonstrating that the deterministic SIAC filter is ineffective when applied to indirect data. 

We now provide a brief description of how the Bayesian SIAC filter can be applied to image reconstruction, before discussing its result in \cref{fig:deconvolution_2d_MAP}. 
The observational data image $\mathbf{B}$ results from the reference image $\mathbf{U}$ by applying a discrete one-dimensional convolution operator in the two canonical coordinate directions and then adding i.i.d.\ zero-mean normal noise with variance $\sigma^2 = 10^{-5}$. 
We note that for unstructured or non-uniform data, extensions such as the one-dimensional Line SIAC \cite{LSIAC} and related approaches \cite{jallepalli2024non} can be employed; however, a mesh intersection algorithm would need to be implemented in order to apply the filter to non-uniform data.
The corresponding forward model is $\mathbf{B} = \mathbf{A}_{\rm 1d} \mathbf{U} \mathbf{A}_{\rm 1d}^T + \mathbf{E}$ or, equivalently, 
\begin{equation}
    \mathbf{b} = \mathbf{A} \mathbf{u} + \mathbf{e}, 
\end{equation} 
after vectorization. 
Here, $\mathbf{z} = \Vec(\mathbf{Z})$ denotes the $mn\times1$ column vectors obtained by stacking the columns of the $m \times n$ matrix $\mathbf{Z}$ on top of one another, and $\mathbf{A} = \mathbf{A}_{\rm 1d} \otimes \mathbf{A}_{\rm 1d}$ denotes the Kronecker product, which constructs a 2D operator from the 1D operator $\mathbf{A}_{\rm 1d}$ applied in both directions. 
Furthermore, $\mathbf{A}_{\rm 1d}$ is obtained by applying the midpoint quadrature to the convolution equation $b(x) = \int \kappa(x-s) u(s) \intd s$. 
We use a Gaussian convolution kernel of the form $\kappa(s) = \exp( - s^2 / (2 \gamma^2) ) / \sqrt{2 \pi \gamma^2}$ with blurring parameter $\gamma = 1.5 \cdot 10^{-2}$. 
The one-dimensional forward operator thus is 
\begin{equation}
	[\mathbf{A}_{\rm 1d}]_{ij} = h \kappa( h[i-j] ), \quad i,j=1,\dots,n,
\end{equation}
where $h=1/n$ is the distance between consecutive grid points on the same vertical/horizontal line. 
Note that $\mathbf{A}_{\rm 1d}$ has full rank but quickly becomes ill-conditioned.

\cref{fig:deconvolution_2d_MAP} shows the MAP estimate of the Bayesian SIAC filter posterior, efficiently computed using the BCD approach described in \Cref{sec:MAP}. 
The reconstruction in \cref{fig:deconvolution_2d_MAP} from the noisy, blurred data in \cref{fig:deconvolution_2d_blurred} demonstrates the effectiveness and flexibility of the Bayesian SIAC filter in handling general data models. 
While we focus on the MAP estimate here, the Bayesian formulation in principle supports full UQ via posterior sampling. 
However, in this setting, the posterior is $65{,}536$-dimensional ($256 \times 256$ pixel values plus two hyper-parameters), and accurate UQ would require generating a prohibitively large number of samples. 
This motivates future exploration of sampling-efficient strategies such as low-rank surrogate models to make full UQ more tractable in such large-scale settings.

%% file: 7_summary.tex
\section{Summary} 
\label{sec:summary} 

We introduced the Bayesian SIAC filter, a hierarchical Bayesian generalization of the deterministic SIAC filter.
The proposed framework supports general (non-nodal) data models and enables rigorous UQ, thereby broadening the applicability of SIAC filtering beyond its traditional scope.
We also developed structure-exploiting algorithms for efficient MAP estimation and MCMC sampling for linear data models with additive Gaussian noise.
Computational experiments---including signal denoising, post-processing of numerical DG solutions to the linear advection equation, and image deblurring---demonstrated that the Bayesian SIAC filter can deliver point estimates with accuracy comparable to, and in some cases exceeding, that of the deterministic SIAC filter. 
At the same time, our results indicate that for post-processing well-resolved DG approximations of smooth solutions to the linear advection equation, the deterministic SIAC filter remains the method of choice. 
Whether this conclusion also extends to under-resolved or discontinuous solutions remains an open question for future work. 
Nevertheless, the Bayesian SIAC filter substantially broadens the scope of the traditional SIAC methodology: 
It applies to general (non-nodal) data models and, critically, provides rigorous UQ. 
We anticipate that this framework will open the door to applying SIAC techniques to a much wider range of problems.

Future work will focus on extending the Bayesian SIAC filter to nonlinear data models and to piecewise smooth functions featuring jump discontinuities and steep gradients.
While the Bayesian formulation generalizes naturally to non-linear data models---at least formally---via appropriate changes in the likelihood function, achieving efficient MAP estimation and MCMC sampling will require the development of new computational strategies.

%% file: acknowledgements.tex
JG acknowledges support by the DOD (ONR MURI) \#N00014-20-1-2595, the Swedish Research Council (VR) Starting Grant \#2025-05370, the Zenith Career Development Grant \#26.07, and the National Academic Infrastructure for Supercomputing in Sweden (NAISS) grants \#2025/22-1599 and \#2024/22-1207.
TL was supported by DOD ONR MURI grant \#N00014-20-1-2595 and AFOSR grant \#F9550-22-1-0411.
JR and RS were supported by the Swedish Research Council (VR grant 2022-03528). 

Furthermore, JG, TL, and JR acknowledge support from the National Science Foundation (NSF) under Grant \#DMS-1929284 while in residence at the Institute for Computational and Experimental Research in Mathematics (ICERM) in Providence, RI, USA, during the ``Empowering a Diverse Computational Mathematics Research Community" topical workshop.

%% file: appA_notation.tex
\section{Notation} 
\label{app:notation} 

We provide some notation and meanings in \cref{tab:notation}.

\begin{table}[h]
	\centering
	\begin{tabular}{@{}c  c  c@{}}
	\toprule
	symbol & space & meaning
	\\\midrule
    $u$ & $\R^d \to \R$ & underlying continuous function of interest \\
	$\mathbf{u}$ & $\R^N$ & unknown parameter vector \\
    $\mathbfcal{A}$ & $\R^N \to \R^M$ & forward operator \\
    $\mathbf{A}$ & $\R^{M \times N}$ & matrix representation of linear forward operator \\
	$\mathbf{b}$ & $\R^M$ & observational data vector \\
	$\mathbf{e}$ & $\R^M$ & observational noise vector \\
    $\sigma^2$ & $\R_{>0}$ & true noise variance \\
    $\mathcal{F}$ & $u \mapsto \mathcal{F}[u]$ & SIAC filter \\
    $\mathbf{F}$ & $\R^{N \times N}$ & matrix representation of SIAC filter \\
    $\mathbf{I}$ & $\R^{N \times N}$ & identity matrix \\
    $\alpha$ & $\R_{>0}$ & noise precision parameter \\
    $\beta$ & $\R_{>0}$ & prior precision parameter \\
    $\boldsymbol{\theta}$ & $\R_{>0}^2$ & collection of hyper-parameters, i.e., $\boldsymbol{\theta} = [\alpha, \beta]$ \\
    $c_{\alpha}$, $d_{\alpha}$ & $\R_{>0}$, $\R_{>0}$ & parameters of the gamma hyper-prior for $\alpha$ \\
    $c_{\beta}$, $d_{\beta}$ & $\R_{>0}$, $\R_{>0}$ & parameters of the gamma hyper-prior for $\beta$ \\ 
    $k$ & $\N$ & degree of the piecewise polynomial approximation of $u$ \\ 
    $J$ & $\N$ & \#cells for the piecewise polynomial approximation of $u$ \\
	\bottomrule
	\end{tabular} 
    \caption{Table of notation.}
	\label{tab:notation}
\end{table}

%% file: appB_SIACimpl.tex
\section{Implementation of the SIAC filter}
\label{sec:SIAC-imple}

We describe how to obtain a matrix representation $\mathbf{F} \in \R^{N \times N}$ of the SIAC filter $\mathcal{F} \equiv K^{(2k+1,k+1)}$ so that $\mathbf{F} \mathbf{u} = \mathcal{F}[u](\mathbf{x})$, assuming that $u$ is given---or approximated---by a piecewise polynomial of degree $k$ on a periodic mesh consisting of $J$ cells. 
Here $\mathbf{x} = [x_1,\dots,x_N]^T \in \R^N$ are grid points, $\mathbf{u} = [u_1,\dots,u_N]^T \in \R^N$ with $u_n = u(x_n)$ are nodal values of $u$ at the grid points, and $\mathcal{F}[u](\mathbf{x}) = [\mathcal{F}[u](x_1),\dots,\mathcal{F}[u](x_N)]^T \in \R^N$ are the function values of the SIAC-filtered function $\mathcal{F}[u]$. 

 We enumerate the grid points and nodal values of $u$ according to the mesh as $\mathbf{x} = [x_1^1,\dots,x_{k+1}^1,\dots,x_1^J,\dots,x_{k+1}^J]^T$ and $\mathbf{u} = [u_1^1,\dots,u_{k+1}^1,\dots,u_1^J,\dots,u_{k+1}^J]^T$, where $N = J(k+1)$, $x_p^j$ is the $p$th local grid point in the $j$th cell, and $u_p^j$ is the corresponding function value of $u$.  
Next, let $\{ \phi_p^j \}_{p=1,j=1}^{k+1,J}$ be a basis of local Lagrange polynomials that satisfy the cardinal property $\phi_p^j(x_q^i) = 1$ if $(p,j) = (q,i)$ and $\phi_p^j(x_q^i) = 0$ otherwise. 
We can then write $u$ as $u(x) = \sum_{j=1}^J\, \sum_{p=1}^{k+1} u_p^j \phi_p^j(x)$. 
Moreover, the SIAC-filtered function $\mathcal{F}[u]$ can then be expressed as
\begin{equation}\label{eq:repr_Fu}
    \mathcal{F}[u](x) = \sum_{j=1}^J\, \sum_{p=1}^{k+1} u_p^j \mathcal{F}[\phi_p^j](x).
\end{equation}
If we now evaluate $\mathcal{F}[u]$ at the grid points $\mathbf{x}$, \cref{eq:repr_Fu} implies 
\begin{equation}\label{eq:matrix_F}
    \mathcal{F}[u](\mathbf{x}) = \mathbf{F} \mathbf{u}.
\end{equation}
Here, $\mathbf{F} \in \R^{N \times N}$ is the desired matrix representation of the SIAC filter, and its $n$-th row is $\mathbf{F}[n,:] = [\mathcal{F}[\phi_1^1](x_n),\dots,\mathcal{F}[\phi_{k+1}^1](x_n),\dots,\mathcal{F}[\phi_1^J](x_n),\dots,\mathcal{F}[\phi_{k+1}^J](x_n)]$, so that $\mathcal{F}[u](x_n) = \mathbf{F}[n,:] \mathbf{u}$.

Recall that we formulate the SIAC filter using the $\ell$-th order B-splines $B^{(\ell)}$ discussed in \Cref{sec:prelim}, in which case each segment of the SIAC-filtered function $\mathcal{F}[u]$ is of degree $k+\ell$. 
In our implementation, we choose $\ell = k$, resulting in $\mathcal{F}[u]$ being of degree $2k$. 
Furthermore, the computation of $\mathbf{F}$ becomes particularly efficient if the cells are uniform and the same local grid points are used in each element. 
In this case, $\mathcal{F}[\phi_p^j]$ will be the same in different elements. 
For further details, we refer to the publicly available Julia code we used to generate our computational experiments, including the construction of $\mathbf{F}$.

%% file: literature.bib
@misc{glaubitz2025BayesianRepro,
  title={Reproducibility repository for 
    ``{T}he {B}ayesian {SIAC} filter''},
  author={Glaubitz, Jan and Li, Tongtong and Ryan, Jennifer and
          Stuhlmacher, Roman},
  year={2025},
  howpublished={\url{https://github.com/RomanStuhlmacher/paper-2025-Bayesian-SIAC-Filter}}, 
  doi={10.5281/zenodo.19495438}
}

@article{bramble1977higher,
  title={Higher order local accuracy by averaging in the finite element method},
  author={Bramble, James H and Schatz, Alfred H},
  journal={Mathematics of Computation},
  volume={31},
  number={137},
  pages={94--111},
  year={1977}
}

@article{mock1978computation,
  title={Computation of discontinuous solutions of linear hyperbolic equations},
  author={Mock, Michael S and Lax, Peter D},
  journal={Communications on Pure and Applied Mathematics},
  volume={31},
  year={1978}
}

@article{cockburn2003enhanced,
  title={Enhanced accuracy by post-processing for finite element methods for hyperbolic equations},
  author={Cockburn, Bernardo and Luskin, Mitchell and Shu, Chi-Wang and S{\"u}li, Endre},
  journal={Mathematics of Computation},
  volume={72},
  number={242},
  pages={577--606},
  year={2003}
}

@article{EdohPicklo,
author = {Picklo, Matthew and Edoh, Ayaboe},
year = {2025},
month = {06},
pages = {},
title = {Entropy Correction with {SIAC} Filters for High-Order {DG} Methods},
volume = {104},
journal = {Journal of Scientific Computing}
}

@article{SIAC-MRA,
  title={Capitalizing on superconvergence for more accurate multi-resolution discontinuous {G}alerkin methods},
  author={Ryan, Jennifer K},
  journal={Communications on Applied Mathematics and Computation},
  volume={4},
  number={2},
  pages={417--436},
  year={2022},
  publisher={Springer}
}

@article{galindoolartevm,
  title={Superconvergence and accuracy enhancement of discontinuous {G}alerkin solutions for {V}lasov--{M}axwell equations},
  author={Galindo-Olarte, Andr{\'e}s and Huang, Juntao and Ryan, Jennifer and Cheng, Yingda},
  journal={BIT Numerical Mathematics},
  volume={63},
  number={4},
  pages={52},
  year={2023},
  publisher={Springer}
}

@article{jallepalli2019adaptive,
  title={Adaptive characteristic length for {L-SIAC} filtering of {FEM} data},
  author={Jallepalli, Ashok and Haimes, Robert and Kirby, Robert M},
  journal={Journal of Scientific Computing},
  volume={79},
  number={1},
  pages={542--563},
  year={2019},
  publisher={Springer}
}

@article{jallepalli2024non,
  title={Non-uniform knot ({NUK}) {SIAC} post-processing of flow fields produced through unstructured grid adaptation and optimization},
  author={Jallepalli, Ashok and Galbraith, Marshall and Haimes, Robert and Kirby, Robert M},
  journal={Journal of Computational Physics},
  volume={514},
  pages={113238},
  year={2024},
  publisher={Elsevier}
}

@book{wahlbin2006superconvergence,
  title={Superconvergence in Galerkin finite element methods},
  author={Wahlbin, Lars},
  year={2006},
  publisher={Springer}
  }

@article{king2012smoothness,
  title={Smoothness-increasing accuracy-conserving ({SIAC}) filtering for discontinuous {G}alerkin solutions: improved errors versus higher-order accuracy},
  author={King, James and Mirzaee, Hanieh and Ryan, Jennifer K and Kirby, Robert M},
  journal={Journal of Scientific Computing},
  volume={53},
  number={1},
  pages={129--149},
  year={2012},
  publisher={Springer}
}

@incollection{li2015siac,
  title={{SIAC} filtering for nonlinear hyperbolic equations},
  author={Li, Xiaozhou and Ryan, Jennifer K},
  booktitle={Interdisciplinary Topics in Applied Mathematics, Modeling and Computational Science},
  pages={285--291},
  year={2015},
  publisher={Springer}
}

@article{adjerid2007discontinuous,
  title={The discontinuous Galerkin method for two-dimensional hyperbolic problems. Part I: Superconvergence error analysis},
  author={Adjerid, Slimane and Baccouch, Mahboub},
  journal={Journal of Scientific Computing},
  volume={33},
  number={1},
  pages={75--113},
  year={2007},
  publisher={Springer}
}

@article{picklo,
author = {Picklo, Matthew J. and Ryan, Jennifer K.},
title = {Enhanced Multiresolution Analysis for Multidimensional Data Utilizing Line Filtering Techniques},
journal = {SIAM Journal on Scientific Computing},
volume = {44},
number = {4},
pages = {A2628-A2650},
year = {2022}
}

@article{LSIAC,
    author = {Julia Docampo Sánchez and Jennifer K. Ryan and Masha Mirzargar and Robert M. Kirby},
    title = {Multi-dimensional Filtering: Reducing the dimension through rotation.},
    journal = {SIAM Journal on Scientific Computing},
    year = {2017},
    volume = {39},
    pages = {A2179–A2200}
}

@article{PICKLO2024JCP,
  title={Denoising Particle-In-Cell data via {S}moothness-{I}increasing {A}ccuracy-{C}onserving filters with application to {B}ohm speed computation.},
  author={Picklo, Matthew J and Tang, Qi and Zhang, Yanzeng and Ryan, Jennifer K and Tang, Xian-Zhu},
  journal={Journal of Computational Physics},
  volume={502},
  pages={112790},
  year={2024},
  publisher={Elsevier}
}

@article{ryan2005extension,
  title={Extension of a post processing technique for the discontinuous {G}alerkin method for hyperbolic equations with application to an aeroacoustic problem},
  author={Ryan, Jennifer and Shu, Chi-Wang and Atkins, Harold},
  journal={SIAM Journal on Scientific Computing},
  volume={26},
  number={3},
  pages={821--843},
  year={2005},
  publisher={SIAM}
}

@article{curtis2008postprocessing,
  title={Postprocessing for the discontinuous {G}alerkin method over nonuniform meshes},
  author={Curtis, Sean and Kirby, Robert M and Ryan, Jennifer K and Shu, Chi-Wang},
  journal={SIAM Journal on Scientific Computing},
  volume={30},
  number={1},
  pages={272--289},
  year={2008},
  publisher={SIAM}
}

@article{mirzaee2010quantification,
  title={Quantification of errors introduced in the numerical approximation and implementation of {S}moothness-{I}ncreasing {A}ccuracy-{C}onserving {(SIAC)} filtering of discontinuous {G}alerkin {(DG)} fields},
  author={Mirzaee, Hanieh and Ryan, Jennifer K and Kirby, Robert M},
  journal={Journal of Scientific Computing},
  volume={45},
  number={1},
  pages={447--470},
  year={2010},
  publisher={Springer}
}

@article{van2011position,
  title={Position-dependent {S}moothness-{I}increasing {A}ccuracy-{C}onserving {(SIAC)} filtering for improving discontinuous {G}alerkin solutions},
  author={van Slingerland, Paulien and Ryan, Jennifer K and Vuik, Cornelis},
  journal={SIAM Journal on Scientific Computing},
  volume={33},
  number={2},
  pages={802--825},
  year={2011},
  publisher={SIAM}
}

@article{mirzaee2011smoothness,
  title={{S}moothness-{I}ncreasing {A}ccuracy-{C}onserving {(SIAC)} postprocessing for discontinuous {G}alerkin solutions over structured triangular meshes},
  author={Mirzaee, Hanieh and Ji, Liangyue and Ryan, Jennifer K and Kirby, Robert M},
  journal={SIAM Journal on Numerical Analysis},
  volume={49},
  number={5},
  pages={1899--1920},
  year={2011},
  publisher={SIAM}
}

@article{mirzaee2012efficient,
  title={Efficient implementation of {S}moothness-{I}ncreasing {A}ccuracy-{C}onserving {(SIAC)} filters for discontinuous {G}alerkin solutions},
  author={Mirzaee, Hanieh and Ryan, Jennifer K and Kirby, Robert M},
  journal={Journal of Scientific Computing},
  volume={52},
  number={1},
  pages={85--112},
  year={2012},
  publisher={Springer}
}

@article{ji2014superconvergent,
  title={Superconvergent error estimates for position-dependent {S}moothness-{I}increasing {A}ccuracy-{C}onserving {(SIAC)} post-processing of discontinuous {G}alerkin solutions},
  author={Ji, Liangyue and Van Slingerland, Paulien and Ryan, Jennifer and Vuik, Kees},
  journal={Mathematics of Computation},
  volume={83},
  number={289},
  pages={2239--2262},
  year={2014}
}

@article{ryan2015one,
  title={One-sided position-dependent {S}moothness-{I}increasing {A}ccuracy-{C}onserving {(SIAC)} filtering over uniform and non-uniform meshes},
  author={Ryan, Jennifer K and Li, Xiaozhou and Kirby, Robert M and Vuik, Kees},
  journal={Journal of Scientific Computing},
  volume={64},
  number={3},
  pages={773--817},
  year={2015},
  publisher={Springer}
}

@article{suarez2017regularization,
  title={Regularization of singularities in the weighted summation of {D}irac-delta functions for the spectral solution of hyperbolic conservation laws},
  author={Suarez, Jean-Piero and Jacobs, Gustaaf B},
  journal={Journal of Scientific Computing},
  volume={72},
  number={3},
  pages={1080--1092},
  year={2017},
  publisher={Springer}
}

@article{wissink2018shock,
  title={Shock regularization with {S}moothness-{I}ncreasing {A}ccuracy-{C}onserving {D}irac-delta polynomial kernels},
  author={Wissink, BW and Jacobs, Gustaaf B and Ryan, Jennifer K and Don, Wai-Sun and Van Der Weide, ETA},
  journal={Journal of Scientific Computing},
  volume={77},
  number={1},
  pages={579--596},
  year={2018},
  publisher={Springer}
}

@misc{tanaka2019investigation,
    author = {Reiko Tanaka},
    title = {Investigation of Filters used in Large Eddy Simulation},
    note = {Bachelor thesis, University of Dusseldorf},
    year = {2019}
}

@article{bohm2019multi,
  title={Multi-element {SIAC} filter for shock capturing applied to high-order discontinuous {G}alerkin spectral element methods},
  author={Bohm, Marvin and Schermeng, Sven and Winters, Andrew R and Gassner, Gregor J and Jacobs, Gustaaf B},
  journal={Journal of Scientific Computing},
  volume={81},
  number={2},
  pages={820--844},
  year={2019},
  publisher={Springer}
}

@article{li2019smoothness,
  title={{S}moothness-{I}increasing {A}ccuracy-{C}onserving {(SIAC)} filtering for discontinuous {G}alerkin solutions over nonuniform meshes: superconvergence and optimal accuracy},
  author={Li, Xiaozhou and Ryan, Jennifer K and Kirby, Robert M and Vuik, Kees},
  journal={Journal of Scientific Computing},
  volume={81},
  number={3},
  pages={1150--1180},
  year={2019},
  publisher={Springer}
}

@article{docampo2020enhancing,
  title={Enhancing accuracy with a convolution filter: What works and why!},
  author={Docampo-S{\'a}nchez, J and Jacobs, GB and Li, X and Ryan, Jennifer K},
  journal={Computers \& Fluids},
  volume={213},
  pages={104727},
  year={2020},
  publisher={Elsevier}
}

@article{thomee1984galerkin,
  title={Galerkin finite element methods for parabolic problems},
  author={Thom{\'e}e, Vidar},
  journal={Lecture notes in mathematics},
  volume={1054},
  year={1984},
  publisher={Springer}
}

@book{gaul2014modified,
  title={Modified {B}ayesian {K}riging for noisy response problems and {B}ayesian confidence-based reliability-based design optimization},
  author={Gaul, Nicholas John},
  year={2014},
  publisher={The University of Iowa}
}

@article{sen2015evaluation,
  title={Evaluation of convergence behavior of metamodeling techniques for bridging scales in multi-scale multimaterial simulation},
  author={Sen, Oishik and Davis, Sean and Jacobs, Gustaaf and Udaykumar, HS},
  journal={Journal of Computational Physics},
  volume={294},
  pages={585--604},
  year={2015},
  publisher={Elsevier}
}

@inproceedings{gaul2015modified,
  title={Modified {B}ayesian {K}riging for noisy response problems for reliability analysis},
  author={Gaul, Nicholas J and Cowles, Mary Kathryn and Cho, Hyunkyoo and Choi, KK and Lamb, David},
  booktitle={International Design Engineering Technical Conferences and Computers and Information in Engineering Conference},
  volume={57083},
  pages={V02BT03A060},
  year={2015},
  organization={American Society of Mechanical Engineers}
}

@book{tikhonov2013numerical,
  title={Numerical Methods for the Solution of Ill-Posed Problems},
  author={Tikhonov, Andrei Nikolaevich and Goncharsky, AV and Stepanov, VV and Yagola, Anatoly G}, 
  series={Mathematics and Its Applications},
  volume={328},
  year={2013},
  publisher={Springer Science \& Business Media}
}

@book{vogel2002computational,
  title={Computational Methods for Inverse Problems},
  author={Vogel, Curtis R},
  year={2002},
  publisher={SIAM}
}

@book{hansen2006deblurring,
  title={Deblurring Images: Matrices, Spectra, and Filtering},
  author={Hansen, Per Christian and Nagy, James G and O'leary, Dianne P},
  year={2006},
  publisher={SIAM}
}

@book{kaipio2006statistical,
  title={Statistical and Computational Inverse Problems},
  author={Kaipio, Jari and Somersalo, Erkki},
  volume={160},
  year={2006},
  publisher={Springer Science \& Business Media}
}

@book{hansen2010discrete,
  title={Discrete Inverse Problems: Insight and Algorithms},
  author={Hansen, Per Christian},
  year={2010},
  publisher={SIAM}
}

@article{stuart2010inverse,
  title={Inverse problems: A {B}ayesian perspective},
  author={Stuart, Andrew M},
  journal={Acta Numerica},
  volume={19},
  pages={451--559},
  year={2010},
  publisher={Cambridge University Press}
}

@book{hansen2021computed,
  title={Computed Tomography: Algorithms, Insight, and Just Enough Theory},
  author={Hansen, Per Christian and J{\o}rgensen, Jakob and Lionheart, William RB},
  year={2021},
  publisher={SIAM}
}

@book{calvetti2023bayesian,
  title={Bayesian Scientific Computing},
  author={Calvetti, Daniela and Somersalo, Erkki},
  volume={215},
  year={2023},
  publisher={Springer Nature}
}

@article{tipping2001sparse,
  title={Sparse {B}ayesian learning and the relevance vector machine},
  author={Tipping, Michael E},
  journal={Journal of machine learning research},
  volume={1},
  number={Jun},
  pages={211--244},
  year={2001}
}

@article{zhang2011clarify,
  title={Clarify some issues on the sparse {B}ayesian learning for sparse signal recovery},
  author={Zhang, Zhilin and Rao, Bhaskar D},
  FJournal = {University of California, San Diego, Technical Report},
  Journal = {UCSD, Tech. Rep.},
  year={2011},
  publisher={Citeseer}
}

@article{bardsley2012mcmc,
  title={{MCMC}-based image reconstruction with uncertainty quantification},
  author={Bardsley, Johnathan M},
  journal={SIAM Journal on Scientific Computing},
  volume={34},
  number={3},
  pages={A1316--A1332},
  year={2012},
  publisher={SIAM}
}

@article{calvetti2018bayes,
  title={Bayes meets {K}rylov: Statistically inspired preconditioners for {CGLS}},
  author={Calvetti, Daniela and Pitolli, Francesca and Somersalo, Erkki and Vantaggi, Barbara},
  journal={SIAM Review},
  volume={60},
  number={2},
  pages={429--461},
  year={2018},
  publisher={SIAM}
}

@article{calvetti2019hierachical,
  title={Hierachical {B}ayesian models and sparsity: $\ell_2$-magic},
  author={Calvetti, Daniela and Somersalo, Erkki and Strang, A},
  journal={Inverse Problems},
  volume={35},
  number={3},
  pages={035003},
  year={2019},
  publisher={IOP Publishing}
}

@article{calvetti2020sparse,
  title={Sparse reconstructions from few noisy data: analysis of hierarchical {B}ayesian models with generalized gamma hyperpriors},
  author={Calvetti, Daniela and Pragliola, Monica and Somersalo, Erkki and Strang, Alexander},
  journal={Inverse Problems},
  volume={36},
  number={2},
  pages={025010},
  year={2020},
  publisher={IOP Publishing}
}

@article{calvetti2020sparsity,
  title={Sparsity promoting hybrid solvers for hierarchical {B}ayesian inverse problems},
  author={Calvetti, Daniela and Pragliola, Monica and Somersalo, Erkki},
  journal={SIAM Journal on Scientific Computing},
  volume={42},
  number={6},
  pages={A3761--A3784},
  year={2020},
  publisher={SIAM}
}

@article{glaubitz2023generalized,
  title={Generalized sparse {B}ayesian learning and application to image reconstruction},
  author={Glaubitz, Jan and Gelb, Anne and Song, Guohui},
  journal={SIAM/ASA Journal on Uncertainty Quantification},
  volume={11},
  number={1},
  pages={262--284},
  year={2023},
  publisher={SIAM}
}

@article{xiao2023sequential,
  title={Sequential image recovery using joint hierarchical {B}ayesian learning},
  author={Xiao, Yao and Glaubitz, Jan},
  journal={Journal of Scientific Computing},
  volume={96},
  number={1},
  pages={4},
  year={2023},
  publisher={Springer}
}

@article{LiGelb2025,
  title={A {B}ayesian framework for spectral reprojection},
  author={Li, Tongtong and Gelb, Anne},
  journal={Journal of Scientific Computing},
  volume={102},
  number={3},
  pages={78},
  year={2025},
  publisher={Springer}
}

@article{glaubitz2023leveraging,
  title={Leveraging joint sparsity in hierarchical {B}ayesian learning},
  author={Glaubitz, Jan and Gelb, Anne},
  journal={SIAM/ASA Journal on Uncertainty Quantification},
  volume={12},
  number={2},
  pages={442--472},
  year={2024},
  publisher={SIAM}
}

@article{lindbloom2024generalized,
  title={Efficient sparsity-promoting MAP estimation for {B}ayesian linear inverse problems},
  author={Lindbloom, Jonathan and Glaubitz, Jan and Gelb, Anne},
  journal={Inverse Problems},
  volume={41},
  number={2},
  pages={025001},
  year={2025},
  publisher={IOP Publishing}
}

@article{lindbloom2025priorconditioned,
  title={Priorconditioned sparsity-promoting projection methods for deterministic and {B}ayesian linear inverse problems},
  author={Lindbloom, Jonathan and Pasha, Mirjeta and Glaubitz, Jan and Marzouk, Youssef},
  journal={arXiv preprint arXiv:2505.01827},
  year={2025}
}

@article{glaubitz2025efficient,
  title={Efficient sampling for sparse {B}ayesian learning using hierarchical prior normalization},
  author={Glaubitz, Jan and Marzouk, Youssef},
  journal={arXiv preprint arXiv:2505.23753},
  year={2025}
}

@book{nocedal2006numerical,
  title={Numerical Optimization},
  author={Nocedal, Jorge and Wright, Stephen J},
  year={2006},
  publisher={Springer}
}

@article{wright2015coordinate,
  title={Coordinate descent algorithms},
  author={Wright, Stephen J},
  journal={Mathematical Programming},
  volume={151},
  number={1},
  pages={3--34},
  year={2015},
  publisher={Springer}
}

@book{beck2017first,
  title={First-Order Methods in Optimization},
  author={Beck, Amir},
  year={2017},
  publisher={SIAM}
}

@article{gelman1992inference,
  title={Inference from iterative simulation using multiple sequences},
  author={Gelman, Andrew and Rubin, Donald B},
  journal={Statistical Science},
  volume={7},
  number={4},
  pages={457--472},
  year={1992},
  publisher={Institute of Mathematical Statistics}
}

@article{brooks1998general,
  title={General methods for monitoring convergence of iterative simulations},
  author={Brooks, Stephen P and Gelman, Andrew},
  journal={Journal of Computational and Graphical Statistics},
  volume={7},
  number={4},
  pages={434--455},
  year={1998},
  publisher={Taylor \& Francis}
}

@article{haario2001adaptive,
  title={An adaptive {M}etropolis algorithm},
  author={Haario, Heikki and Saksman, Eero and Tamminen, Johanna},
  journal={Bernoulli},
  pages={223--242},
  year={2001},
  publisher={JSTOR}
}

@book{gelman2003bayesian,
  title={Bayesian Data Analysis},
  author={Gelman, Andrew and Carlin, John B and Stern, Hal S and Rubin, Donald B},
  year={2003},
  publisher={Chapman and Hall/CRC}
}

@article{wolff2004monte,
  title={Monte {C}arlo errors with less errors},
  author={Wolff, Ulli and Alpha Collaboration},
  journal={Computer Physics Communications},
  volume={156},
  number={2},
  pages={143--153},
  year={2004},
  publisher={Elsevier}
}

@article{andrieu2008tutorial,
  title={A tutorial on adaptive {MCMC}},
  author={Andrieu, Christophe and Thoms, Johannes},
  journal={Statistics and Computing},
  volume={18},
  pages={343--373},
  year={2008},
  publisher={Springer}
}

@article{atchade2010limit,
  title={Limit theorems for some adaptive {MCMC} algorithms with subgeometric kernels},
  author={Atchad{\'e}, Yves and Fort, Gersende},
  journal={Bernoulli},
  volume={16},
  number={1},
  year={2010},
  publisher={Bernoulli Society for Mathematical Statistics and Probability}
}

@article{vehtari2021rank,
  title={Rank-normalization, folding, and localization: An improved $\hat{R}$ for assessing convergence of {MCMC} (with discussion)},
  author={Vehtari, Aki and Gelman, Andrew and Simpson, Daniel and Carpenter, Bob and B{\"u}rkner, Paul-Christian},
  journal={Bayesian Analysis},
  volume={16},
  number={2},
  pages={667--718},
  year={2021},
  publisher={International Society for Bayesian Analysis}
}

@book{saad2003iterative,
  title={Iterative Methods for Sparse Linear Systems},
  author={Saad, Yousef},
  year={2003},
  publisher={SIAM}
}

@article{chung2024computational,
  title={Computational methods for large-scale inverse problems: a survey on hybrid projection methods},
  author={Chung, Julianne and Gazzola, Silvia},
  journal={Siam Review},
  volume={66},
  number={2},
  pages={205--284},
  year={2024},
  publisher={SIAM}
}

@book{hesthaven2008nodal,
  title={Nodal Discontinuous Galerkin Methods: Algorithms, Analysis, and Applications},
  author={Hesthaven, Jan S and Warburton, Tim},
  year={2008},
  publisher={Springer}
}

@article{bezanson2017julia,
  title={Julia: A fresh approach to numerical computing},
  author={Bezanson, Jeff and Edelman, Alan and Karpinski, Stefan and Shah, Viral B},
  journal={SIAM Review},
  volume={59},
  number={1},
  pages={65--98},
  year={2017},
  publisher={SIAM}
}

@article{verner2010numerically,
  title={Numerically optimal {R}unge--{K}utta pairs with interpolants},
  author={Verner, James H},
  journal={Numerical Algorithms},
  volume={53},
  number={2},
  pages={383--396},
  year={2010},
  publisher={Springer}
}

@misc{AdaptiveMCMC,
  title = {{AdaptiveMCMC.jl}: Adaptive {M}arkov chain {M}onte {C}arlo in {J}ulia},
  author = {Vihola, Matti},
  year = {2024}, 
  month={01}, 
  version={0.1.4},
  howpublished = {\url{https://github.com/mvihola/AdaptiveMCMC.jl}}
}
